\documentclass{amsart}
\usepackage{a4wide, graphicx, amsmath, amsfonts, amsthm, amssymb, latexsym, tikz, color}   

\usepackage[colorlinks]{hyperref}
\hypersetup{linkcolor=blue, urlcolor=blue, citecolor=red}

\newtheorem{main}{Theorem}

\newtheorem{thm}{Theorem}[section]
\newtheorem{cor}[thm]{Corollary}
\newtheorem{lem}[thm]{Lemma}
\newtheorem{prop}[thm]{Proposition}

\newtheorem{question}[thm]{Question}

\newtheorem{defn}[thm]{Definition}
\newenvironment{de}{\begin{defn} \rm}{\end{defn}}

\newcommand{\proofref}[1]{\noindent {\emph{Proof of Theorem} \ref{#1}.\ }}

\newcommand{\set}[2]{\{#1:#2\}}
\newcommand{\genset}[1]{\langle#1\rangle}
\newcommand{\fin}{\mathfrak{F}}
\newcommand{\sym}{\operatorname{Sym}}
\newcommand{\nat}{\mathbb{N}}
\newcommand{\N}{\mathbb{N}}

\newcommand{\astab}{\operatorname{AStab}}
\newcommand{\stab}{\operatorname{Stab}}

\newcommand{\supp}{\operatorname{supp}}

\newcommand{\mutt}[2]{\mathcal{C}(#1,  #2)}

\renewcommand{\P}{\mathcal{P}}
\newcommand{\F}{\mathcal{F}}

\newcommand{\n}{\{0,1,\ldots, n-1\}}

\usetikzlibrary{decorations.markings}
\tikzset{middlearrow/.style={
        decoration={markings,
            mark= at position 0.5 with {\arrow{#1}} ,
        },
        postaction={decorate}
    }
}
\tikzset{->-/.style={decoration={
  markings,
  mark=at position #1 with {\arrow{>}}},postaction={decorate}}}

\newcommand{\chip}[2]{
 \draw[line width=0.6mm](#1-.15,#2)--(#1+.15,#2);
}

\newcommand{\chipcolor}[3]{
 \draw[line width=0.6mm,#3](#1-.15,#2)--(#1+.15,#2);
}

\newcommand{\smallchip}[2]{
 \draw[thick](#1-.1,#2)--(#1+.1,#2);
}

\newcommand{\verticalline}[1]{
 \draw[line width=0.6mm](#1,0)--(#1,6);
 \chip{#1}{0}
 \chip{#1}{6}
}

\newcommand{\verticallinearrowleft}[4]{
 \draw[<->](#1-.3,#2)--(#1-.3,#3);
 \draw(#1-.3,{0.5*(#2+#3)}) node [left] {{\small #4}};
}

\newcommand{\verticallinearrowlefthigh}[4]{
 \draw[<->](#1-.3,#2)--(#1-.3,#3);
 \draw(#1-.3,{(0.25*#2)+(0.75*#3)}) node [left] {{\small #4}};
}

\newcommand{\verticallinearrowleftlow}[4]{
 \draw[<->](#1-.3,#2)--(#1-.3,#3);
 \draw(#1-.3,{(0.75*#2)+(0.25*#3)}) node [left] {{\small #4}};
}

\newcommand{\verticallinearrowright}[4]{
 \draw[<->](#1+.3,#2)--(#1+.3,#3);
 \draw(#1+.3,{0.5*(#2+#3)}) node [right] {{\small #4}};
}

\newcommand{\verticallinearrowrighthigh}[4]{
 \draw[<->](#1+.3,#2)--(#1+.3,#3);
 \draw(#1+.3,{(0.25*#2)+(0.75*#3)}) node [right] {{\small #4}};
}

\newcommand{\verticallinearrowrightlow}[4]{
 \draw[<->](#1+.3,#2)--(#1+.3,#3);
 \draw(#1+.3,{(0.75*#2)+(0.25*#3)}) node [right] {{\small #4}};
}

\newcommand{\functionline}[3]{
 \draw[->-=0.5] (#1,#2) to [out=0,in=180] (#1+3,#3); 
}

\newcommand{\functionlinedash}[3]{
 \draw[->-=0.5,dashed] (#1,#2) to [out=0,in=180] (#1+3,#3); 
}

\newcommand{\functionlineleft}[3]{
 \draw[->-=0.15] (#1,#2) to [out=0,in=180] (#1+3,#3); 
}

\newcommand{\functionlineright}[3]{
 \draw[->-=0.85] (#1,#2) to [out=0,in=180] (#1+3,#3); 
}

\newcommand{\functionlabel}[2]{
 \draw(#1+1.5,0)node[below]{$\phantom{fg}#2\phantom{fg}$};
 \draw(#1+1.5,6)node[above]{$\phantom{fg}#2\phantom{fg}$};
}

\newcommand{\pt}[2]{
 \fill(#1,#2)circle(0.08);
}

\newcommand{\colorinterval}[4]{
 \draw[line width=0.6mm,#4](#1,#2)--(#1,#3);
 \chipcolor{#1}{#2}{#4};
 \chipcolor{#1}{#3}{#4};
}

\newcommand{\functionlinealeft}[3]{
 \draw[->-=0.3] (#1,#2) to [out=0,in=180] (#1+3,#3); 
}

\newcommand{\functionlinearight}[3]{
 \draw[->-=0.7] (#1,#2) to [out=0,in=180] (#1+3,#3); 
}

\newcommand{\functionlinealeftdash}[3]{
 \draw[->-=0.3,dashed] (#1,#2) to [out=0,in=180] (#1+3,#3); 
}

\newcommand{\functionlinearightdash}[3]{
 \draw[->-=0.7,dashed] (#1,#2) to [out=0,in=180] (#1+3,#3); 
}

\begin{document}

\title[Maximal subsemigroups of $\Omega^\Omega$]{Maximal subsemigroups of the semigroup of all mappings on an 
infinite set}

\author{J. East, J. D. Mitchell, and Y. P\'eresse}
\begin{abstract}
In this paper we classify the maximal subsemigroups of the \emph{full transformation semigroup} $\Omega^\Omega$,
 which consists of all mappings on the 
infinite set $\Omega$, containing certain subgroups of the symmetric group $\sym(\Omega)$ on $\Omega$. 
In 1965 Gavrilov showed that there are five maximal subsemigroups of $\Omega^\Omega$ containing $\sym(\Omega)$
when $\Omega$ is countable and in  2005 Pinsker extended Gavrilov's result to sets of arbitrary cardinality.

We classify the maximal subsemigroups of $\Omega^\Omega$ on a set $\Omega$ of arbitrary infinite cardinality containing one 
of the following subgroups of  $\sym(\Omega)$: 
the pointwise stabiliser  of a non-empty finite subset of $\Omega$, 
the stabiliser of an ultrafilter on $\Omega$, or the stabiliser 
of a partition  of $\Omega$ into finitely many subsets of equal cardinality.  
If $G$ is any of these subgroups, then we deduce a characterisation of the 
mappings $f,g\in \Omega^\Omega$ such that the 
semigroup generated by $G\cup \{f,g\}$ equals $\Omega^\Omega$. 
\end{abstract}

\maketitle

%%%%%%%%%%%%%%%%%%%%%%%%%
%%%%%%%%%%%%%%%%%%%%%%%%%

\section{Introduction} 
A subgroup $H$ of a group $G$ is  a \emph{maximal subgroup} if $H\not =G$ and the subgroup  generated by $H$ and 
$g$ equals $G$ for all $g\in G\setminus H$. The definition of a \emph{maximal subsemigroup} of a semigroup is 
analogous: a subsemigroup $T$ of a semigroup (or group) $S$  is  a \emph{maximal subsemigroup} if $T\not =S$ and the 
subsemigroup $\genset{T, s}$  generated by $T$ and $s$ equals $S$ for all $s\in S\setminus T$. 

Let $\Omega$ denote an arbitrary (finite or infinite) set, let $\Omega^\Omega$ denote the 
semigroup of mappings from $\Omega$ to itself, and let $\sym(\Omega)$ denote the symmetric group on $\Omega$. 
In this paper we are interested in those maximal subsemigroups of $\Omega^\Omega$ that contain 
certain subgroups of 
$\sym(\Omega)$. The maximal subgroups of finite symmetric groups, having been investigated by O'Nan and Scott, see 
\cite{Scott1980aa},
 Aschbacher and Scott \cite{Aschbacher1985aa}, and Liebeck, Praeger and Saxl \cite{Liebeck1988aa},  
 are, in some sense, known. When $\Omega$ is finite, it is  easy to see   that a maximal 
 subsemigroup of 
 $\Omega^\Omega$ is either the union of a maximal subgroup of 
 the symmetric group and $\Omega^\Omega\setminus \sym(\Omega)$; or it is the union of $\sym(\Omega)$ and the 
 mappings with at most $|\Omega|-2$ points in their images.  In general,  
 the  maximal 
 subsemigroups of an arbitrary finite semigroup are  determined, roughly speaking, by their maximal 
 subgroups; see Graham, Graham, and Rhodes  \cite{Graham1968aa}. 
  
Maximal subgroups of $\sym(\Omega)$  have also been extensively studied when $\Omega$ is infinite;
see \cite{Ball1966aa, Ball1968aa, Baumgartner1993aa, Biryukov2000aa, Brazil1994aa, 
Covington1996aa,Macpherson1993aa, Macpherson1990aa,  Mishkin1995aa,Richman1967aa} and the references 
therein. It seems 
extremely unlikely that a complete description, in any sense, of maximal subgroups of $\sym(\Omega)$ exists for infinite 
$\Omega$. 
Maximal subsemigroups of $\Omega^\Omega$ when $\Omega$ is infinite have been considered to a lesser degree. 
The maximal subsemigroups of $\Omega^\Omega$ 
containing the symmetric group were classified by Gavrilov in \cite{Gavrilov1965aa} for countable $\Omega$ and 
Pinsker \cite[Theorem 1.4]{Pinsker2005aa} for arbitrary infinite $\Omega$; these are the only results regarding 
 maximal subsemigroups of  $\Omega^\Omega$, when $\Omega$ is infinite, of which we are aware. We state and prove Gavrilov and 
Pinsker's theorem  (Theorem \ref{dicks}) since elements of the proof are required  later on,  for 
the sake of completeness, and for the convenience of the reader.
Maximal subsemigroups of other infinite semigroups of mappings have been considered. For example, Levi and Wood 
\cite{Levi1984aa} and Hotzel \cite{Hotzel1995aa} considered maximal subsemigroups of Baer-Levi semigroups, and 
Shneperman \cite{Shneperman1993aa} 
considered the maximal subsemigroups of the endomorphism monoid of a finite dimension complex vector space that are
maximal with respect to being compact.  

The subsemigroups of $\Omega^\Omega$ form an algebraic lattice with $2^{|\Omega|}$ compact elements under 
inclusion. 
The study of maximal subsemigroups of $\Omega^\Omega$ belongs to the wider study of this lattice.  
Pinsker and Shelah \cite{Pinsker2011aa} prove that every alegbraic lattice with at most $2^{|\Omega|}$ compact 
elements 
can be embedded into the subsemigroup lattice of $\Omega^\Omega$. There are $2^{2^{|\Omega|}}$ distinct 
subsemigroups of $\Omega^\Omega$.   There are even $2^{2^{\kappa}}$ subsemigroups between $\sym(\Omega)$ and 
any maximal subsemigroup of $\Omega^\Omega$
that contains $\sym(\Omega)$ where $|\Omega|=\aleph_{\alpha}$ and $\kappa=\max\{\alpha, \aleph_0\}$;  for 
further details 
see Pinsker \cite{Pinsker2005ab}.  We show, as a consequence of Theorem \ref{ultra}, that there are also 
$2^{2^{|\Omega|}}$ maximal subsemigroups of $\Omega^\Omega$. 
The  maximal subsemigroups of the maximal subsemigroups described by Gavrilov \cite{Gavrilov1965aa} are classified 
in \cite{Dougall2012aa}; perhaps 
surprisingly there are only countably many such semigroups. In further contrast to Pinsker and Shelah's result 
\cite{Pinsker2011aa}, 
there are only 38 subsemigroups 
between the intersection $S_1\cap S_2\cap S_3(\aleph_0)\cap S_4(\aleph_0)\cap S_5$ of the maximal subsemigroups described by Gavrilov \cite{Gavrilov1965aa} and $\Omega^\Omega$; see Mitchell and Jonu\v sas \cite{Jonusas2013aa}. 

Another natural question to ask about the subsemigroup lattice of $\Omega^\Omega$ is whether or not every 
subsemigroup is contained in a maximal one.  In \cite{Baumgartner1993aa} it is shown that under certain set theoretic 
assumptions 
there exists a subgroup of $\sym(\Omega)$ that is not contained in a maximal subgroup; it seems likely that the 
analogous 
result holds for $\Omega^\Omega$.  There are several results in the literature concerning sufficient conditions for a 
subgroup of $\sym(\Omega)$ to lie in a maximal subgroup; see \cite{Macpherson1990aa} and \cite{Macpherson1990ab}. 
In Section \ref{containment} we explore 
the analogous problem for subsemigroups of $\Omega^\Omega$. 

In this paper we classify the maximal 
subsemigroups of $\Omega^\Omega$, where $\Omega$ is any infinite set, containing certain subgroups of 
$\sym(\Omega)$, which we define in the next section.   In particular, we classify the maximal subsemigroups of 
$\Omega^\Omega$ 
containing one of the following groups: the symmetric group $\sym(\Omega)$ (Theorem \ref{dicks}), the pointwise 
stabiliser of a non-empty finite subset of $\Omega$ (Theorem \ref{stab_finite}),  the stabiliser of an ultrafilter on $\Omega$ 
(Theorem \ref{ultra}), or the  stabiliser of a finite partition of $\Omega$ (Theorem \ref{almost_stab}). 
For each of these subgroups,  we obtain a characterisation of those  pairs of 
elements  that, together with the subgroup, generate 
$\Omega^\Omega$; see Corollaries \ref{hhr_new}, \ref{stab_gen_pairs}, \ref{ultra_pairs} 
and \ref{almost_stab_pairs}.  Such a classification in the case that $G=\sym(\Omega)$ and $|\Omega|$ is a regular 
cardinal  was originally given in  \cite[Theorem 3.3]{Howie1998aa}.
As previously mentioned the classification of maximal subsemigroup of $\Omega^\Omega$ containing $\sym(\Omega)$ is 
originally due to Gavrilov  \cite{Gavrilov1965aa}  and 
Pinsker \cite{Pinsker2005aa}.

The paper is organised as follows: in Section \ref{statements} 
we state the main theorems of the paper.  In Section \ref{containment}, we give several sufficient conditions 
for a subsemigroup of $\Omega^\Omega$ to be contained in a maximal subsemigroup, and also give a new proof of the 
result of Macpherson and Praeger \cite{Macpherson1990ab} which states that every subgroup of $\sym(\Omega)$ that is 
not highly transitive  is contained in a maximal subgroup.  
In Section \ref{corollaries}, we state and prove  Corollaries \ref{hhr_new}, \ref{stab_gen_pairs}, \ref{ultra_pairs} 
and \ref{almost_stab_pairs}.
In Section \ref{mundane}, we give several 
technical results which underpin the proofs of the main  results in the paper. In Sections \ref{proof_dicks}, \ref{proof_stab_finite}, \ref{proof_ultra}, and \ref{proof_almost_stab} we
 give the respective proofs of the four main 
theorems from Section \ref{statements}.  In Section \ref{subsemigroups}, 
we  show that the setwise stabiliser of a non-empty finite set, 
the almost stabiliser 
of a finite partition, and the stabiliser of an ultrafilter are  maximal subsemigroups (and not just maximal subgroups as is 
already well-known) of the symmetric group. 

We end this section by asking the three most interesting questions, in our eyes at least, arising from our consideration of 
maximal subsemigroups of $\Omega^\Omega$. 

\begin{question}
 Let $G$ be a maximal subsemigroup of $\sym(\Omega)$. Then does there exist a maximal subsemigroup $M$ of 
 $\Omega^\Omega$ such that $M\cap \sym(\Omega)=G$? 
 \end{question}

The intersection of every known example of a maximal subsemigroup of $\Omega^\Omega$ with $\sym(\Omega)$ is 
either a maximal subsemigroup of $\sym(\Omega)$ or $\sym(\Omega)$ itself, which prompts the following question. 

\begin{question}\label{inter}
Does there exist a maximal subsemigroup of $\Omega^\Omega$ that does not contain a
maximal subsemigroup of $\sym(\Omega)$?
\end{question}

We suspect that the answer to Question \ref{inter} is yes. A step in the other direction would, perhaps, 
be a positive answer to the 
following question.

\begin{question}
 Does every maximal subsemigroup of $\Omega^\Omega$ 
have  non-trivial intersection with $\sym(\Omega)$?
\end{question}

%%%%%%%%%%%%%%%%%%%%%%%%%
%%%%%%%%%%%%%%%%%%%%%%%%%

\section{Statements of the Main Theorems}\label{statements}
Throughout the paper we write functions to the right of their argument and compose from left to right. 
If $\alpha\in \Omega$, $f\in\Omega^\Omega$ and $\Sigma\subseteq \Omega$, then  
$\alpha f^{-1}=\set{\beta\in \Omega}{\beta f=\alpha}$, $\Sigma f=\set{\alpha f}{\alpha\in \Sigma}$, and $f|_{\Sigma}$ 
denotes the \emph{restriction} of $f$ 
to $\Sigma$. 
We denote $\set{f\in \Omega^\Omega}{|\Omega f|<|\Omega|}$ by $\mathfrak{F}$.  Since $\mathfrak{F}$ is an ideal of
$\Omega^\Omega$, if $S$ is any subsemigroup of $\Omega^\Omega$, then so is $S\cup \fin$. Hence if $S$ is maximal, then 
either $\fin\subseteq S$ or $S\cup \fin=\Omega^\Omega$. In the latter case, $\Omega^\Omega\setminus \fin$ is a subset of $S$. But 
$\Omega^\Omega\setminus \fin$  is also a generating set for $\Omega^\Omega$ and so $S=\Omega^\Omega$,
which contradicts the assumption that $S$ is a maximal subsemigroup. Hence $\fin$ 
 is contained in every maximal subsemigroup of $\Omega^\Omega$. 

Let $\Sigma$ be any subset of $\Omega$ and let $f:\Sigma\to \Omega$ be arbitrary. 
If $\Gamma\subseteq \Sigma$ such that  $f|_{\Gamma}$ is injective and $\Gamma f=\Sigma f$, 
then we will refer to $\Gamma$ as a \emph{\hypertarget{transversal}{transversal}} of $f$. 
We require the following parameters of $f$ to  state our main theorems: 
\begin{eqnarray*}
d(f)& =& |\Omega\setminus\Sigma f|\\
c(f)& =&|\Sigma\setminus \Gamma|\text{, where }\Gamma\text{ is any transversal of }f\\
k(f, \mu)& =& |\set{\alpha\in\Omega}{|\alpha f^{-1}|\geq \mu}|\text{, where }\mu\leq |\Omega|.
\end{eqnarray*}
 The parameters $d(f),c(f)$, and $k(f, |\Omega|)$ were termed the \emph{defect}, \emph{collapse}, and \emph{infinite 
 contraction index}, respectively, of $f$ in~\cite{Howie1998aa}. 

As usual, we will think of a cardinal $\kappa$ as the set of all ordinals strictly less than $\kappa$.  Recall that a cardinal 
$\kappa$ is \emph{singular} if there exists a cardinal $\lambda<\kappa$ and a family of sets $\Sigma_\mu$ 
($\mu\in\lambda$) such that $|\Sigma_\mu|<\kappa$ for each $\mu<\lambda$, yet 
$\big|\bigcup_{\mu<\lambda}\Sigma_\mu\big|=\kappa$; otherwise, $\kappa$ is \emph{regular}. We denote the successor 
to any cardinal $\kappa$ by $\kappa^+$. 

A subset $\Sigma$ of an infinite set $\Gamma$ is a \emph{moiety} of $\Gamma$ if 
$|\Sigma|=|\Gamma\setminus \Sigma|=|\Gamma|$. 

\subsection{The symmetric group}

 %%%%%%%%%%%%%%%%%%%%%%%%%

\begin{main}[Gavrilov \cite{Gavrilov1965aa}, Pinsker \cite{Pinsker2005aa}]\label{dicks} 
% Note that we are not insulting Gavrilov or Pinsker, we simply named the semigroups D, I, C, K, S when we first
% proved this theorem, before we knew it already existed. 
Let $\Omega$ be any infinite set.  If $|\Omega|$ is a regular cardinal, then the maximal subsemigroups of $\Omega^\Omega$ containing $\sym(\Omega)$ are: 
\begin{align*}
S_1&=\set{f \in \Omega^\Omega}{c(f)= 0 \text{ or }d(f)>0};\\
S_2&=\set{f \in \Omega^\Omega}{c(f)>0 \text{ or }d(f)=0};\\
S_3(\mu)&=\set{f \in \Omega^\Omega}{c(f)< \mu \text{ or }d(f)\geq \mu};\\
S_4(\mu)&=\set{f \in \Omega^\Omega}{c(f)\geq\mu \text{ or }d(f)<\mu};\\
S_5&=\set{f \in \Omega^\Omega}{k(f, |\Omega|)<|\Omega|};\\ 
\intertext{where $\mu$ is any infinite cardinal not greater than $|\Omega|$.}
\intertext{If $|\Omega|$ is a singular cardinal, then the maximal subsemigroups of $\Omega^\Omega$ containing $\sym(\Omega)$ are $S_1$, $S_2$, $S_3(\mu)$, $S_4(\mu)$ where $\mu$ is any infinite cardinal not greater than $|\Omega|$,  and: 
}
S_5^{\prime}&=\set{f \in \Omega^\Omega}{(\exists \nu<|\Omega|)\,(k(f, \nu)<|\Omega|)}.
\end{align*}
\end{main}

%%%%%%%%%%%%%%%%%%%%%%%%%

The countable case of Theorem \ref{dicks} was first proved by Gavrilov \cite{Gavrilov1965aa}. The full version of Theorem~\ref{dicks} 
given above was first proved by  Pinsker \cite[Theorem 1.4]{Pinsker2005aa}.  We independently proved Theorem 
\ref{dicks} whilst unaware of the work of Gavrilov and Pinsker. We thank Martin Goldstern and Lutz Heindorf for bringing 
these references to our attention.
 A full proof of Theorem \ref{dicks} is included in Section \ref{proof_dicks} for the convenience of the reader and the sake 
 of completeness.

%%%%%%%%%%%%%%%%%%%%%%%%%

%%%%%%%%%%%%%%%%%%%%%%%%%

\subsection{The pointwise stabiliser of a finite set}
If $G$ is a group acting on a set $\Omega$ and $\Sigma$ is any subset of $\Omega$, then we denote the 
\emph{\hypertarget{pt_stab}{pointwise stabiliser}} of $\Sigma$ under $G$ by $G_{(\Sigma)}$ and the 
\emph{\hypertarget{set_stab}{setwise stabiliser}} of $\Sigma$ under $G$ by $G_{\{\Sigma\}}$.  
%For the sake of 
%convenience, we denote the pointwise stabiliser $\sym(\Omega)_{(\Omega\setminus \Sigma)}$ by $\sym(\Sigma)$. 
In \cite{Ball1966aa}, it is shown that if $\Sigma$ is a non-empty finite subset of $\Omega$, then  
$\sym(\Omega)_{\{\Sigma\}}$ is a maximal subgroup of $\sym(\Omega)$. 

%%%%%%%%%%%%%%%%%%%%%%%%%

\begin{main}  \label{stab_finite}
Let $\Omega$ be any infinite set and let $\Sigma$ be a non-empty finite subset of $\Omega$.
Then the maximal subsemigroups of $\Omega^\Omega$ containing the pointwise stabiliser 
$\sym(\Omega)_{(\Sigma)}$   but not 
$\sym(\Omega)$ are:
$$F_1(\Gamma, \mu)=\set{f\in \Omega^\Omega}{d(f)\geq \mu \text{ or } \Gamma \not \subseteq \Omega f \text{ or } 
(\Gamma f^{-1}\subseteq \Gamma \text{ and } c(f)<\mu)} \cup \fin;$$
$$F_2(\Gamma, \nu)=\set{f\in \Omega^\Omega}{c(f)\geq \nu \text{ or } |\Gamma f| < |\Gamma| \text{ or } 
(\Gamma f =  \Gamma \text{ and } d(f)<\nu)} \cup \fin$$
where $\Gamma$ is a non-empty subset of $\Sigma$ and $\mu$ and $\nu$ are  infinite cardinals with 
$\mu \leq |\Omega|^+$ and either:
$|\Gamma|= 1$ and $\nu= |\Omega|^{+}$; or $|\Gamma|\geq 2$ and $\nu\leq |\Omega|^{+}$. 
\end{main}

If $\mu\leq |\Omega|$  and $f\in \fin$, then
 $d(f)=|\Omega|=c(f)$, and so $``\cup\fin"$ could be omitted from the definition of 
$F_1(\Gamma, \mu)$ and $F_2(\Gamma, \mu)$ in these cases.

%%%%%%%%%%%%%%%%%%%%%%%%%
%If $\Gamma$ is any finite subset of $\Omega$ and $\nu$ is any infinite cardinal such that $\nu\leq |\Omega|$, then 
%$F_1(\Gamma, \nu)$ and $F_2(\Gamma, \nu)$ are semigroups. 

If $|\Gamma|=1$, then $F_2(\Gamma,\nu)$ is properly contained in 
$S_4(\nu)$ for all $\nu\leq |\Omega|$. In particular, $F_2(\Gamma,\nu)$ is not maximal in this case. 
When $\mu$ or $\nu$ equals $|\Omega|^+$,  we obtain the following simpler definitions of the semigroups in Theorem 
\ref{stab_finite}: 
\begin{align*}
F_1(\Gamma, |\Omega|^+)&=\set{f\in \Omega^\Omega}{\Gamma \not \subseteq \Omega f \text{ or }
 \Gamma f^{-1}\subseteq \Gamma}\cup \fin\\
F_2(\Gamma, |\Omega|^+)&=\set{f\in \Omega^\Omega}{|\Gamma f| < |\Gamma| \text{ or } \Gamma f =  \Gamma} \cup \fin.
\intertext{In particular, if $\Gamma=\{\gamma\}$, then }
F_1(\Gamma, |\Omega|^+)&=\set{f\in \Omega^\Omega}{\gamma \not \in \Omega f \text{ or } 
\gamma f^{-1}=\{\gamma\}}\cup\mathfrak{F}\\
F_2(\Gamma,  |\Omega|^{+})&=\set{f\in \Omega^\Omega}{\gamma f =  \gamma} \cup \fin.
\end{align*}

If $\Gamma$ is any finite subset of $\Omega$, then the intersection of $F_1(\Gamma, \mu)$ or $F_2(\Gamma, \mu)$ 
with $\sym(\Omega)$ is the setwise stabiliser 
$\sym(\Omega)_{\{\Gamma\}}$. Thus every maximal subsemigroup of $\Omega^\Omega$ containing the pointwise 
stabiliser of a finite subset $\Sigma$ of $\Omega$ also contains the setwise stabiliser of some subset $\Gamma$ of 
$\Sigma$. 
Since $\sym(\Omega)_{\{\Sigma\}}$ is a maximal subgroup of $\sym(\Omega)$, 
it follows that the maximal subsemigroups of $\Omega^\Omega$ containing 
$\sym(\Omega)_{\{\Sigma\}}$   but not 
$\sym(\Omega)$ are those listed in Theorem \ref{stab_finite} where $\Gamma=\Sigma$.

%%%%%%%%%%%%%%%%%%%%%%%%%
%%%%%%%%%%%%%%%%%%%%%%%%%

\subsection{The stabiliser of an ultrafilter}

A set of subsets $\mathcal{F}$ of $\Omega$ is called a \emph{filter} if:
\begin{itemize}
\item[\rm (i)] $\emptyset\not\in\F$;
\item[\rm (ii)] if $\Sigma\in \mathcal{F}$ and $\Sigma\subseteq \Gamma\subseteq \Omega$, then $\Gamma\in \mathcal{F}$;
\item[\rm (iii)] if $\Sigma, \Gamma\in \mathcal{F}$, then $\Sigma\cap \Gamma\in \mathcal{F}$.
\end{itemize} 
A filter is called an \emph{ultrafilter} if it is maximal with respect to containment among filters on $\Omega$. Equivalently, a filter $\mathcal{F}$ is an ultrafilter if, for every $\Sigma\subseteq\Omega$, either $\Sigma\in\mathcal{F}$ or $\Omega\setminus\Sigma\in\mathcal{F}$.  An ultrafilter $\F$ on $\Omega$ is \emph{principal} if there exists $\alpha\in\Omega$ such that $\F=\set{\Sigma\subseteq \Omega}{\alpha\in \Sigma}$. An ultrafilter $\F$ is \emph{uniform} if $|\Sigma|=|\Omega|$ for all $\Sigma\in \F$. 
The \emph{stabiliser} of a filter $\F$ in $\sym(\Omega)$ is defined to be 
$$\set{f\in \sym(\Omega)}{(\forall \Sigma\subseteq \Omega)(\Sigma\in \F \leftrightarrow \Sigma f\in \mathcal{F})}$$
and is denoted by $\sym(\Omega)_{\{\F\}}$. The stabiliser of an ultrafilter is the 
union of the pointwise stabilisers of the sets in the filter; i.e.
$$\sym(\Omega)_{\{\F\}}=\bigcup_{\Sigma\in \F} \sym(\Omega)_{(\Sigma)};$$
see \cite[Theorem 6.4]{Macpherson1990aa}.
It is shown in \cite[Theorem 6.4]{Macpherson1990aa} and \cite{Richman1967aa} that the stabiliser  
$\sym(\Omega)_{\{\F\}}$ of any ultrafilter is a maximal subgroup of the symmetric group. 

Let $\F$ be any filter on $\Omega$ and let $\mu$ be an infinite cardinal. Then we define the following 
subsemigroups of $\Omega^\Omega$:
\begin{align*}
U_1(\F, \mu)&=\set{f\in \Omega^\Omega}{(d(f)\geq\mu)\text{ or }(\Omega f\not\in \F) \text{ or }(c(f)< \mu \text{ and } 
(\forall \Sigma \not\in \F)(\Sigma f\not\in \F))}\cup \fin;\\
U_2(\F, \mu)&=\set{f\in \Omega^\Omega}{(c(f)\geq\mu)\text{ or }(\forall \Sigma \in \F)( c(f|_\Sigma)>0)\text{ or }(d(f)< \mu
\text{ and }(\forall \Sigma \in \F)(\Sigma f\in \F))}\cup \fin.
\end{align*}
If $\mu\leq|\Omega|$ and $f\in \fin$, then $d(f)=|\Omega|=c(f)$, and so $``\cup\fin"$ could be omitted from the definition of 
$U_1(\F, \mu)$ and $U_2(\F, \mu)$ in these cases.

If $\Gamma$ is any subset of $\Omega$, then the collection $\F$ of subsets of $\Omega$ containing $\Gamma$ is a filter. 
In this case, the stabiliser of $\F$ in $\sym(\Omega)$ and the setwise stabiliser of $\Gamma$ in $\sym(\Omega)$ coincide.  
In the following lemma, we show that $U_1(\F, \mu)$ and $U_2(\F, \mu)$ are the generalisations of the semigroups in Theorem 
\ref{stab_finite} to arbitrary filters. 

\begin{lem}
Let $\Gamma$ be a non-empty 
finite subset of $\Omega$ and let $\F$ be the filter consisting of  subsets of $\Omega$ containing $\Gamma$.
Then $F_1(\Gamma, \mu)=U_1(\F, \mu)$ and $F_2(\Gamma, \mu)=U_2(\F, \mu)$ for all infinite cardinals $\mu$.
\end{lem}
\proof
Suppose that $f\in \Omega^\Omega$. 
Clearly, $\Omega f\not\in \F$ if and only if $\Gamma\not\subseteq \Omega f$. Also $\Sigma f\not\in \F$ for all 
$\Sigma\not\in \F$ if and only if $\Gamma\not\subseteq \Sigma f$ for all $\Sigma \subseteq \Omega$ such that 
$\Gamma\not\subseteq \Sigma$ if and only if $\Gamma f^{-1}\subseteq \Gamma$. Therefore $F_1(\Gamma, \mu)=U_1(\F, \mu)$.

It is straightforward to show that $c(f|_\Sigma)>0$ for all $\Sigma\in \F$ if and only if no transversal of $f$ 
belongs to $\F$ if and only if $\Gamma$ is not a subset of any transversal of $f$ if and only if $|\Gamma f|<|\Gamma|$. 
Suppose that $|\Gamma f|=|\Gamma|$. Then $\Sigma f\in \F$ for all $\Sigma \in \F$ if and only if 
$\Gamma \subseteq \Sigma f$ for all $\Sigma \subseteq \Omega$
such that $\Gamma\subseteq \Sigma$ if and only if $\Gamma f\subseteq \Gamma$ if and only if $\Gamma f=\Gamma$.  
Thus $F_2(\Gamma, \nu)=U_2(\F, \nu)$, as required.\qed\vspace{\baselineskip}

The semigroups in Theorem 
\ref{stab_finite} contain  not only the pointwise stabiliser, but the setwise stabiliser of a finite set. It follows that  the maximal 
subsemigroups of $\Omega^\Omega$ containing the stabiliser of a filter generated by a finite set, in particular principal ultrafilters, 
have already been classified in Theorem \ref{stab_finite}. 
For the sake of convenience, we state the analogue of Theorem \ref{stab_finite} in terms of filters. 

%%%%%%%%%%%%%%%%%%%%%%%%%

\begin{cor}\label{filter_finite}
Let $\Omega$ be any infinite set, let $\Gamma$ be a non-empty finite subset of $\Omega$, and let
$\F$ be the filter consisting of subsets of $\Omega$ containing $\Gamma$. 
Then the maximal subsemigroups of $\Omega^\Omega$ containing 
$\sym(\Omega)_{\{\F\}}$ but not $\sym(\Omega)$ are:
$F_1(\Gamma, \mu)=U_1(\F, \mu)$
and 
$F_2(\Gamma, \nu)=U_2(\F, \nu)$
where $\mu$ and $\nu$ are  infinite cardinals with 
$\mu \leq |\Omega|^+$ and either:
$|\Gamma|= 1$ and $\nu= |\Omega|^{+}$; or $|\Gamma|\geq 2$ and $\nu\leq |\Omega|^{+}$.  
\end{cor}

%%%%%%%%%%%%%%%%%%%%%%%%%

If $|\Gamma|=1$, then $\F$ in Corollary \ref{filter_finite} is a principal ultrafilter. Replacing this principal ultrafilter   
by a non-principal ultrafilter yields the following theorem, which is similar to Corollary \ref{filter_finite}; 
the main difference being the possible values that the cardinals $\mu$ and $\nu$ can have. 

%%%%%%%%%%%%%%%%%%%%%%%%%

\begin{main}\label{ultra}
Let $\Omega$ be any infinite set, let $\F$ be a non-principal ultrafilter on $\Omega$, and let $\kappa(\geq \aleph_0)$ be 
the least 
cardinality of a subset of $\Omega$ in $\F$.  Then the maximal subsemigroups of $\Omega^\Omega$ containing 
$\sym(\Omega)_{\{\F\}}$ but not $\sym(\Omega)$ are $U_1(\F, \mu)$ and $U_2(\F, \mu)$ where 
 $\mu$ is an infinite cardinal such that $\kappa<\mu\leq |\Omega|^+$. 
\end{main}

Suppose that $\F$ is a non-principal ultrafilter. 
If $f\in \Omega^\Omega$ such that $\Omega f\not\in \F$, then 
$\Omega\setminus \Omega f\in \F$ and so $d(f)=|\Omega\setminus \Omega f|\geq \kappa.$ Hence 
 if $\mu\leq \kappa$, then
 $$U_1(\F, \mu)=\set{f\in \Omega^\Omega}{(d(f)\geq\mu)\text{ or }(c(f)< \mu \text{ and } 
(\forall \Sigma \not\in \F)(\Sigma f\not\in \F))}\subsetneq S_3(\mu),$$
and $U_1(\F, \mu)$ is not maximal in this case. 
If $f\in \Omega^\Omega$ is such that $c(f|_{\Sigma})>0$ for all $\Sigma\in \F$, then no transversal of $f$ 
belongs to $\F$. Hence the complement of any transversal of $f$ belongs to $\F$, and so $c(f)\geq \kappa$. 
In particular, if $\mu\leq \kappa$, then 
$$U_2(\F, \mu)=\set{f\in \Omega^\Omega}{(c(f)\geq\mu)\text{ or }(d(f)< \mu\text{ and }(\forall \Sigma \in \F)(\Sigma f\in \F)))}\subsetneq S_4(\mu)$$
and so $U_2(\F, \mu)$ is also not maximal in this case. 

If $\F$ in Theorem \ref{ultra} is an uniform ultrafilter, then $\kappa=|\Omega|$ and so there is only one possible value 
for 
$\mu$, namely $|\Omega|^+$, and the conditions on $U_1(\F, |\Omega|^+)$ and $U_2(\F, |\Omega|^+)$  
become much simpler:
\begin{align*}
U_1(\F, |\Omega|^+)&=\set{f\in \Omega^\Omega}{(\forall \Sigma \not\in \F)(\Sigma f\not\in \F)}\cup \fin\\
U_2(\F, |\Omega|^+)&=\set{f\in \Omega^\Omega}{(\forall \Sigma \in \F)( c(f|_\Sigma)>0)\text{ or }(\forall \Sigma \in \F)(\Sigma f\in \F)}\cup \fin.
\end{align*}

%%%%%%%%%%%%%%%%%%%%%%%%%

There are $2^{|\Omega|}$ elements in $\sym(\Omega)$ and by Posp\'i\u sil's Theorem \cite[Theorem 7.6]{Jech2003aa} 
there are $2^{2^{|\Omega|}}$  ultrafilters on $\Omega$. Hence  there are $2^{2^{|\Omega|}}$ non-conjugate maximal 
subsemigroups of $\Omega^\Omega$. 

%%%%%%%%%%%%%%%%%%%%%%%%%

While the semigroups in Corollary \ref{filter_finite} and Theorem \ref{ultra} 
have the 
same definitions in terms of their respective filters, neither result appears to be a corollary of the other. 
We were unable to formulate a more general theorem having Corollary \ref{filter_finite} and Theorem \ref{ultra}, let alone 
Theorem \ref{stab_finite}, as special cases.

%%%%%%%%%%%%%%%%%%%%%%%%%
\subsection{The stabiliser of a finite partition}\label{subsection_almost}

Let $n\geq 2$ and let $\mathcal{P}=\{\Sigma  _0, \Sigma_1, \ldots, \Sigma_{n-1}\}$ be a partition of $\Omega$  such that $|\Sigma_0|=\cdots =|\Sigma_{n-1}|=|\Omega|$. 
 We will refer to such a partition $\mathcal{P}$ as a \emph{finite partition} of $\Omega$. 
The \emph{stabiliser} of a finite partition $\mathcal{P}=\{\Sigma_0, \Sigma_1, \ldots, \Sigma_{n-1}\}$ is defined by
$$\stab(\P)=\set{f\in\sym(\Omega)}{(\forall\ i)(\exists\ j)(\Sigma_if= \Sigma_j)}$$
and the \emph{\hypertarget{almost_stab}{almost stabiliser}} of $\P$ is defined by 
$$\astab(\mathcal{P})=\set{f\in\sym(\Omega)}{(\forall\ i)(\exists\ j)(|\Sigma_if\setminus\Sigma_j|+|\Sigma_j\setminus\Sigma_if| <|\Omega|)}.$$
Of course, $\stab(\P)$ is a subgroup of $\astab(\P)$ and so $\stab(\P)$ is not a maximal subgroup of $\sym(\Omega)$. On 
the other hand, it was shown in \cite{Richman1967aa} (and \cite{Macpherson1990aa} independently)  that $\astab(\P)$ is 
a maximal subgroup of $\sym(\Omega)$. 

Let 
$f\in \Omega^\Omega$. Then define the binary relation $\rho_f$ on $\{0,1,\ldots, n-1\}$ by 
\begin{equation}\label{rho_f}
\rho_f=\set{(i,j)}{|\Sigma_if\cap \Sigma_j|=|\Omega|}.
\end{equation} 
If $\sigma$ is a binary relation on a set $\Omega$, then  $\sigma^{-1}=\set{(i,j)}{(j,i)\in\sigma}$ and $\sigma$ is 
\emph{total} if for all $\alpha\in \Omega$ there exists $\beta\in \Omega$ such that $(\alpha, \beta)\in \sigma$.  We will 
write $\sym(n)$ for the symmetric group on the set $n=\{0,1,\ldots,n-1\}$.

%%%%%%%%%%%%%%%%%%%%%%%%%

\begin{main} \label{almost_stab} 
Let $\Omega$ be any infinite set and let $\mathcal{P}=\{\Sigma_0, \Sigma_1, \ldots,\Sigma_{n-1}\}$, $n\geq 2$, be a finite partition of $\Omega$. Then the maximal subsemigroups of $\Omega^\Omega$ containing $\stab(\P)$ but not $\sym(\Omega)$ are:
\begin{align*}
A_1(\P)&=\set{f\in \Omega^\Omega}{\rho_f\in \sym(n)\text{ or }\rho_f\text{ is not total }};\\
A_2(\P)&=\set{f\in \Omega^\Omega}{\rho_f\in \sym(n)\text{ or }\rho_f^{-1}\text{ is not total }}.
\end{align*}
\end{main}

If $\P$ is any finite partition of $\Omega$, then the intersection of $A_1(\P)$ and $A_2(\P)$ with $\sym(\Omega)$ 
is the almost stabiliser $\astab(\P)$ of $\P$. Thus every maximal subsemigroup of $\Omega^\Omega$ containing the  
stabiliser of $\P$ also contains the almost stabiliser of $\P$.  
%%%%%%%%%%%%%%%%%%%%%%%%%
%%%%%%%%%%%%%%%%%%%%%%%%%

\section{Containment}\label{containment}

In this section we consider the question of when a subsemigroup of $\Omega^\Omega$ is contained in a maximal 
subsemigroup.  The analogous question has been considered for subgroups of the symmetric group; see, for example,
 \cite{Baumgartner1993aa, Macpherson1990aa, Macpherson1990ab}. 
The proposition below is of particular interest here. 
 In \cite{Baumgartner1993aa} it is shown that under certain set theoretic assumptions 
there exists a subgroup of $\sym(\Omega)$ that is not contained in a maximal subgroup. However, such examples are
difficult to find, and, roughly speaking,
if a subgroup of $\sym(\Omega)$ is large or small enough, then it is contained in a maximal subgroup.

It will be convenient to use the following notion: if $S$ is a semigroup and $T$ is subset of $S$, then the \emph{relative
rank} of $T$ in $S$ is the least cardinality of a subset $U$ of $S$ such that $\genset{T, U}=S$. 

 Part  (i) of the following proposition is a special case of Lemma 6.9 in Macpherson and Neumann \cite{Macpherson1990aa}, and parts
  (ii) and (iii)  are Theorems 1.5 and 1.6 in Macpherson and Praeger
 \cite{Macpherson1990ab}.
 
\begin{prop}\label{sym_contain}
Let $G$ be a subgroup of $\sym(\Omega)$ satisfying any of the following:
\begin{itemize}
\item[\rm (i)]  $G$ has finite relative rank in $\sym(\Omega)$;
\item[\rm (ii)] $|G|\leq |\Omega|$;
\item[\rm (iii)]  $|\Omega|$ is countable and there exists $t\in\nat$ such that $G$ has infinitely many orbits on $\Omega^t$.
\end{itemize}
Then $G$ is contained in a maximal subgroup of $\sym(\Omega)$. 
\end{prop}

\begin{prop}[Bergman-Shelah, Section 5 in \cite{Bergman2006ac}]\label{bergmanshelah}
Let $\Omega$ be countably infinite and let 
$G$ be a subgroup of $\sym(\Omega)$ such that  $G_{(\Sigma)}$  has an infinite orbit for all finite 
$\Sigma\subseteq \Omega$. Then $G$ has finite relative rank in $\sym(\Omega)$ and hence is contained in a maximal 
subgroup.
\end{prop}

We give an analogue of Proposition \ref{sym_contain}(i) and (ii) for subsemigroups of $\Omega^\Omega$. 

\begin{prop}\label{contain}
Let $S$ be a subsemigroup of $\Omega^\Omega$ satisfying either of the following:
\begin{itemize}
\item[\rm (i)]  $S$ has finite relative rank in $\Omega^\Omega$;
\item[\rm (ii)] $|S|\leq |\Omega|$.
\end{itemize}
Then $S$ is contained in a maximal subsemigroup of $\Omega^\Omega$. 
\end{prop}
\proof  
\noindent{\bf (i).} This is a straightforward consequence of Zorn's Lemma, analogous to the proof of 
 Proposition \ref{sym_contain}(i). \vspace{\baselineskip}

\noindent{\bf (ii).}  
Let $\iota$ be the cardinality of the set of injective elements of $S$ and let $\set{f_{\alpha}}{\alpha<\iota}$ be those injective 
elements. Using transfinite induction for all ordinals $\alpha<\iota$ we may define 
 $$x_{\alpha}, y_{\alpha}\in \Omega f_{\alpha}\setminus \set{x_\beta,y_\beta}{\beta<\alpha}$$ such that  
 $x_{\alpha}\not=y_{\alpha}$. 
Let $T=\set{f\in \Omega^\Omega}{x_{\alpha}f=y_{\alpha}f \:(\forall \alpha<\iota)}$. Then $\genset{S, T}$ is a proper 
subsemigroup of 
$\Omega^\Omega$, since every injective function in $\genset{S, T}$ belongs to $S$ and $|S|\leq |\Omega|$.  Also if 
$\Sigma$ is a transversal of any $f\in T$ such that $|\Omega f|=|\Omega|$, then 
$\set{g|_{\Sigma}}{g\in T}=\Omega^\Sigma$. Hence if $h$ is any injective function in $\Omega^\Omega$ such that 
$\Omega h=\Sigma$, then  $\genset{S, T, h}=\Omega^\Omega$. Hence $\genset{S,T}$, and so $S$, are contained in a 
maximal subsemigroup 
of $\Omega^\Omega$ by part (i). 
\qed\vspace{\baselineskip}

A subgroup $G$ of $\sym(\Omega)$ is \emph{highly transitive} if for all $n\in\N$  and for all 
$(\alpha_1, \alpha_2, \ldots, \alpha_n),$ $(\beta_1, \beta_2,\ldots, \beta_n)\in \Omega^n$, there exists $g\in G$ such that 
$$(\alpha_1g, \alpha_2g, \ldots, \alpha_ng)= (\beta_1, \beta_2,\ldots, \beta_n).$$

We give a new proof of the next theorem using Propositions \ref{sym_contain} and \ref{bergmanshelah}.

\begin{thm}[Macpherson \& Praeger \cite{Macpherson1990ab}]\label{not_dense}
Let $\Omega$ be countably infinite and let $G$ be a subgroup of $\sym(\Omega)$ that is not highly transitive. Then $G$ 
is contained in a maximal subgroup of $\sym(\Omega)$. 
\end{thm} 
\proof
If $G$ is any subgroup of $\sym(\Omega)$, then $G$ satisfies one of the following conditions:
\begin{enumerate}
\item[\rm (a)]  $G_{(\Sigma)}$  has an infinite orbit for all finite $\Sigma\subseteq \Omega$;
\item[\rm (b)]  there exists finite $\Sigma\subseteq \Omega$ such that every orbit of $G_{(\Sigma)}$ is finite.
\end{enumerate}

Suppose that $G$ is a subgroup of $\sym(\Omega)$ that is not highly transitive. 
 If $G$ satisfies (a), then, by Proposition \ref{bergmanshelah}, $G$ is contained in a maximal 
subgroup.

If $G$ satisfies (b), then we may assume without loss of generality that $\Sigma=\{0,1,\ldots, m-1\}$. Since every orbit of 
$G_{(\Sigma)}$ is finite, every orbit of $G$ on $\Omega^{m+1}$ contains only finitely many tuples of the form 
$(0,1,\ldots, m-1, n)$ where $n\in\nat$. But there are infinitely many such tuples and so $G$ has infinitely many orbits on 
$\Omega^{m+1}$. Thus Proposition \ref{sym_contain}(iii) implies that $G$ is contained in a maximal subgroup of 
$\sym(\Omega)$. \qed

%%%%%%%%%%%%%%%%%%%%%%%%%
%%%%%%%%%%%%%%%%%%%%%%%%%

\section{Generating pairs}\label{corollaries}

In \cite[Theorem 3.3]{Howie1998aa} it is shown that $\sym(\Omega)$ has relative rank $2$ in $\Omega^\Omega$; 
that is, there exist 
$f, g\in \Omega^\Omega$ such that $\genset{\sym(\Omega), f,g}=\Omega^\Omega$. Those pairs $f,g\in \Omega^\Omega$ 
satisfying this property are completely classified  in the case that $|\Omega|$ is a 
regular cardinal; see 
\cite[Theorem 4.1]{Howie1998aa}. In this section, we recover this classification as a corollary to Theorem \ref{dicks}, and
 extend it to sets of arbitrary cardinality. Furthermore,  we obtain  analogous results where $\sym(\Omega)$ is replaced
 by the stabiliser of a finite set, an ultrafilter, or a finite partition. 
We require the following straightforward lemma to obtain the corollaries in this section. 

%%%%%%%%%%%%%%%%%%%%%%%%%

\begin{lem}\label{supereasy}
Let $G$ be a subgroup  of $\sym(\Omega)$ containing $\sym(\Omega)_{(\Sigma)}$ for some $\Sigma\subseteq \Omega$ 
such that $|\Omega\setminus \Sigma|=|\Omega|$ and let $H$ be any subset of $\Omega^\Omega$. 
Then $\genset{G, H}=\Omega^\Omega$ if and only if $H$ is not contained in any maximal subsemigroup of 
$\Omega^\Omega$ that contains $G$. 
\end{lem}
\proof
If $H$ is a subset of a maximal subsemigroup of $\Omega^\Omega$ containing $G$, then 
$\genset{G, H}$ is contained in that semigroup, and so 
$\genset{G, H}\not=\Omega^\Omega$.
For the converse,  \cite[Lemma 2.4]{Macpherson1990aa} states that if $U$ is any subgroup of $\sym(\Omega)$ containing $\sym(\Omega)_{(\Gamma)}$ for some moiety $\Gamma$ of $\Omega$, then there exists $x\in \sym(\Omega)$ such that 
$\genset{U, x, x^{-1}}=\sym(\Omega)$. 
 It follows  that $G$ has 
finite relative rank
in $\sym(\Omega)$. Hence, since $\sym(\Omega)$ has finite relative rank in $\Omega^\Omega$ 
(by \cite[Theorem 3.3]{Howie1998aa} as stated above),   
any  subsemigroup of 
$\Omega^\Omega$ containing $G$ has finite relative rank in $\Omega^\Omega$. 
It follows by Proposition \ref{contain}(i) that any proper subsemigroup of $\Omega^\Omega$ containing $G$ is
contained in a maximal subsemigroup of $\Omega^\Omega$. Therefore if $H$ is not contained in any 
maximal subsemigroup containing $G$, then 
$\genset{G, H}=\Omega^\Omega$.
\qed\vspace{\baselineskip}

The following corollary of Theorem \ref{dicks} and Lemma \ref{supereasy} extends 
 \cite[Theorem 4.1]{Howie1998aa}. 

%%%%%%%%%%%%%%%%%%%%%%%%%

\begin{cor}\label{hhr_new}
Let $\Omega$ be any infinite set and let $f,g\in\Omega^\Omega$. 
Then $\genset{\sym(\Omega), f, g}=\Omega^\Omega$ if and only if (up to renaming $f$ and $g$) $f$ is injective, 
$d(f)=|\Omega|$, $g$ is surjective, and either: 
\begin{enumerate}
\item[\rm (i)]  $|\Omega|$ is regular and $k(g, |\Omega|)=|\Omega|$; or 
\item[\rm (ii)] $|\Omega|$ is singular and $k(g, \nu)=|\Omega|$ for all $\nu<|\Omega|$.
\end{enumerate}
\end{cor}
\proof  
 By Lemma \ref{supereasy}, it suffices to show that none of the 
maximal subsemigroups in Theorem \ref{dicks} contains both $f$ and $g$ if and only if (up to renaming $f$ and $g$) $f$ 
is injective, 
$d(f)=|\Omega|$, $g$ is surjective, and either: 
\begin{enumerate}
\item[\rm (i)]  $|\Omega|$ is regular and $k(g, |\Omega|)=|\Omega|$; or 
\item[\rm (ii)] $|\Omega|$ is singular and $k(g, \nu)=|\Omega|$ for all $\nu<|\Omega|$.
\end{enumerate}

For the direct implication, if $\{f,g\}$ is not contained in $S_1\cup S_2$, then (up to renaming $f$ and $g$) $f$ is injective 
and $g$ is 
surjective. Hence $g\in S_4(\mu)$ and so $f\not\in S_4(\mu)$ for all $\mu$. It follows that  $d(f)=|\Omega|$. 
 Regardless of the cardinality of $\Omega$, $f$ belongs to $S_5'\subseteq S_5$. So, if $|\Omega|$ is regular, 
 then $S_5$ is maximal, $g\not\in S_5$, and $k(g, |\Omega|)=|\Omega|$. Similarly, if $|\Omega|$ is singular, then 
 $k(g, \nu)=|\Omega|$ for all $\nu<|\Omega|$. 
 
 For the converse implication, it is easy to verify that 
$f\not\in S_2\cup S_4(\mu)$ and $g\not\in  S_1\cup S_3(\mu)$ for all infinite cardinals $\mu$ not greater than $|\Omega|$. 
If $|\Omega|$ is regular, then $g\not \in S_5$ and if  
 $|\Omega|$ is singular, then $g\not\in S_5'$. 
 \qed\vspace{\baselineskip}

Analogous to Corollary \ref{hhr_new}, we can deduce from Theorem \ref{stab_finite} a characterisation of those $f,g\in \Omega^\Omega$ that together with the pointwise stabiliser of a finite set  generate $\Omega^\Omega$.

\begin{cor}\label{stab_gen_pairs} 
Let $\Omega$ be any infinite set, let $\Sigma$ be a non-empty finite subset of $\Omega$, and let $f,g\in \Omega^\Omega$.
Then the following are equivalent:
\begin{enumerate}
\item[\rm (I)] $\genset{\sym(\Omega)_{(\Sigma)},f,g}=\Omega^\Omega$;
\item[\rm (II)] $\genset{\sym(\Omega)_{\{\Gamma\}},f,g}=\Omega^\Omega$ for all $\Gamma\subseteq \Sigma$;
\item[\rm (III)] $f$ and $g$ satisfy the conditions of Corollary \ref{hhr_new} and for all  non-empty 
$\Gamma\subseteq \Sigma$ one of the following holds:
\begin{enumerate}
\item[\rm (i)] $\Gamma f \not \subseteq \Gamma$ and $\Gamma g^{-1} \not \subseteq \Gamma$;
\item[\rm (ii)] $\Gamma g \not \subseteq \Gamma$, $\Gamma g^{-1} \not \subseteq \Gamma$, and 
$|\Gamma g|=|\Gamma|$;
\item[\rm (iii)] $\Gamma f \not \subseteq \Gamma$, $\Gamma f^{-1} \not \subseteq \Gamma$, and 
$\Gamma \subseteq \Omega f$.
\end{enumerate}
\end{enumerate}
\end{cor}
\proof

%%%%%%%%%%%%%%%%%%%%%%%%%

(I) $\Rightarrow$ (II) This implication follows immediately since 
$\sym(\Omega)_{(\Sigma)}\subseteq \sym(\Omega)_{\{\Gamma\}}$ for all $\Gamma\subseteq \Sigma$.
\vspace{\baselineskip}

%%%%%%%%%%%%%%%%%%%%%%%%%

(II) $\Rightarrow$ (III) Let $\Gamma$ be any non-empty subset of $\Sigma$.  Since 
$\genset{\sym(\Omega), f, g}\supseteq \genset{\sym(\Omega)_{\{\Gamma\}},f,g}=\Omega^\Omega$, 
clearly $f$ and $g$ satisfy the conditions of 
Corollary \ref{hhr_new} and $\{f, g\}$ is not contained in any proper 
 subsemigroup of $\Omega^\Omega$ containing $\sym(\Omega)_{\{\Gamma\}}$. In particular, $\{f, g\}$ is not a subset of 
$F_1(\Gamma, |\Omega|^+)$ or $F_2(\Gamma, |\Omega|^+)$.  
If $f\not\in F_2(\Gamma, |\Omega|^+)$ and $g\not\in F_1(\Gamma, |\Omega|^+)$, then $\Gamma f\not \subseteq \Gamma$ and 
 $\Gamma g^{-1}\not \subseteq \Gamma$ and so (i) holds. 
If $g\not \in F_2(\Gamma, |\Omega|^+)$, then $|\Gamma g|=|\Gamma|$ and $\Gamma g\not\subseteq \Gamma$. But $g$ is 
surjective and so $\Gamma g^{-1}\not\subseteq \Gamma$, and so (ii) holds. 
If $f\not\in F_1(\Gamma, |\Omega|^+)$, then $\Gamma\subseteq \Omega f$ and 
$\Gamma f^{-1}\not \subseteq \Gamma$. Hence, since $f$ is injective, $\Gamma f\not\subseteq \Gamma$ and 
(iii) holds. \vspace{\baselineskip}

%%%%%%%%%%%%%%%%%%%%%%%%%

(III) $\Rightarrow$ (I). 
%The pointwise stabiliser $\sym(\Omega)_{(\Sigma)}$ has finite relative rank in $\sym(\Omega)$  by Proposition \ref{}. 
 Again by Lemma \ref{supereasy} it suffices to show that  none of the  maximal subsemigroups in Theorems 
 \ref{dicks} and \ref{stab_finite} contain both $f$ and 
$g$.

Since $f$ and $g$ satisfy the conditions of Corollary \ref{hhr_new}, it follows that they are not contained in any of the semigroups from Theorem \ref{dicks}.  Moreover, the same conditions
imply that  $f\not\in F_2(\Gamma, \mu)$ and $g\not\in F_1(\Gamma, \mu)$ for all $\mu\leq |\Omega|$. 
If (i) holds, then $f\not\in F_2(\Gamma, |\Omega|^{+})$ and $g\not\in F_1(\Gamma, |\Omega|^{+})$. If (ii) holds, then $g\not\in F_1(\Gamma, |\Omega|^{+})\cup F_2(\Gamma, |\Omega|^{+})$, and if (iii) holds, $f\not\in 
F_1(\Gamma, |\Omega|^{+})\cup F_2(\Gamma, |\Omega|^{+})$.
\qed \vspace{\baselineskip}
 
 %%%%%%%%%%%%%%%%%%%%%%%%%
 
 In the next corollary we characterise the pairs of functions that together with the stabiliser of an ultrafilter generate 
 $\Omega^\Omega$. The statement of this result is similar to that of Corollary \ref{stab_gen_pairs}.
 
\begin{cor} \label{ultra_pairs}
Let $\Omega$ be any infinite set, let $\F$ be an ultrafilter on $\Omega$, and let $f,g\in \Omega^\Omega$.  
Then $\genset{\sym(\Omega)_{\{\F\}}, f,g}=\Omega^\Omega$ if and only if $f$ and $g$ satisfy the conditions of Corollary 
\ref{hhr_new} and there exist $\Sigma\in \F$ and $\Gamma\not\in \F$ such that one of the following holds:
\begin{enumerate}
\item[\rm(i)] $\Sigma f\not\in \F$ and $\Gamma g\in \F$;
\item[\rm(ii)] $\Sigma g\not\in \F$, $c(g|_{\Sigma})=0$, and $\Gamma g\in \F$;
\item[\rm(iii)] $\Sigma f\not\in \F$ and $\Gamma f\in \F$;
\end{enumerate}
\end{cor}
\proof 
If $\F$ is a principal ultrafilter, say generated by $\{\alpha\}$, then $\sym(\Omega)_{\{\F\}}=\sym(\Omega)_{(\{\alpha\})}$ and  the result follows by Corollary \ref{stab_gen_pairs}. 

Suppose that $\F$ is a non-principal ultrafilter. Recall that $\sym(\Omega)_{\{\F\}}$ contains the pointwise stabiliser of any $\Sigma\in \F$.
Let $\kappa$ denote the least cardinality of a set in $\F$. If $\kappa<|\Omega|$, then there exists $\Sigma\in \F$ such that 
$|\Sigma|=\kappa$ and so $|\Omega\setminus \Sigma|=|\Omega|$. Suppose that $\kappa=|\Omega|$.  Then if $\Sigma\in \F$ is such 
that  $|\Omega\setminus \Sigma|<|\Omega|$ and $\Gamma$ is a moiety of 
$\Sigma$ (and hence in $\Omega$), then either $\Gamma\in\F$ or $\Sigma\setminus\Gamma\in \F$ (since otherwise 
$\Omega\setminus \Sigma\in \F$ and $|\Omega\setminus \Sigma|<\kappa$, which is a contradiction). In either case, it follows that
$\sym(\Omega)_{\{\F\}}$ contains the pointwise stabiliser of some $\Sigma\in \F$ such that $|\Omega\setminus \Sigma|=|\Omega|$. 

Therefore by Lemma \ref{supereasy} and Theorem \ref{ultra}, it follows that 
$\genset{\sym(\Omega)_{\{\F\}}, f,g}=\Omega^\Omega$ if and only if $\{f, g\}$ is 
not a subset of  $U_1(\F, \mu)\cup U_2(\F,\mu)$ for any cardinal $\mu$ such that $\kappa< \mu\leq |\Omega|^+$. \vspace{\baselineskip}

($\Rightarrow$) Since 
$\genset{\sym(\Omega), f, g}\supseteq \genset{\sym(\Omega)_{\{\F\}},f,g}=\Omega^\Omega$, 
clearly $f$ and $g$ satisfy the conditions of 
Corollary \ref{hhr_new}.  From the discussion above, it follows that, in particular, 
$\{f, g\}\not\subseteq U_1(\F, |\Omega|^+)\cup U_2(\F,|\Omega|^+).$
 If $g\not\in U_1(\F, |\Omega|^+)$ and $f\not\in U_2(\F, |\Omega|^+)$, then there exists $\Sigma\in \F$ and 
$\Gamma\not\in\F$ such that $\Sigma f\not\in \F$ and $\Gamma g\in \F$, in which case (i) holds. 
If $f\not\in U_1(\F, |\Omega|^+)$, then there exists $\Gamma\not\in\F$ such that $\Gamma f\in \F$. It follows that 
$\Omega \setminus\Gamma\in \F$  and $(\Omega\setminus \Gamma)f\subseteq (\Omega\setminus \Gamma f)\not\in \F$,
and so $(\Omega\setminus \Gamma)f\not\in \F$, which implies  (iii) holds. 
If $g\not\in U_2(\F, |\Omega|^+)$, then there exists $\Sigma\in \F$ such that $\Sigma g\not\in \F$ and $c(g|_{\Sigma})=0$. 
But $g$ is surjective and so $(\Omega\setminus \Sigma)g\supseteq \Omega\setminus (\Sigma g)\in \F$. Thus  
 $\Omega\setminus \Sigma\not\in \F$ but $(\Omega\setminus \Sigma)\in \F$ and so (ii) holds. \vspace{\baselineskip}

($\Leftarrow$)  If $\mu\leq |\Omega|$, then, since $f$ and $g$ satisfy the conditions of Corollary \ref{hhr_new}, 
it follows that $f\not\in U_2(\F, \mu)$ and $g\not\in U_1(\F, \mu)$. Hence it suffices to show that 
$\{f,g\}\not\subseteq U_1(\F, |\Omega|^+)\cup U_2(\F,|\Omega|^+)$ if one of  (i), (ii), or (iii) holds.
It is easy verify that if (i) holds, then $f\not\in  U_2(\F, |\Omega|^+)$ and $g\not\in  U_1(\F,|\Omega|^+)$;
and if (ii) or (iii) holds, then $g\not\in U_1(\F, |\Omega|^+)\cup U_2(\F,|\Omega|^+)$ or 
$f\not\in U_1(\F, |\Omega|^+)\cup U_2(\F,|\Omega|^+)$, respectively.
\qed\vspace{\baselineskip}

As above, Theorem \ref{almost_stab} can be used 
to characterise those $f,g\in \Omega^\Omega$ that together with either $\stab(\P)$ or $\astab(\P)$ generate 
$\Omega^\Omega$. 

%%%%%%%%%%%%%%%%%%%%%%%%%

\begin{cor}\label{almost_stab_pairs}
Let $\Omega$ be any infinite set, let $\mathcal{P}=\{\Sigma_0, \Sigma_1, \ldots, \Sigma_{n-1}\}$, $n\geq 2$, be a finite partition of $\Omega$, and let $f,g\in \Omega^\Omega$. Then the following are equivalent:
\begin{enumerate}
\item[\rm (I)] $\genset{\stab(\P), f, g}=\Omega^\Omega$; 
\item[\rm (II)]  $\genset{\astab(\P), f, g}=\Omega^\Omega$;
\item[\rm (III)] $f$ and $g$ satisfy the conditions of Corollary \ref{hhr_new} and  one of the following holds:
\begin{enumerate}
\item[\rm (i)] $\rho_f,\rho_g\not\in\sym(n)$;
\item[\rm (ii)] $\rho_f\not\in \sym(n)$ and $\rho_f^{-1}$ is total;
\item[\rm (iii)] $\rho_g\not\in \sym(n)$ and $\rho_g$ is total.
\end{enumerate}
\end{enumerate}
\end{cor}
\proof 
(I) $\Rightarrow$ (II). This implication follows immediately since $\stab(\P)$ is a subgroup of $\astab(\P)$.
\vspace{\baselineskip}

(II) $\Rightarrow$ (III). If $\Sigma=\Sigma_1\cup\cdots \cup \Sigma_{n-1}$, then $\stab(\P)$, and hence $\astab(\P)$, contains the pointwise stabiliser of $\Sigma$
in $\sym(\Omega)$. 
Hence by Lemma \ref{supereasy}, $\genset{\astab(\P), f, g}=\Omega^\Omega$ 
implies that $\{f,g\}$ is not a subset of $A_1(\P)$ or $A_2(\P)$.  If $f\not\in A_1(\P)$ and $g\not\in A_2(\P)$, then 
$\rho_f, \rho_g\not \in \sym(n)$ and (i) holds.  If $f\not\in A_2(\P)$, then $\rho_f\not\in \sym(n)$ and $\rho_f^{-1}$ 
is total and (ii) holds. If $g\not\in A_1(\P)$, then $\rho_g\not\in \sym(n)$ and $\rho_g$ is total and we are in case 
(iii).\vspace{\baselineskip}

(III) $\Rightarrow$ (I).
Again by Lemma \ref{supereasy},  to prove that $\genset{\stab(\P), f, g}=\Omega^\Omega$, it suffices to 
show that  none of the  maximal subsemigroups in Theorems \ref{dicks} and \ref{almost_stab} contain both $f$ and 
$g$. 

Since $f$ and $g$ satisfy the conditions of Corollary \ref{hhr_new}, it follows that they are not contained in any of the 
semigroups from Theorem \ref{dicks}. If (i) holds, then $f\not\in A_1(\P)$ and $g\not\in A_2(\P)$; if (ii) holds, then, since $f$ 
is injective, $\rho_f$ is total and so $f\not\in A_1(\P)\cup A_2(\P)$; and if (iii) holds, then, since $g$ is surjective, 
$\rho_g^{-1}$ is total and so $g\not \in A_1(\P)\cup A_2(\P)$.
\qed

\setcounter{main}{0}

%%%%%%%%%%%%%%%%%%%%%%%%%%%%
%%%%%%%%%%%%%%%%%%%%%%%%%%%%

\section{Inverses and parameters of mappings}\label{mundane}

In this section we present several technical results, which we will use repeatedly throughout the paper. 

We begin by considering the semigroup theoretic inverses of mappings in $\Omega^\Omega$. Roughly 
speaking, the proofs of the main theorems are in two parts and Corollary \ref{dual_cor} will imply that
 one part is a corollary of the other. More precisely, the majority of the proof of, say, Theorem \ref{stab_finite}, consists of 
 showing the  following. 
 If $U$ is a subsemigroup of $\Omega^\Omega$ that is not 
contained in any of the semigroups listed in Theorems \ref{dicks} or \ref{stab_finite}, but that does contain the 
 stabiliser of 
a non-empty finite subset of $\Omega$, then  $U=\Omega^\Omega$. 
The stabiliser contains the symmetric group on an infinite subset $\Sigma$ of $\Omega$. The 
two parts of the proof, referred to above, are to construct an injective mapping in $U$ with image contained in $\Sigma$ 
and a surjective mapping in $U$ mapping $\Sigma$ onto $\Omega$.  Using Corollary \ref{dual_cor},
the existence of the surjective mapping is a consequence of the existence of the injective mapping. The proofs of 
Theorems \ref{dicks}, \ref{ultra}, and \ref{almost_stab} follow a similar strategy. 

If $S$ is a semigroup and $s\in S$, then $t\in S$ is an \emph{inverse} of $s$ if $sts=s$ and $tst=t$. 
Clearly, $t$ is an inverse for $s$ if and only if $s$ is an inverse for $t$. 
If $f,f'\in \Omega^\Omega$, then $f'$ is an inverse for $f$ if and only if 
$\Omega f'$ is a transversal of $f$ and $ff'$ is the identity on $\Omega f'$. Note that if $f,f'\in \Omega^\Omega$ 
are inverses, then $c(f)=d(f')$.

In general, the composition $g'f'$ of inverses of $g$ and $f$ is not an inverse of the composite $fg$. 
However, for certain composites $g'f'$ is an inverse of $fg$. 

%%%%%%%%%%%%%%%%%%%%%%%%%

\begin{lem}\label{dual_lemma}
Let $u_0, u_1, \ldots, u_n\in \Omega^\Omega$ be arbitrary and let $u_i'$ be an inverse of $u_i$ for all $i\in \{0,1,\ldots, n\}.$ 
If $\Omega u_0' u_1'\cdots u_{i-1}' \subseteq \Omega u_{i}$
for all $ i \in\{1,\ldots, n\}$, then
$u_0'u_1'\cdots u_n'$ and $u_n\cdots u_1 u_0$ are inverses.
\end{lem}
\proof
We show that $u_n\cdots u_1 u_0$ is an inverse of $u_0'u_1'\cdots u_n'$ by showing that 
$\Omega u_n \cdots u_1 u_0$ is a transversal of $u_0'u_1' \cdots u_n'$ and $u_0'u_1' \cdots u_n'u_n \cdots u_1u_0$ 
is the identity on $\Omega u_n \cdots u_1 u_0$.

Since $\Omega u_0' u_1'\cdots u_{i-1}'$ is contained in the transversal $\Omega u_{i}$ of $u_{i}'$ for all 
$i \in\{1,\ldots, n\}$, it follows that $u_1'\cdots u_n'$ is injective on $\Omega u_0'$. Hence the transversal 
$\Omega u_0$ of $u_0'$ is also a transversal of $u_0'u_1' \cdots u_n'$. 

If $x\in \Omega$, then $ x u_0' u_1'\cdots u_{i-1}' \in \Omega u_{i}$ for all $i\in\{1,\ldots, n\}$. 
Since $u_{i}'u_{i}$ is the identity on $\Omega u_{i}$, it follows that 
$ x u_0' u_1'\cdots u_{i-1}' u_{i}' u_{i}=x u_0' u_1'\cdots u_{i-1}'$ for all $x\in \Omega$. 

Applying this $n$ times we obtain
 $$x u_0'u_1'\cdots u_n' u_n \cdots u_1 u_0 = x u_0' u_0$$ 
for all $x\in \Omega$.
In particular,  if $x \in \Omega u_0$, then $x u_0'u_1'\cdots u_n' u_n \cdots u_1 u_0 = x$.

Certainly, $\Omega u_n\cdots u_0\subseteq \Omega u_0$ and since $u_0'u_1'\cdots u_n' u_n \cdots u_1 u_0$
is the identity on $\Omega u_0$, it follows that 
$\Omega u_n\cdots u_0\supseteq \Omega u_0'u_1'\cdots u_n' u_n \cdots u_1 u_0\supseteq \Omega u_0$.
\qed\vspace{\baselineskip}

%%%%%%%%%%%%%%%%%%%%%%%%%

\begin{de}\label{assign}
Let $V\subseteq \Omega^\Omega$ and let $\Lambda:V \to \mathcal{P}(\Omega)$ be such that $\Lambda(v)$ is a transversal of $v$ for all $v\in V$. We refer to such a $\Lambda$ as an \emph{assignment of transversals} for $V$. Then  the set of products $v_0v_1\cdots v_n\in \genset{V}$ such that $v_i\in V$ and $\Omega v_0\cdots v_{i-1}\subseteq \Lambda(v_i)$ for all 
$i\in \{1, \ldots, n\}$ is denoted by $\mutt{V}{\Lambda}$. 
\end{de}

%%%%%%%%%%%%%%%%%%%%%%%%%

\begin{cor}\label{dual_cor}
Let $U\subseteq \Omega^\Omega$, let $u'\in \Omega^\Omega$ be an inverse of $u$ for every $u\in U$,  let 
$U'=\set{u'}{u\in U}$, and let $\Lambda:U'\to \mathcal{P}(\Omega)$  be defined by $\Lambda(u')=\Omega u$. 
Then every element of $\mutt{U'}{\Lambda}$ has an inverse in $\genset{U}$. 
\end{cor}
\proof
Let $u_0'u_1' \cdots, u_n'\in \mutt{U'}{\Lambda}$. Then $\Omega u_0' u_1'\cdots u_{i-1}' \subseteq \Lambda(u_i')=\Omega u_{i}$ for all $ i \in\{1,\ldots, n\}$. Thus, by Lemma \ref{dual_lemma}, 
$u_n\cdots u_0\in \genset{U}$ is an inverse of $u_0'\cdots u_n'$. 
\qed\vspace{\baselineskip}

%%%%%%%%%%%%%%%%%%%%%%%%%

We will make repeated use of the following lemma,  which is similar to Lemma 2.1 in \cite{Howie1998aa}.

%%%%%%%%%%%%%%%%%%%%%%%%%

\begin{lem}\label{tech}
Let $\Omega$ be any infinite set, let $\mu$ be an infinite cardinal such that $\mu\leq |\Omega|$, and let $f,g\in\Omega^\Omega$. Then the following hold:
\begin{itemize}
\item[\rm (i)] if $\mu$ is a regular cardinal, then $k(fg, \mu)\leq k(f, \mu)+k(g, \mu)$;
\item[\rm (ii)] $d(g)\leq d(fg)\leq d(f)+d(g)$;
\item[\rm (iii)] if $g$ is injective (i.e. $c(g)=0$), then $d(fg)=d(f)+d(g)$;
\item[\rm (iv)] $c(f)\leq c(fg)\leq c(f)+c(g)$;
\item[\rm (v)]  if $f$ is surjective (i.e. $d(f)=0$), then $c(fg)=c(f)+c(g)$;
\item[\rm (vi)] if $c(g)<\mu\leq d(f)$, then $d(fg)\geq \mu$;
\item[\rm (vii)]  if $d(f)<\mu\leq c(g)$, then $c(fg)\geq \mu$.
\end{itemize}
\end{lem} 
\proof 
{\bf (i).}  Let $\alpha\in\Omega$. Then 
$$\alpha (fg)^{-1}=\bigcup_{\beta\in\alpha g^{-1}} \beta f^{-1}.$$
If $|\alpha g^{-1}|<\mu$ and $|\beta f^{-1}|<\mu$ for all $\beta\in \alpha g^{-1}$, then, since $\mu$ is regular, $|\alpha (fg)^{-1}|<\mu$. Hence 
\begin{align*}
k(fg, \mu) &=|\set{\alpha\in\Omega}{|\alpha(fg)^{-1}|\geq \mu}| \\
&\leq |\set{\alpha\in\Omega}{(\exists \beta\in \alpha g^{-1})\,(|\beta f^{-1}|\geq \mu)}|+|\set{\alpha\in \Omega}{|\alpha g^{-1}|\geq \mu}| \\
&\leq k(f, \mu)+k(g, \mu),
\end{align*}
as required.
\vspace{\baselineskip}

{\bf (ii).} It is straightforward to see that
$$\Omega\setminus \Omega g\subseteq \Omega\setminus \Omega fg\subseteq (\Omega\setminus \Omega f)g\cup (\Omega\setminus \Omega g)$$
and so $d(g)\leq d(fg)\leq |(\Omega\setminus \Omega f)g|+d(g)\leq d(f)+d(g)$. \vspace{\baselineskip}

{\bf (iii).} If $c(g)=0$, then 
$$\Omega\setminus \Omega fg= (\Omega\setminus \Omega g)\cup (\Omega g\setminus \Omega fg)=(\Omega\setminus \Omega f)g\cup (\Omega\setminus \Omega g)$$
and $|(\Omega\setminus \Omega f)g|=|\Omega\setminus \Omega f|=d(f).$ Hence $d(fg)=d(f)+d(g)$, as required.
\vspace{\baselineskip}

{\bf (iv).} Let $\Sigma\subseteq \Omega$ be a transversal of $f$. Then there exists $\Sigma^{\prime}\subseteq \Sigma$ such that $\Sigma^{\prime}$ is a transversal of $fg$. Hence $c(f)\leq c(fg)$. Also $c(fg)=|\Omega\setminus \Sigma^{\prime}|=|\Omega\setminus \Sigma|+|\Sigma\setminus \Sigma^{\prime}|=c(f)+|\Sigma\setminus \Sigma^{\prime}|$, and so it suffices to show that $|\Sigma\setminus \Sigma^{\prime}|\leq c(g)$. Let $\Gamma$ be any transversal of $g$ such that $\Sigma^{\prime}f\subseteq \Gamma$. If $\alpha\in \Sigma\setminus \Sigma^{\prime}$, then there exists $\beta\in \Sigma^{\prime}$ such that $(\alpha)fg=(\beta)fg$. Since $f$ is injective on $\Sigma$,  $\alpha f\not=\beta f$. But $\beta f\in\Sigma^{\prime}f\subseteq \Gamma$ and so $\alpha f\not\in \Gamma$.
Thus $(\Sigma \setminus \Sigma^{\prime})f\subseteq \Omega\setminus \Gamma$ and so $|\Sigma \setminus \Sigma^{\prime}|=|(\Sigma \setminus \Sigma^{\prime})f|\leq |\Omega\setminus \Gamma|=c(g)$, as required.\vspace{\baselineskip}

{\bf (v).} Let $\Sigma, \Sigma'$, and $\Gamma$ be as in part (iv).  If $d(f)=0$, then $\Sigma f=\Omega$. But we saw in part (iv) that $(\Sigma\setminus \Sigma^{\prime})f\subseteq \Omega\setminus \Gamma$ and $\Sigma^{\prime}f\subseteq \Gamma$ and so in this case $(\Sigma\setminus \Sigma^{\prime})f= \Omega\setminus \Gamma$ and $\Sigma^{\prime}f= \Gamma$. Therefore $c(fg)=|\Omega\setminus \Sigma^{\prime}|=|\Omega\setminus \Sigma|+ |\Sigma\setminus \Sigma^{\prime}|=|\Omega\setminus \Sigma|+ |\Omega\setminus \Gamma|=c(f)+c(g)$. \vspace{\baselineskip}

{\bf (vi).} %Since $\Omega\setminus \Omega fg\subseteq(\Omega\setminus \Omega g)\cup (\Omega\setminus\Omega f)g$, it suffices to show that $|(\Omega\setminus \Omega f)g|\geq \mu$.
 If $\Sigma$ is any transversal of $g$, then, by assumption, $|\Omega\setminus \Sigma|=c(g)<\mu$ and $|\Omega\setminus \Omega f|=d(f)\geq \mu$.  Hence $|\Sigma\cap (\Omega\setminus \Omega f)|\geq \mu$. If $\alpha\in \Sigma\cap (\Omega\setminus \Omega f)$ is such that $\alpha g\in\Omega fg$, then there exists $\beta\in \Omega f$ such that $\alpha g=\beta g$. So, since $c(g)<\mu$, 
 $$|\set{\alpha\in\Sigma\cap (\Omega\setminus \Omega f)}{\alpha g\in \Omega fg}|< \mu.$$
 Therefore $|\set{\alpha\in\Sigma\cap (\Omega\setminus \Omega f)}{\alpha g\not\in \Omega fg}|\geq \mu$ and so
 $$|\Omega\setminus \Omega fg|\geq |\set{\alpha\in\Sigma\cap (\Omega\setminus \Omega f)}{\alpha g\not\in \Omega fg}g|=|\set{\alpha\in\Sigma\cap (\Omega\setminus \Omega f)}{\alpha g\not\in \Omega fg}|\geq \mu,$$
 as required.
\vspace{\baselineskip}

{\bf (vii).}  As in the proof of (iv), let $\Sigma$ be a transversal of $f$, let $\Sigma^{\prime}\subseteq \Sigma$ be a 
transversal of $fg$, and let $\Gamma$ be a transversal of $g$ such that $\Sigma^{\prime}f\subseteq \Gamma$. By 
assumption, $|\Omega\setminus \Sigma f|=|\Omega\setminus \Omega f|=d(f)<\mu$ and 
$|\Omega\setminus \Gamma|=c(g)\geq \mu$. Hence $|\Sigma f\cap (\Omega\setminus \Gamma)|\geq \mu$.  Since 
$\Sigma^{\prime}f\subseteq \Gamma$ 
and, again as in the proof of (iv), $(\Sigma\setminus \Sigma^{\prime})f\subseteq \Omega \setminus \Gamma$, it follows 
that $\Sigma f\cap (\Omega\setminus \Gamma)=(\Sigma\setminus \Sigma^{\prime})f$. Thus
$$\mu\leq |\Sigma f\cap (\Omega\setminus \Gamma)|=|(\Sigma\setminus \Sigma^{\prime})f|=|\Sigma\setminus\Sigma^{\prime}|\leq |\Omega\setminus \Sigma^{\prime}|=c(fg),$$
as required.\qed

%%%%%%%%%%%%%%%%%%%%%%%%%%%%
%%%%%%%%%%%%%%%%%%%%%%%%%%%%

\section{The symmetric group -- the proof of Theorem \ref{dicks}}\label{proof_dicks}

In this section, we give the proof of Theorem  \ref{dicks}.  
We  require the following result from \cite[Theorem 3.3]{Howie1998aa}.

%%%%%%%%%%%%%%%%%%%%%%%%%

\begin{thm}\label{hhr}
Let $\Omega$ be an infinite set and let $f, g\in \Omega^\Omega$ such that $f$ is injective, $g$ is surjective, and  
$d(f)=k(g, |\Omega|)=|\Omega|$. Then $\genset{\sym(\Omega), f, g}=\Omega^\Omega$. 
\end{thm}

Recall that a subset $\Sigma$ of an infinite set $\Gamma$ is a moiety of $\Gamma$ if 
$|\Sigma|=|\Gamma\setminus \Sigma|=|\Gamma|$. 

%%%%%%%%%%%%%%%%%%%%%%%%%

\begin{lem} \label{s3}
Let $\Omega$ be any set of  singular cardinality and let $g\in \Omega^\Omega$ such that $k(g, \mu)=|\Omega|$ for all $\mu<|\Omega|$. Then there exists $a\in \sym(\Omega)$ such that $k(gag, |\Omega|)=|\Omega|$. 
\end{lem}
\proof 
Since $|\Omega|$ is singular, there exist $\kappa<|\Omega|$ and  $\Omega_{\mu}\subseteq \Omega$  such that $|\Omega_{\mu}|<|\Omega|$ for all $\mu<\kappa$ and $\Omega=\bigcup_{\mu<\kappa} \Omega_{\mu}$.  Let $\Sigma$ be a moiety of $\set{\alpha\in\Omega}{|\alpha g^{-1}|\geq \kappa}$, let $\set{\beta(\alpha, \mu)\in \Omega}{\mu<\kappa}\subseteq \alpha g^{-1}$ for all $\alpha\in \Sigma$, where $\beta(\alpha,\mu)\not=\beta(\alpha,\nu)$ if $\mu\not=\nu$, and let $\Sigma^{\prime}=\bigcup_{\alpha\in \Sigma}\set{\beta(\alpha, \mu)\in\Omega}{\mu<\kappa}$. 

We next show that there exists a moiety $\Gamma$ of $\Omega$ such that 
$|\set{\alpha\in\Gamma}{|\alpha g^{-1}|\geq \mu}|=|\Omega|$
for all $\mu<|\Omega|$.
In fact, if $\Omega$  is arbitrarily partitioned into moieties  $\Gamma_1$ and $\Gamma_2$, then one or the other of these sets has the required property. To see this, suppose that there exists $\nu<|\Omega|$  such that 
$|\set{\alpha\in \Gamma_1}{|\alpha g^{-1}|\geq \nu}<|\Omega|$. If $\mu$ is a cardinal such that $\nu\leq \mu< |\Omega|$, then $k(g, \mu)=|\set{\alpha\in \Omega}{|\alpha g^{-1}|\geq \mu}|=|\Omega|.$
But 
$$\set{\alpha\in \Omega}{|\alpha g^{-1}|\geq \mu}=\set{\alpha\in \Gamma_1}{|\alpha g^{-1}|\geq \mu}\cup \set{\alpha\in \Gamma_2}{|\alpha g^{-1}|\geq \mu}$$
and so $|\set{\alpha\in \Gamma_2}{|\alpha g^{-1}\geq \mu}|=|\Omega|$.  So we now fix $\Gamma$ with the above property.

Assume that $\Sigma\times \kappa$ is well-ordered. We define, by transfinite recursion, distinct $\gamma(\alpha, \mu)\in \Gamma$ such that $|\gamma(\alpha, \mu) g^{-1}|\geq |\Omega_\mu|$ for all $(\alpha, \mu)\in \Sigma\times \kappa$ as follows.  Let $(\alpha, \mu)\in \Sigma\times \kappa$ and let $$\Gamma^{\prime}=\set{\gamma(\beta, \nu)}{(\beta, \nu)<(\alpha, \mu)}.$$ Then $|\set{\gamma\in \Gamma}{|\gamma g^{-1}|\geq |\Omega_\mu|}|=|\Omega|$ and $|\Gamma^{\prime}|<|\Sigma\times \kappa|=|\Omega|$. 
So, we may choose $\gamma(\alpha, \mu)$ to be any element in the set $\set{\gamma\in \Gamma\setminus \Gamma^{\prime}}{|\gamma g^{-1}|\geq |\Omega_\mu|}$, which is of cardinality $|\Omega|$. 

\begin{figure}
\begin{center}
 \begin{tikzpicture}[scale=1]
 \verticalline0
 \verticalline3
 \verticalline6
 \verticalline9
 \smallchip01
 \smallchip03
 \smallchip31
 \smallchip33
 \smallchip62
 \smallchip63
 \smallchip64
 \smallchip65
 \smallchip91
 \smallchip93
 \pt3{1.5}
 \pt6{3.5}
 \pt92
 \functionline01{1.5}
 \functionline03{1.5}
 \functionlinedash3{1.5}{3.5}
 \functionline632
 \functionline642
 \functionline653
 \functionline621
 \verticallinearrowleft006{$\Omega$}
 \verticallinearrowright013{$\geq\!|\Omega_\mu|$}
 \verticallinearrowleft313{$\Gamma$}
 \verticallinearrowlefthigh625{$\Sigma g^{-1}$}
 \verticallinearrowright913{$\Sigma$}
 \functionlabel0{g}
 \functionlabel3{a}
 \functionlabel6{g}
 \draw(3,1.5)node[below right]{{\small $\gamma(\alpha,\mu)$}};
 \draw(6,3.5)node[right]{{\small $\beta(\alpha,\mu)$}};
 \draw(9,2)node[below left]{{\small $\alpha$}};
 \end{tikzpicture}
\caption{The composite $gag$ in the proof of Lemma \ref{s3}.}\label{fig1}
\end{center}
\end{figure}
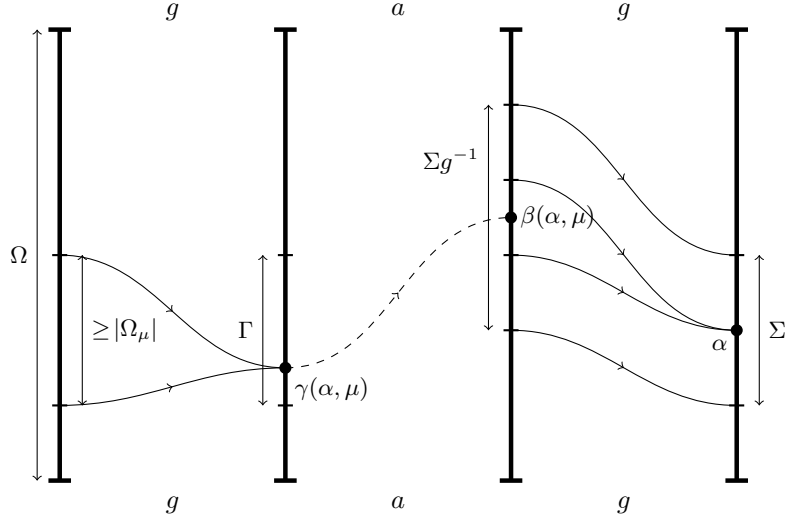

Since $\Gamma$ and $\Sigma'$ are moieties, there exists $a\in\sym(\Omega)$ such that 
$$(\gamma(\alpha, \mu))a=\beta(\alpha, \mu)$$
for all $(\alpha, \mu)\in \Sigma\times \kappa$. Therefore, for any $\alpha\in\Sigma$,
$$\alpha(gag)^{-1}\supseteq \bigcup_{\mu<\kappa}\beta(\alpha, \mu)a^{-1}g^{-1}=\bigcup_{\mu<\kappa}\gamma(\alpha, \mu)g^{-1}$$
and so $|\alpha(gag)^{-1}|\geq|\bigcup_{\mu<\kappa}\Omega_{\mu}|=|\Omega|$; see Figure \ref{fig1}. 
Since $|\Sigma|=|\Omega|$, it follows that $k(gag, |\Omega|)=|\Omega|$, as required. \qed\vspace{\baselineskip}

%%%%%%%%%%%%%%%%%%%%%%%%%%%%

\begin{main}[Gavrilov \cite{Gavrilov1965aa}, Pinsker \cite{Pinsker2005aa}]
Let $\Omega$ be any infinite set.  If $|\Omega|$ is a regular cardinal, then the maximal subsemigroups of $\Omega^\Omega$ containing $\sym(\Omega)$ are: 
\begin{align*}
S_1&=\set{f \in \Omega^\Omega}{c(f)= 0 \text{ or }d(f)>0};\\
S_2&=\set{f \in \Omega^\Omega}{c(f)>0 \text{ or }d(f)=0};\\
S_3(\mu)&=\set{f \in \Omega^\Omega}{c(f)< \mu \text{ or }d(f)\geq \mu};\\
S_4(\mu)&=\set{f \in \Omega^\Omega}{c(f)\geq\mu \text{ or }d(f)<\mu};\\
S_5&=\set{f \in \Omega^\Omega}{k(f, |\Omega|)<|\Omega|};\\ 
\intertext{where $\mu$ is any infinite cardinal not greater than $|\Omega|$.}
\intertext{If $|\Omega|$ is a singular cardinal, then the maximal subsemigroups of $\Omega^\Omega$ containing $\sym(\Omega)$ are $S_1$, $S_2$, $S_3(\mu)$, $S_4(\mu)$ where $\mu$ is any infinite cardinal not greater than $|\Omega|$,  and: 
}
S_5^{\prime}&=\set{f \in \Omega^\Omega}{(\exists \nu<|\Omega|)\,(k(f, \nu)<|\Omega|)}.
\end{align*}
\end{main}

%%%%%%%%%%%%%%%%%%%%%%%%%

If $\Omega$ is any infinite set, then Lemma \ref{tech} can be used to show that $S_1$, $S_2$,  $S_3(\mu)$, and $S_4(\mu)$  are semigroups for all infinite $\mu\leq |\Omega|$. In particular, parts 
(ii), (iii), and (iv) show this for $S_1$;
(ii), (iv), (v) show this for $S_2$;
(ii), (iv), (vi) show this for $S_3(\mu)$;  
(ii), (iv), (vii) show this for $S_4(\mu)$. 
It is also straightforward to verify that none of $S_1$, $S_2$, $S_3(\mu)$, $S_4(\nu)$, with $\mu,\nu\leq|\Omega|$  infinite cardinals, is contained in any of the others. %Note that $S_3(\mu)$ and $S_4(\mu)$ are not contained in $S_3(\nu)$ or $S_4(\nu)$ for any infinite cardinals $\mu,\nu\leq |\Omega|$ such that $\mu\not=\nu$. 

If $|\Omega|$ is regular, then Lemma \ref{tech}(i) shows that $S_5$ is a semigroup. If $|\Omega|$ is singular, then $S_5$ is a generating set for $\Omega^\Omega$, and, in particular,  not a semigroup. Regardless of the nature of $|\Omega|$, $S_5^{\prime}$ is contained in  $S_5$. However, $S_5$ and $S_5^{\prime}$ are not contained in, and do not contain, any of $S_1$, $S_2$,  $S_3(\mu)$, and $S_4(\mu)$ for any $\mu$.
%is S_5^{\prime} always a semigroup? JDM

To show that $S_5^{\prime}$ is a semigroup in the case that $|\Omega|$ is a singular cardinal,  let $f, g\in S_5^{\prime}$. Then there exists $\mu, \nu<|\Omega|$ such that $k(f, \mu)<|\Omega|$ and $k(g, \nu)<|\Omega|$. Let $\kappa=\max\{\mu, \nu\}$. If $\kappa^{+}$ denotes the successor of $\kappa$, then $\kappa^{+}<|\Omega|$ since $|\Omega|$ is singular. Since $k(f,\kappa^{+}), k(g, \kappa^{+})<|\Omega|$ and $\kappa^{+}$ is regular, it follows, by Lemma \ref{tech}(i),  that $k(fg,\kappa^{+})\leq k(f, \kappa^{+})+k(g, \kappa^{+})<|\Omega|$. Hence $fg\in S_5^{\prime}$ and $S_5^{\prime}$ is a semigroup.

We require Lemmas \ref{sym_inj} and \ref{sym_surj} below to complete the proof of Theorem \ref{dicks}; they are stated in far greater generality than required in this section because we will use them again in later sections. 

%%%%%%%%%%%%%%%%%%%%%%%%%

If $a\in \sym(\Omega)$, then we denote the set $\set{\alpha\in\Omega}{\alpha a\not=\alpha}$ by $\supp(a)$ and refer to 
this set as the \emph{support} of $a$.  

\begin{lem} \label{sym_inj}
Let $U$ be a subset of $\Omega^\Omega$, which is not contained in 
$S_2$ or $S_4(\mu)$ for any infinite $\mu\leq |\Omega|$,  let $\Lambda$ be any assignment of transversals for $U$ 
(as defined in Definition \ref{assign}), and let $\kappa$ be any  cardinal such that $\aleph_0\leq \kappa\leq |\Omega|$. 
If $U$ contains an injective $f$ and 
 every  $a\in \sym(\Omega)$ with $\supp(a)\subseteq \Omega f$ and $|\supp(a)|< \kappa$, then there exists an injective 
 $f^*\in \mutt{U}{\Lambda}$ such that $d(f^*)\geq \kappa$ and $\Omega f^*\subseteq \Omega f$. 
\end{lem}
\proof
We prove by transfinite induction that for each cardinal $\mu\leq\kappa$,
\begin{equation}\label{inj} \text{ there exists } f_{\mu}\in \mutt{U}{\Lambda} \text{ such that } f_{\mu}\text{ is injective, }d(f_{\mu})\geq \mu, \text{ and }\Omega f_{\mu}\subseteq \Omega f. 
\end{equation}
If  there exists $f_{\mu}\in \mutt{U}{\Lambda}$  such that  $f_{\mu}$  is injective and $d(f_{\mu})\geq \mu$, then,  since $f_{\mu}$ and $f$ are injective, $f_{\mu}f\in \mutt{U}{\Lambda}$, $f_{\mu}f$ is injective, $d(f_{\mu} f)= d(f_{\mu})+d(f)\geq \mu$ and $\Omega f_{\mu}\subseteq \Omega f$, i.e. \eqref{inj} holds for $f_{\mu}f$. 
Hence it suffices to show that  there exists $f_{\mu}\in \mutt{U}{\Lambda}$  such that  $f_{\mu}$  is injective and $d(f_{\mu})\geq \mu$.

Since $U\not\subseteq S_2$, there exists an injective $h_0\in U\subseteq \mutt{U}{\Lambda}$ such that $d(h_0)>0$. Since $h_0$ is injective, $h_0^n$ belongs to $\mutt{U}{\Lambda}$ and, by Lemma \ref{tech}(iii),  $d(h_0^n)\geq n$  for all $n\in\nat$ . Thus \eqref{inj} holds for all finite $\mu$.

\begin{figure}
\begin{center}
 \begin{tikzpicture}[scale=1]
 \verticalline0
 \verticalline3
 \verticalline6
 \verticalline9

 \chip9{3.5}
 \verticallinearrowleft006{$\Omega$}
 \verticallinearrowleft3{0.5}{3.5}{$\Omega f_{c(h_1)}^2$}
  \verticallinearrowleft356{$\Omega f_{c(h_1)}^2$}
 \verticallinearrowlefthigh3{3.5}5{$ c(h_1)\!\leq$}
 \verticallinearrowright656{$c(h_1)$}
 \verticallinearrowright605{$\Lambda(h_1)$}
 \verticallinearrowright9{3.5}6{$\geq\!\mu$}
  \verticallinearrowright90{3.5}{$\Omega h_1$}
    \verticallinearrowright3{0.5}6{$\Omega f$}
 \functionline066
 \functionline055
 \functionline05{3.5}
 \functionline00{0.5}
 \functionlinealeft365
  \functionlinearight365
 \functionlinealeft354
  \functionlinearight354
 \functionlinealeftdash345
  \functionlinearightdash345
 \functionlinealeftdash356
  \functionlinearightdash356
  \functionlinedash300
  \functionline3{0.5}{0.5}
  \functionline3{3.5}{3.5}
  \functionlinedash344
% \functionlineleft36{1.5}
 %\functionlineleft3{4.5}0
 %\functionlineleft3{1.5}6
 %\functionlineleft30{4.5}
 %\functionlineright36{1.5}
 %\functionlineright3{4.5}0
 %\functionlineright3{1.5}6
 %\functionlineright30{4.5}
\functionline65{3.5}
\functionline642
\functionline6{3.5}1
\functionline6{0.5}{0.5}
\functionlinedash600

 \functionlabel0{f_{c(h_1)}^2}
 \functionlabel3{a^{-1}}
 \functionlabel6{h_1}
  \colorinterval3{4}5{blue}
   \colorinterval656{blue}
 \colorinterval006{red}
 \colorinterval056{red}
  \colorinterval356{red}
    \colorinterval3{0.5}{3.5}{red}
    \colorinterval645{red}
    \colorinterval6{0.5}{3.5}{red}
    
  \colorinterval9{0.5}1{red}

 \colorinterval92{3.5}{red}

 \end{tikzpicture}
\caption{The composition $f_{c(h_1)}^2a^{-1} h_1$ from the proof of Lemma \ref{sym_inj}.}\label{fig2}
\end{center}
\end{figure}
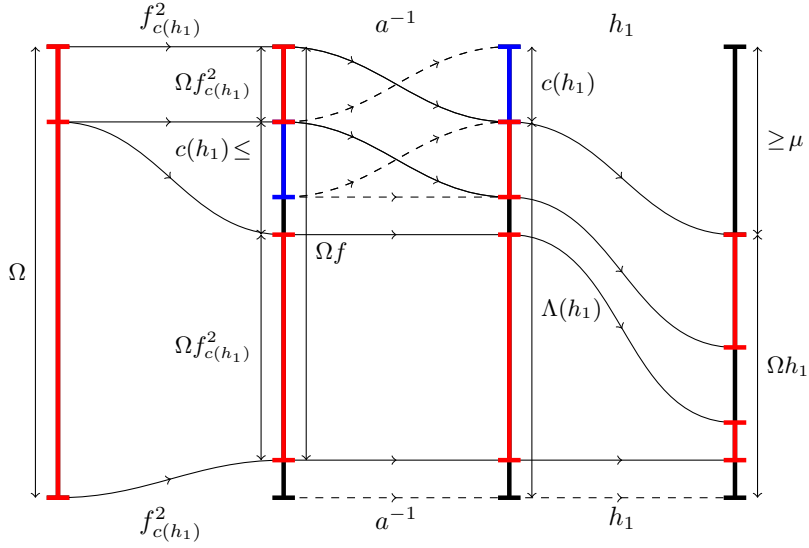

Let $\mu$ be any cardinal such that $\aleph_0\leq \mu\leq \kappa$ and assume that (\ref{inj}) holds for  every cardinal strictly less than $\mu$.  
Since $U\not\subseteq S_4(\mu)$, there exists $h_1\in U$ such that $c(h_1)<\mu\leq d(h_1)$.  By our inductive hypothesis, there exists an injective $f_{c(h_1)}\in \mutt{U}{\Lambda}$ such that $d(f_{c(h_1)})\geq c(h_1)$ and $\Omega f_{c(h_1)}\subseteq \Omega f$. 
Since $f_{c(h_1)}$ is injective, 
$|(\Omega\setminus \Omega f_{c(h_1)})f_{c(h_1)}|=|\Omega\setminus \Omega f_{c(h_1)}|$ and so 
$$|\Omega\setminus \Omega f_{c(h_1)}|=|(\Omega\setminus \Omega f_{c(h_1)})f_{c(h_1)}|=|\Omega f_{c(h_1)}\setminus \Omega f_{c(h_1)}^2|\leq |\Omega f\setminus \Omega f_{c(h_1)}^2|.$$
It follows that 
$$|\Omega f_{c(h_1)}^2\setminus \Lambda(h_1)|\leq |\Omega \setminus \Lambda(h_1)|=c(h_1)\leq 
d(f_{c(h_1)})=|\Omega\setminus \Omega f_{c(h_1)}|\leq |\Omega f\setminus \Omega f_{c(h_1)}^2|.$$
Thus there is $a\in \sym(\Omega)$ such that 
$$\big(\Omega f_{c(h_1)}^2\setminus \Lambda(h_1)\big) a\subseteq \Omega f\setminus \Omega f_{c(h_1)}^2$$
and $\supp(a)\subseteq\big(\Omega f_{c(h_1)}^2\setminus \Lambda(h_1)\big)\cup\big(\Omega f_{c(h_1)}^2\setminus \Lambda(h_1)\big)a$;  see Figure \ref{fig2}.
Hence, since $\Omega f_{c(h_1)}\subseteq \Omega f$, it follows that $\supp(a)\subseteq\ \Omega f$ and  
$$|\supp(a)|\leq 2|\Omega f_{c(h_1)}^2\setminus \Lambda(h_1)|\leq 2|\Omega \setminus \Lambda(h_1)|<\mu\leq \kappa.$$
In particular, $a,a^{-1}\in U$. 

From the definition of $a$, it follows that 
$$\Omega f_{c(h_1)}^2\setminus \Lambda(h_1)\subseteq \big(\Omega f\setminus \Omega f_{c(h_1)}^2\big)a^{-1}$$
 and so $\Omega f_{c(h_1)}^2a^{-1}\subseteq \Lambda(h_1)$; see Figure \ref{fig2}.
 This shows that $f^2_{c(h_1)}a^{-1}h_1\in \mutt{U}{\Lambda}$ is injective, and $d(f^2_{c(h_1)}a^{-1}h_1)\geq d(h_1)\geq\mu$ by Lemma \ref{tech}(ii).  It follows that $f_{\mu}:=f^2_{c(h_1)}a^{-1}h_1\in \mutt{U}{\Lambda}$, and so \eqref{inj} holds for $\mu$.  
\qed
%%%%%%%%%%%%%%%%%%%%%%%%%
\begin{lem} \label{sym_surj} 
Let $U$ be a subset of $\Omega^\Omega$, which is not contained in $S_1$ or $S_3(\mu)$ for any infinite 
$\mu\leq |\Omega|$, and let $\kappa$ be any  cardinal such that $\aleph_0\leq \kappa\leq |\Omega|$. If there exists a 
surjective $g\in U$ and a transversal $\Gamma$ of $g$ such that $U$ contains 
every $a\in \sym(\Omega)$ with 
$\supp(a)\subseteq \Gamma$ and $|\supp(a)|< \kappa$,  then there exists  $g^*\in \genset{U}$ such that 
$c(g^*)\geq\kappa$ and $\Gamma g^*=\Omega$. 
\end{lem}
\proof 
Let $g'\in \Omega^\Omega$ be any inverse for $g$ such that $\Omega g'=\Gamma$ and   let $u'\in \Omega^\Omega$ be 
an arbitrary inverse for $u$ for all  
$u\in U\setminus\{g\}$.
We denote $\set{u'\in \Omega^\Omega}{u\in U}$ by $U'$ and we set $\Lambda: U'\to \P(\Omega)$ to be the assignment 
of transversals for $U'$ defined by $\Lambda(u')=\Omega u$.  Recall that $c(u)=d(u')$  and  $d(u)=c(u')$ 
for all $u\in U$. 

We prove that $U'$, $g'$, and $\Lambda$ satisfies the conditions of Lemma \ref{sym_inj}. Since 
$U\not\subseteq S_1$, $U\not\subseteq S_3(\mu)$, it follows that $U'\not\subseteq S_2$ and 
$U'\not\subseteq S_4(\mu)$ for all 
infinite  $\mu\leq \kappa$. Since $g$ is surjective, $g'$ is injective and by assumption $\Omega g'=\Gamma$. 
In particular, $U'$ contains every $a'=a^{-1}\in \sym(\Omega)$ where $\supp(a)\subseteq \Omega g'$ and 
$|\supp(a)|\leq\kappa$. Thus by Lemma \ref{sym_inj} there exists an injective $f^*\in \mutt{U'}{\Lambda}$ such 
that $d(f^*)\geq \kappa$ and $\Omega f^*\subseteq \Omega g'$.  By Corollary \ref{dual_cor}, $\genset{U}$ 
contains an inverse $g^*$ of $f^*$.  Therefore $c(g^*)=d(f^*)\geq \kappa$ and $\Omega f^*\subseteq\Gamma$  is a 
transversal of $g^*$, and in particular $\Gamma g^*=\Omega$. 
\qed\vspace{\baselineskip}

%%%%%%%%%%%%%%%%%%%%%%%%%

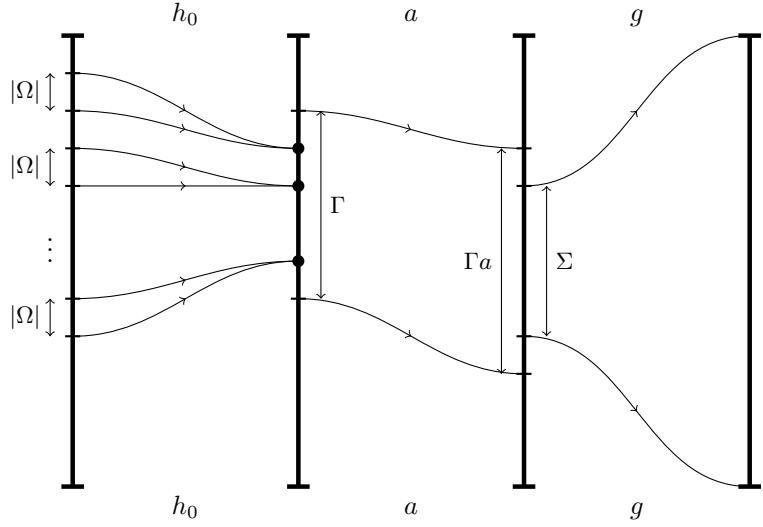
\begin{figure}
\begin{center}
 \begin{tikzpicture}[scale=1]
 \verticalline0
 \verticalline3
 \verticalline6
 \verticalline9
 \smallchip0{5.5}
 \smallchip0{5}
 \smallchip0{4.5}
 \smallchip04
 \smallchip0{2}
 \smallchip0{2.5}
 \smallchip35
 \smallchip3{2.5}
 \smallchip6{4.5}
 \smallchip64
 \smallchip62
 \smallchip6{1.5}
 \verticallinearrowleft0{5}{5.5}{$|\Omega|$}
 \verticallinearrowleft04{4.5}{$|\Omega|$}
 \verticallinearrowleft0{2}{2.5}{$|\Omega|$}
 \verticallinearrowright3{2.5}5{$\Gamma$}
 \verticallinearrowleft6{1.5}{4.5}{$\Gamma a$}
 \verticallinearrowright624{$\Sigma$}
 \pt3{4.5}
 \pt3{4}
 \pt3{3}
 \functionline0{5.5}{4.5}
 \functionline05{4.5}
 \functionline0{4}{4}
 \functionline0{4.5}{4}
 \functionline0{2.5}{3}
 \functionline023
 \functionline35{4.5}
 \functionline3{2.5}{1.5}
 \functionline646
 \functionline620
 \functionlabel0{h_0}
 \functionlabel3{a}
 \functionlabel6{g}
 \draw(-.3,3.25)node{$\vdots$};
%
% \chip03
% \verticallinearrowleft003{$\Sigma$}
 %
 %
 %
 \end{tikzpicture}
\caption{The composite $h_0ag$ in the proof of Theorem \ref{dicks}.}\label{fig3a}
\end{center}
\end{figure}

\proofref{dicks}
Let $M$ be a subsemigroup of $\Omega^\Omega$ containing $\sym(\Omega)$. We first prove that if $M$ is not 
contained in any of $S_1$, $S_2$, $S_3(\mu)$, $S_4(\mu)$, or $S_5$ where $\mu$ is any infinite cardinal not greater than 
$|\Omega|$, then $M=\Omega^\Omega$.  By 
Lemmas \ref{sym_inj} and~\ref{sym_surj}, there exist  $f, g\in M$ such that $f$ is injective, $d(f)=|\Omega|$, $g$ 
is surjective, and $c(g)=|\Omega|$.  By Theorem \ref{hhr}, it suffices to show that there exists a surjective 
$h\in M$  such that $k(h, |\Omega|)=|\Omega|$. 
Since $M \not \subseteq S_5$, there exists $h_0 \in M$ such that $k(h_0, |\Omega|)=|\Omega|$. Let 
$\Gamma=\set{\alpha\in\Omega}{|\alpha h_0^{-1}|=|\Omega|}$. Then $|\Gamma|=|\Omega|$.  Let 
$a\in \sym(\Omega)$ be any element such that $\Gamma a$ contains a transversal $\Sigma$ of $g$. So, if 
$\alpha\in \Omega$, then there exists $\beta\in\Sigma$ such that $\beta g=\alpha$ and so 
$\alpha(h_0ag)^{-1}=\alpha g^{-1}a^{-1}h_0^{-1}\supseteq \beta a^{-1}h_0^{-1}$. But $\beta a^{-1}\in \Gamma$ 
and so $|\beta a^{-1}h_0^{-1}|=|\Omega|$. Thus $|\alpha(h_0ag)^{-1}|=|\Omega|$ and, since $\alpha\in\Omega$ 
was arbitrary, it follows that $h_0ag$ is surjective and $k(h_0ag, |\Omega|)=|\Omega|$. So the proof is 
concluded by setting $h=h_0ag$; see Figure \ref{fig3a}.

If $|\Omega|$ is regular, then from the above either $M$ is contained in one of  $S_1$, $S_2$, $S_3(\mu)$, 
$S_4(\mu)$, or $S_5$; or $M=\Omega^\Omega$. It then follows that if $M$ is a maximal subsemigroup of 
$\Omega^\Omega$ containing $\sym(\Omega)$, then $M$ equals one of $S_1$, $S_2$, $S_3(\mu)$, 
$S_4(\mu)$, or $S_5$. On the other hand, if $M$ is one of the semigroups $S_1$, $S_2$, $S_3(\mu)$, 
$S_4(\mu)$, or $S_5$, then, since none of these semigroups is contained in any other, it follows that $M$ is a 
maximal subsemigroup of $\Omega^\Omega$. 

Suppose that $|\Omega|$ is singular. 
If $M$ is not contained in any of the semigroups $S_1$, $S_2$, $S_3(\mu)$, $S_4(\mu)$, or $S_5'$, then, by 
Lemma \ref{s3}, $M$ is also not contained in $S_5$ and so, from the above,  $M=\Omega^\Omega$.  Hence as 
in the case that $|\Omega|$ is regular, it follows that $M$ is a maximal subsemigroup of $\Omega^\Omega$ if 
and only if $M$ equals one of $S_1$, $S_2$, $S_3(\mu)$, $S_4(\mu)$, or $S_5'$.  \qed

%%%%%%%%%%%%%%%%%%%%%%%%%%%%
%%%%%%%%%%%%%%%%%%%%%%%%%%%%

\section{Pointwise stabilisers of finite sets -- the proof of Theorem \ref{stab_finite}}\label{proof_stab_finite}

In this section we prove Theorem \ref{stab_finite}.

%%%%%%%%%%%%%%%%%%%%%%%%%
\begin{main} 
Let $\Omega$ be any infinite set and let $\Sigma$ be a non-empty finite subset of $\Omega$.
Then the maximal subsemigroups of $\Omega^\Omega$ containing the pointwise stabiliser 
$\sym(\Omega)_{(\Sigma)}$   but not 
$\sym(\Omega)$ are:
$$F_1(\Gamma, \mu)=\set{f\in \Omega^\Omega}{d(f)\geq \mu \text{ or } \Gamma \not \subseteq \Omega f \text{ or } 
(\Gamma f^{-1}\subseteq \Gamma \text{ and } c(f)<\mu)} \cup \fin;$$
$$F_2(\Gamma, \nu)=\set{f\in \Omega^\Omega}{c(f)\geq \nu \text{ or } |\Gamma f| < |\Gamma| \text{ or } 
(\Gamma f =  \Gamma \text{ and } d(f)<\nu)} \cup \fin$$
where $\Gamma$ is a non-empty subset of $\Sigma$ and $\mu$ and $\nu$ are  infinite cardinals with 
$\mu \leq |\Omega|^+$ and either:
$|\Gamma|= 1$ and $\nu= |\Omega|^{+}$; or $|\Gamma|\geq 2$ and $\nu\leq |\Omega|^{+}$. 
\end{main}

%%%%%%%%%%%%%%%%%%%%%%%%%

Throughout this section we let  $\Sigma$ be a non-empty finite subset of $\Omega$.  We start by showing that 
the sets given in the theorem are actually semigroups. 

%%%%%%%%%%%%%%%%%%%%%%%%%

\begin{prop} 
Let $\mu$ be any infinite cardinal such that $\mu\leq |\Omega|^+$ and let $\Gamma$ be any non-empty subset 
of $\Omega$. Then $F_1(\Gamma,\mu)$ and $F_2(\Gamma, \mu)$ as defined in Theorem \ref{stab_finite} are 
subsemigroups of $\Omega^\Omega$. 
\end{prop}
\proof 
Let $f,g\in F_1(\Gamma, \mu)$. If $f\in \fin$ or $g\in \fin$, then $fg\in \fin$.
If $d(g)\geq \mu$, then Lemma \ref{tech}(ii) implies that $d(fg)\geq d(g)\geq \mu$ and so 
$fg\in F_1(\Gamma, \mu)$. 
If $\Gamma\not\subseteq \Omega g$, then $\Gamma \not\subseteq \Omega fg$ and so  
$fg\in F_1(\Gamma, \mu)$.
 Assume that $\Gamma g^{-1}\subseteq \Gamma$ and $c(g)<\mu$. If $d(f)\geq \mu$, then, by Lemma \ref{tech}
 (vi), $d(fg)\geq \mu$.  If $\Gamma \not\subseteq \Omega f$, then either $\Gamma \not\subseteq \Omega g$ or 
 $\Gamma \subseteq \Omega g$.
In the former case, $\Gamma\not\subseteq \Omega fg$, and in the latter, 
$\Gamma g^{-1}=\Gamma \not\subseteq \Omega f$ and so $\Gamma \not\subseteq \Omega fg$. In either case, 
$fg\in F_1(\Gamma,\mu)$.  If 
$\Gamma f^{-1}\subseteq\Gamma$ and $c(f)<\mu$, then $\Gamma(fg)^{-1}\subseteq\Gamma$ and 
$c(fg)\leq c(f)+c(g)<\mu$ by Lemma~\ref{tech}(iv).  Hence $F_1(\Gamma, \mu)$ is a semigroup.

Let $f, g\in F_2(\Gamma, \mu)$.  If $f\in \fin$ or $g\in \fin$, then $fg\in \fin$. If $c(f)\geq \mu$, then 
$c(fg)\geq c(f)\geq \mu$ by Lemma \ref{tech}(iv) and so $fg\in F_2(\Gamma, \mu)$. If $|\Gamma f|<|\Gamma|$, 
then $|\Gamma fg|<|\Gamma|$ and so $fg\in F_2(\Gamma,\mu)$. Hence we may assume that 
$\Gamma f=\Gamma$ and $d(f)<\mu$. If $c(g)\geq \mu$, then, by Lemma \ref{tech}(vii), $c(fg)\geq \mu$ and so
 $fg\in F_2(\Gamma, \mu)$. If $|\Gamma g|<|\Gamma|$, then $|\Gamma fg|=|\Gamma g|<|\Gamma|$ and 
 $fg\in F_2(\Gamma,\mu)$.  
 If $\Gamma g=\Gamma$ and $d(g)<\mu$, then $\Gamma fg=\Gamma g=\Gamma$ and 
 $d(fg)\leq d(f)+d(g)<\mu$, by Lemma \ref{tech}(ii), and so $fg\in F_2(\Gamma, \mu)$. 
\qed\vspace{\baselineskip}

We require the following two lemmas to prove Theorem \ref{stab_finite}.

%%%%%%%%%%%%%%%%%%%%%%%%%

\begin{figure}
\begin{center}
 \begin{tikzpicture}[scale=1]
 \verticalline0
 \verticalline3
 \verticalline6
 \verticalline9
 \verticalline{12}

 \smallchip33
 \smallchip61
 \smallchip6{0.5}
 \smallchip92
 \smallchip9{3.5}
 \smallchip{12}3
 \smallchip{12}{0.5}

 \chip{12}4
 \chip94
 \smallchip{12}{1.5}
 \smallchip{12}{3.75}
 \verticallinearrowleft006{$\Omega$}
 \verticallinearrowright{0.2}56{$\Gamma$}
 \verticallinearrowrightlow046{$\Sigma$}
 \verticallinearrowright{12.2}56{$\Gamma$}
 \verticallinearrowrightlow{12}46{$\Sigma$}
 \verticallinearrowright334{$\geq\!\kappa$}
 \verticallinearrowleft303{$\Omega g_\kappa\setminus\Sigma$}
 \verticallinearrowright6{0.5}4{$\Lambda(h_\mu)\setminus\Sigma$}
 \verticallinearrowright60{0.5}{$\kappa$}
 \verticallinearrowright902{$\geq\!\mu$}
 \verticallinearrowleft{12}{0.5}{1.5}{$\!\mu\leq$}
 \verticallinearrowright{12}{0.5}{3.75}{$\Omega f_0\setminus \Sigma$}
 %\verticallinearrowright{12}13{$\subseteq\Omega f_0$}
%
 \functionline066
 \functionline000
 \functionline366
 \functionlinedash341
 \functionlinedash330
 \functionline666 %{5.75}
 \functionline966 %{5.75}
 \functionline055
 \functionline655
 \functionline955
 \functionline053
 \functionline334
 \functionline355
 \functionlinedash344
 \functionline301
 \functionline64{3.5}
 \functionlinedash6{0.5}2
 \functionline9{3.5}3
 \functionlinedash90{0.5}
 \functionline61{2.5}
\functionline9{2.5}2
\functionlinedash95{3.75}
\functionlinedash92{1.5}
 \functionlabel0{g_{\kappa}}
 \functionlabel3{a}
 \functionlabel6{h_\mu}
 \functionlabel9{f_0}
 \colorinterval601{blue}
 \colorinterval334{blue}
\colorinterval056{red}
 \colorinterval303{red}
 \colorinterval006{red}
 \colorinterval614{red}
 \colorinterval356{red}
 \colorinterval{12}23{red}
 \colorinterval9{2.5}{3.5}{red}
 \colorinterval656{red}

 \colorinterval956{red}

 \colorinterval{12}56{red}
  \chip04
 \draw(4.5,5.25)node{{\small $\Sigma$ fixed}};
 \draw(4.5,4.75)node{{\small pointwise}};
 \end{tikzpicture}
\caption{The composite $g_{\kappa}ah_{\mu}f_0$ from the proof of Lemma \ref{finite_stab_inj}.}\label{fig3}
\end{center}
\end{figure}
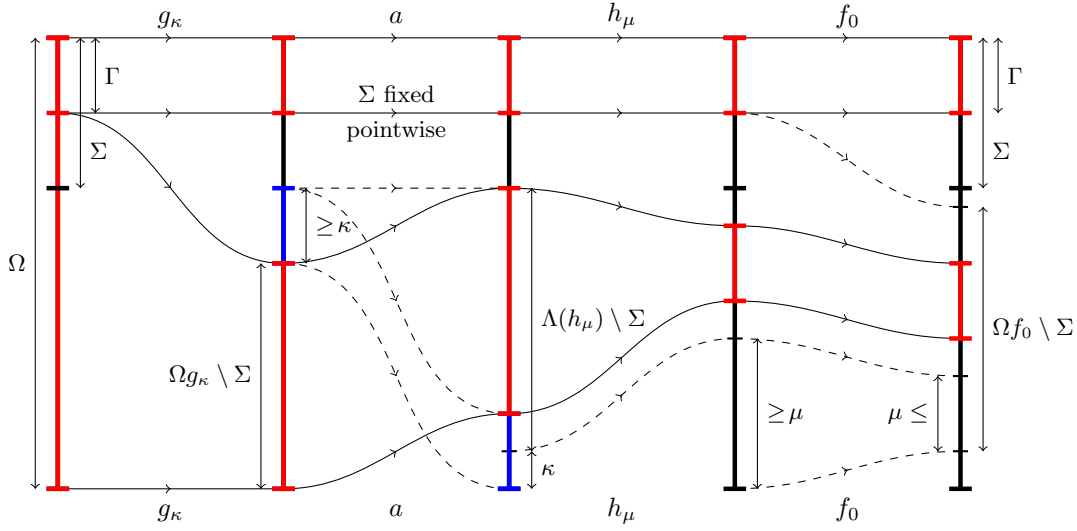

\begin{figure}
\begin{center}
 \begin{tikzpicture}[scale=1]
 \verticalline0
 \verticalline3
 \verticalline6
 \verticalline9
 \smallchip31
 \smallchip6{3.5}
 \smallchip9{4.5}
 \smallchip91
 \smallchip93
 \smallchip9{3.5}

 \verticallinearrowleft006{$\Omega$}
 \verticallinearrowright{0.2}56{$\Gamma$}
 \verticallinearrowrightlow046{$\Sigma$}
 \verticallinearrowright301{$\geq\!\nu$}
 \verticallinearrowleft314{$\Omega g_\nu\setminus\Sigma$}
 \verticallinearrowright60{3.5}{$\Lambda(h_\lambda)\setminus{\Sigma}$}
 \verticallinearrowright6{3.5}4{$\nu$}
 \verticallinearrowleft9{3.5}{4.5}{$\Gamma h_\lambda$}
  \verticallinearrowright946{$\Sigma$}
  \verticallinearrowright913{$(\Omega\setminus \Sigma)h_{\lambda}$}
 \functionline066
 \functionline001
 \functionline054
 \functionline055
 \functionline355
 \functionline366
 \functionlinedash344
 \functionline343
 \functionline310
 \functionlinedash314
 \functionlinedash303
 \functionline66{4.5}
 \functionline65{3.5}
 \functionlinedash6{3.5}3
 \functionline601
 \functionline632
 \functionlabel0{g_{\nu}}
 \functionlabel3{b}
 \functionlabel6{h_\lambda}
 \colorinterval301{blue}
 \colorinterval634{blue}
\colorinterval006{red}
\colorinterval056{red}
 \colorinterval356{red}
 \colorinterval314{red}
 \colorinterval603{red}
 \colorinterval656{red}
\colorinterval9{3.5}{4.5}{red}
\colorinterval912{red}
  \chip04
   \chip94
 \draw(4.5,5.25)node{{\small $\Sigma$ fixed}};
 \draw(4.5,4.75)node{{\small pointwise}};
 \end{tikzpicture}
\caption{The composite $g_{\nu}bh_{\lambda}$ in the proof of Lemma \ref{finite_stab_inj}.}\label{fig4}
\end{center}
\end{figure}
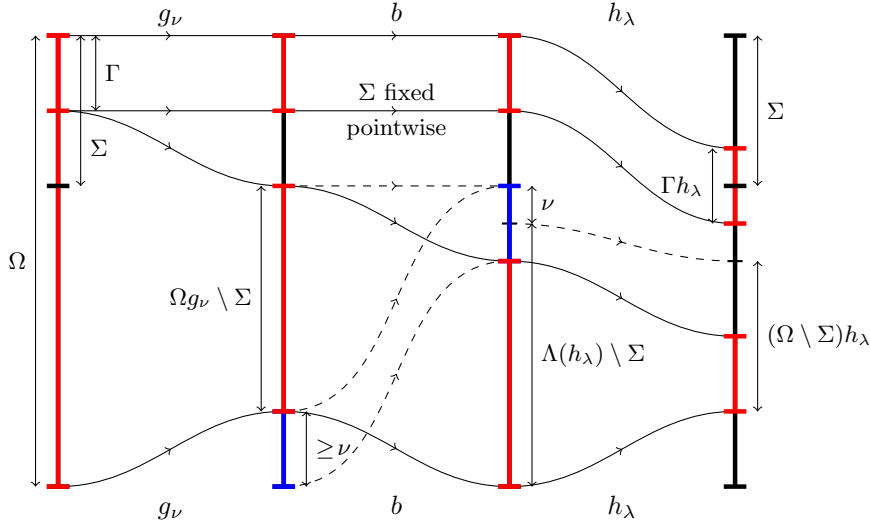

\begin{figure}
\begin{center}
 \begin{tikzpicture}[scale=1]
 \verticalline0
 \verticalline3
 \verticalline6
 \verticalline9
 \smallchip31
 \smallchip32
 \smallchip63
 \smallchip95
 \smallchip91
 \chip04
 \chip34
  \chip64
    \chip94
 \verticallinearrowleft006{$\Omega$}
 \verticallinearrowright046{$\Sigma$}
 \verticallinearrowleft316{$\Omega s^n$}
 \verticallinearrowleft915{$\Omega f_2$}
 \functionline066
 \functionline001
 \functionlinedash3{0.5}3
 \functionline344
 \functionline366
 \functionline665
 \functionlinedash63{4.5}
 \functionline601

 \draw(3,.5)node[left]{{\small $j$}};
 \draw(6,3)node[below right]{{\small $i$}};
 \draw(9,4.5)node[right]{{\small $if_2$}};
 \functionlabel0{s^n}
 \functionlabel3{p}
 \functionlabel6{f_2}
%
 %\colorinterval046{red}
 %\colorinterval346{red}
 %\colorinterval646{red}
 %\colorinterval946{red}
%
 \pt3{0.5}
 \pt63
 \pt9{4.5}
 \draw(4.5,5.25)node{{\small $\Sigma$ fixed}};
 \draw(4.5,4.75)node{{\small pointwise}};
 \end{tikzpicture}
\caption{The composite $s^npf_2$ in the proof of Lemma \ref{finite_stab_inj}.}\label{fig5}
\end{center}
\end{figure}

\begin{lem}\label{finite_stab_inj} Let $\Sigma$ be a finite subset of $\Omega$ and let $U$ be a subset of $\Omega^\Omega$ 
containing $\sym(\Omega)_{(\Sigma)}$ but which is not 
contained in $S_2$ or in $F_2( \Gamma , \mu)$ for any non-empty subset $\Gamma$ of $\Sigma$ and 
any infinite cardinal $\mu \leq |\Omega|^+$.
If $\Lambda$ is an assignment of transversals  (as defined 
in Definition \ref{assign}) for $U$ such that $\Gamma\subseteq \Lambda(u)$ for all $u\in U\setminus F_2(\Gamma,\mu)$, then there exists an injective $f\in \mutt{U}{\Lambda}$ such that $\Omega f \cap \Sigma =\emptyset$.
\end{lem}
\proof
 Since $\Sigma$ is finite (and $1_{\Omega} \in U$), it suffices to show that for every injective $f_0\in \mutt{U}{\Lambda}$ 
 with $\Omega f_0 \cap \Sigma\not = \emptyset$ there exists an injective $f_{1}\in  \mutt{U}{\Lambda}$ 
 with $\Omega f_{1} \cap \Sigma \subsetneq\Omega f_0 \cap \Sigma$.  We will denote $\Omega f_0 \cap \Sigma$ by 
 $\Gamma$.  We start by showing that there exists an injective $f_2\in \mutt{U}{\Lambda}$ such that 
 $\Omega f_2\cap \Sigma\subseteq \Gamma$ and $\Gamma f_2\not=\Gamma$.

If $\Gamma f_0\not=\Gamma$, then let $f_2=f_0$. Hence we may assume that $\Gamma f_0=\Gamma$. Since $U$ is not
contained in $S_2$, there exists an injective $s\in U\setminus S_2$ such that $d(s)>0$. If $\Gamma s\not=\Gamma$, then  
we set $f_2=sf_0$. Thus the final case to consider is when $\Gamma s=\Gamma$. 

For every infinite cardinal $\mu$ with $\mu \leq |\Omega|^+$, let $h_\mu$ be an element of $U\setminus F_2(\Gamma, \mu)$. 
Then the following hold:
$$c(h_\mu)<\mu,\quad \Gamma \subseteq \Lambda(h_{\mu}),$$ 
and either  
$$\Gamma h_{\mu} \not = \Gamma\text{ or }d(h_{\mu})\geq \mu.$$ 

Note that $d(h_{|\Omega|^+})\leq |\Omega|<|\Omega|^+$ and so $\Gamma h_{|\Omega|^+}\not = \Gamma$. Thus we may let $\lambda$ be the least infinite cardinal such that $\Gamma h_{\lambda} \not =\Gamma$. We will  show by transfinite induction that for every cardinal $\mu$ strictly less than $\lambda$ 
\begin{equation}\label{pointwise_induction}
\text{ there exists an injective $g_{\mu} \in \mutt{U}{\Lambda}$ with $\Omega g_{\mu}\cap \Sigma = \Gamma$, $\Gamma g_{\mu}=\Gamma$ and $d(g_{\mu})\geq \mu$.}  
\end{equation}
For any finite $\mu$, we may let $g_{\mu}:=s^{\mu}f_0$. So let $\mu<\lambda$ be an infinite cardinal and assume (\ref{pointwise_induction}) holds for all cardinals strictly less than $\mu$. 
%To avoid double subscripts, let $\kappa:=c(h_{\mu})$. 
By the inductive assumption there exists an injective $g_{c(h_{\mu})}\in \mutt{U}{\Lambda}$ with $\Omega g_{c(h_{\mu})}\cap \Sigma = \Gamma$,
 $\Gamma g_{c(h_{\mu})}=\Gamma$, and
 $d(g_{c(h_{\mu})})\geq c(h_{\mu})$. 
 Hence there exists a bijection $a\in \sym(\Omega)_{(\Sigma)}$ such that  $(\Omega g_{c(h_{\mu})})a\subseteq \Lambda(h_{\mu})$. 
We define $g_{\mu}:=g_{c(h_{\mu})}a h_{\mu} f_0$; see Figure \ref{fig3}.
Then by construction $g_{\mu}\in \mutt{U}{\Lambda}$,
$g_{\mu}$ is  injective, and since $f_0$ is injective $d(g_{\mu})\geq  d(h_{\mu})\geq \mu$. 
Also $\Gamma h_{\mu}=\Gamma$, since $\mu<\lambda$, and $g_{c(h_{\mu})},f,  a$ stabilise $\Gamma$ setwise, and hence $\Gamma g_{\mu} = \Gamma$. Finally,  
$$\Gamma = \Gamma g_{\mu}=\Gamma g_{\mu} \cap \Sigma \subseteq \Omega g_{\mu}\cap \Sigma \subseteq \Omega f_{0}\cap \Sigma = \Gamma.$$ 
 Hence (\ref{pointwise_induction}) holds for all $\mu < \lambda$.

For the same reason as above we will denote $c(h_{\lambda})$ by $\nu$. Since $\nu<\lambda$, there exists an injective $g_{\nu} \in \mutt{U}{\Lambda}$ with $\Omega g_{\nu}\cap \Sigma = \Gamma$, $\Gamma g_{\nu}=\Gamma$ and $d(g_{\nu})\geq \nu$. Let $b\in \sym(\Omega)_{(\Sigma)}$ be such that  $(\Omega g_{\nu})b\subseteq \Lambda(h_{\lambda})$; see Figure \ref{fig4}. Then $g_{\nu}b h_{\lambda}$ is injective, and since $g_{\nu}$ and $b$ stabilise $\Gamma$ setwise, but $h_{\lambda}$ does not, $\Gamma g_{\nu}b h_{\lambda}\not=\Gamma$. Thus we let $f_2=g_{\nu} a h_{\lambda}$, which completes this part of the proof.\vspace{\baselineskip}

We will use the function $f_2$ to prove that there exists an injective $f_{1}\in  \mutt{U}{\Lambda}$ 
 with $\Omega f_{1} \cap \Sigma \subsetneq\Gamma$.
 If $\Omega f_2\cap \Sigma\not=\Gamma$, then setting $f_1=f_2$ concludes the proof. 
 Hence we only have to consider the case when $\Omega f_2\cap \Sigma=\Gamma$. Since $\Gamma f_2\not=\Gamma$, it follows, in this case, that $\Gamma f_2^{-1}\not\subseteq \Gamma=\Omega f_2\cap \Sigma$.  Thus there are two cases to consider: $\Gamma f_2^{-1}\not\subseteq \Omega f_2$ or $\Gamma f_2^{-1}\not\subseteq \Sigma$.  If $\Gamma f_2^{-1}\not\subseteq \Omega f_2$, then $\Omega f_2^2\cap \Sigma\subsetneq \Gamma$, and we set $f_1=f_2^2$. 
 If $\Gamma f_2^{-1}\not\subseteq \Sigma$, then 
 there exists  $i \in \Gamma f_2^{-1} \setminus \Sigma$.  
Since $U$ is not contained in $S_2$, there exists $s\in U\setminus S_2$ such that $s$ is injective and $d(s)>0$. 
It follows from Lemma \ref{tech}(iii) that $d(s^n)>|\Sigma|$ for some $n\in\nat$. Hence there exists 
$j\in \Omega\setminus \Sigma$ such that $ j\not\in \Omega s^n$ and there is $p\in \sym(\Omega)_{(\Sigma)}$ such that 
$(j)p=i$.
In this case, we set $f_1:=s^{n}pf_2$; see Figure \ref{fig5}.  
Then since $j\not\in \Omega s^{n}$ and $pf_2$ is injective, it follows that 
$if_2=(j)pf_2\not\in \Omega s^npf_2=\Omega f_1$. But $if_2\in\Gamma$, and so 
$\Omega f_1\cap \Sigma\subseteq \Gamma\setminus \{if_2\}\subsetneq  \Gamma$, as required. 
 \qed

%%%%%%%%%%%%%%%%%%%%%%%%%

\begin{lem}\label{finite_stab_surj}
 Let $\Sigma$ be a finite subset of $\Omega$ and let $U$ be a subset of $\Omega^\Omega$ 
containing $\sym(\Omega)_{(\Sigma)}$ but which is not 
contained in $S_1$ or in $F_1( \Gamma , \mu)$ for any non-empty subset $\Gamma$ of $\Sigma$ and 
any infinite cardinal $\mu \leq |\Omega|^+$.
Then there exists a surjective $g\in \genset{U}$ such that $(\Omega\setminus \Sigma) g =\Omega$.
\end{lem}
\proof
If $u\in U$ is arbitrary, then we denote an arbitrary inverse for $u$ by $u'$. 
We denote $\set{u'\in \Omega^\Omega}{u\in U}$ by $U'$ and we set $\Lambda: U'\to \P(\Omega)$ to be the 
assignment of transversals for $U'$ defined by $\Lambda(u')=\Omega u$.  

Since $\sym(\Omega)_{(\Sigma)}\subseteq U\not\subseteq S_1\cup  F_1(\Gamma, \mu)$, it follows that 
$\sym(\Omega)_{(\Sigma)}\subseteq U'\not\subseteq S_2\cup F_2(\Gamma, \mu)$ for all non-empty subsets 
$\Gamma$ of $\Sigma$ and for all infinite $\mu\leq |\Omega|^+$. If $u'\not\in F_2(\Gamma, \mu)$ for some 
$u\in U$, then $u\not\in F_1(\Gamma, \mu)$ and so $\Gamma\subseteq \Omega u=\Lambda(u')$. Thus by 
Lemma \ref{finite_stab_inj} there exists an injective $f\in \mutt{U'}{\Lambda}$ such that 
$\Omega f \cap \Sigma=\emptyset$. Then, by Corollary \ref{dual_cor}, $f$ has an inverse $g\in \genset{U}$. 
Then $g$ is surjective and 
$\Omega f$ is a transversal of $g$. In particular, $(\Omega\setminus\Sigma) g=\Omega g=\Omega$, as required.
\qed\vspace{\baselineskip}
%%%%%%%%%%%%%%%%%%%%%%%%%

\proofref{stab_finite}
It is straightforward to verify that none of the semigroups listed in the statement of Theorem \ref{stab_finite} are 
contained in any of the others from that list.  Moreover, none of these semigroups are contained in any of the 
semigroups from Theorem \ref{dicks}.

Let $M$ be a subsemigroup of $\Omega^\Omega$ containing $\sym(\Omega)_{(\Sigma)}$ that is not contained 
in  any of the semigroups  in Theorems \ref{dicks} or \ref{stab_finite}. 
We will prove that $\sym(\Omega)$ is a subsemigroup of $M$ and so Theorem \ref{dicks} implies that 
$M=\Omega^\Omega$.  

Let $\Gamma$ be a finite subset of $\Omega$ and let $\mu$ be an infinite cardinal such that 
$\mu\leq |\Omega|^+$.  If $u\in \Omega^\Omega$ but $u\not\in F_2( \Gamma , \mu)$, then, in particular, $u$ is 
injective on $\Gamma$ and so there exists a transversal of $u$ containing $\Gamma$. In particular, there is an 
assignment  of transversals $\Lambda$ for $M$ such that $\Gamma\subseteq \Lambda(u)$ for all 
$u\in M\setminus F_2(\Gamma,\mu)$.
Hence by Lemma \ref{finite_stab_inj}, there exists an injective $f\in M$ such that 
$\Omega f\cap \Sigma=\emptyset.$
Since $M$ contains all permutations with support contained in $\Omega \setminus \Sigma$, it contains all 
permutations with support contained in $\Omega f$. Thus by Lemma \ref{sym_inj} there exists an injective 
$f^*\in M$ with $d(f^*)=|\Omega|$ and $\Omega f^*\subseteq \Omega f\subseteq \Omega \setminus \Sigma$. 

By Lemma \ref{finite_stab_surj}, there exists a surjective $g\in M$ with a transversal 
$\Gamma\subseteq \Omega\setminus\Sigma$. Clearly $M$ contains every permutation with support contained 
in $\Gamma$. Hence by Lemma \ref{sym_surj} there exists $g^*\in M$ such that $c(g^*)=|\Omega|$ and $\Gamma g^*=\Omega$. 

Since $\Omega f^*$ and $\Gamma$ are moieties of $\Omega\setminus \Sigma$, every bijection from $\Omega f^*$ to $\Gamma$ is a restriction of some element of $\sym(\Omega)_{(\Sigma)}$.
So, if $a\in \sym(\Omega)$ is arbitrary, then, since $f^*$ and $g^*|_{\Gamma}$ are injective,  there  exists $b\in \sym(\Omega)_{(\Sigma)}$ such that $a=f^*bg^*$. 
Therefore $\sym(\Omega)$ is a subsemigroup of $M$ and so, by Theorem \ref{dicks}, $M=\Omega^\Omega$. 

We have shown that if $M$ is a subsemigroup of $\Omega^\Omega$ that contains  $\sym(\Omega)_{(\Sigma)}$, then either $M$ is contained in one of the semigroups from Theorem \ref{dicks}, one of the semigroups  $F_1(\Gamma, \mu)$ or $F_2(\Gamma, \nu)$ from the statement of the theorem, or $M=\Omega^\Omega$.  It follows that if $M$ is maximal, then $M$ is one of these semigroups.
On the other hand,  if  $M$ is one of $F_1(\Gamma, \mu)$ or $F_2(\Gamma, \nu)$ then, since none of these semigroups is contained in any of the others or any of the semigroups in Theorem \ref{dicks}, it follows that $M$ is a maximal subsemigroup of $\Omega^\Omega$.  
\qed

\section{The stabiliser of an ultrafilter - the proof of Theorem \ref{ultra}}\label{proof_ultra}

%%%%%%%%%%%%%%%%%%%%%%%%%

In this section we give the proof of Theorem \ref{ultra}. 

\begin{main}
Let $\Omega$ be any infinite set, let $\F$ be a non-principal ultrafilter on $\Omega$, and let $\kappa(\geq \aleph_0)$ be 
the least 
cardinality of an element of $\F$.  Then the maximal subsemigroups of $\Omega^\Omega$ containing 
$\sym(\Omega)_{\{\F\}}$ but not $\sym(\Omega)$ are:
\begin{align*}
U_1(\F, \mu)&=\set{f\in \Omega^\Omega}{(d(f)\geq\mu)\text{ or }(\Omega f\not\in \F) \text{ or }(c(f)< \mu \text{ and } 
(\forall \Sigma \not\in \F)(\Sigma f\not\in \F))}\cup \fin;\\
U_2(\F, \mu)&=\set{f\in \Omega^\Omega}{(c(f)\geq\mu)\text{ or }(\forall \Sigma \in \F)( c(f|_\Sigma)>0)\text{ or }(d(f)< \mu
\text{ and }(\forall \Sigma \in \F)(\Sigma f\in \F))}\cup \fin.
\end{align*}
for cardinals $\mu$ such that $\kappa< \mu\leq |\Omega|^+$. 
\end{main}

Throughout this section we let $\F$ be an arbitrary non-principal 
ultrafilter on $\Omega$ and let $\kappa$ be the least cardinality of a set belonging to $\F$. Since $\F$ is non-principal, it
follows that $\kappa\geq \aleph_0$. 
A subset $S$ of $\sym(\Omega)$ is \emph{transitive on moieties} of $\Omega$ if for all moieties $\Sigma, \Gamma$ of $\Omega$ there exists $f\in S$ such that $\Sigma f=\Gamma$. 
Recall that $\sym(\Omega)_{\{\F\}}$ is transitive on moieties in $\F$ 
and hence also moieties not in $\F$. Moreover, if $\Gamma,\Sigma\in \F$ such that 
$|\Omega\setminus \Gamma|\geq |\Omega\setminus \Sigma|$, then there exists $a\in \sym(\Omega)_{\{\F\}}$ 
such that $\Gamma a\subseteq \Sigma$. 

The following lemma and its proof are similar to Lemma \ref{finite_stab_inj}.
We use the following observation in the statement and proof of the next lemma.  If $f\in\Omega^\Omega$ but 
$f\not\in U_2(\F, \mu)$, then there exists $\Sigma\in \F$ such that $c(f|_{\Sigma})=0$, in other words $f$ is injective on 
$\Sigma$. It follows that $\Sigma$ is contained in a transversal $\Lambda(f)$ for $f$ and so $\Lambda(f)\in \F$. We have 
shown that every element of $\Omega^\Omega$ which does not belong to $U_2(\F, \mu)$ has a transversal in $\F$. 

\begin{figure}
\begin{center}
 \begin{tikzpicture}[scale=1]
 \verticalline0
 \verticalline3
 \verticalline6
 \verticalline9
 \smallchip34
 \smallchip62
 \smallchip6{1.5}
 \smallchip9{4.5}
 \verticallinearrowleft006{$\Omega$}
 \verticallinearrowlefthigh346{$c(h_\mu)\!\leq$}
 \verticallinearrowright60{1.5}{$c(h_\mu)$}
 \verticallinearrowright6{1.5}6{$\Lambda(h_\mu)$}
 \verticallinearrowright9{4.5}6{$\geq\!\mu$}
  \verticallinearrowleft304{$\Omega g_{c(h_\mu)}$}
  \verticallinearrowright90{4.5}{$\Omega h_{\mu}$}
 \functionline000
 \functionline064
 \functionlinedash362
 \functionlinedash340
 \functionline346
 \functionline302
 \functionline66{4.5}
 \functionlinedash6{1.5}0
 \functionline62{0.5}
 \functionlabel0{g_{c(h_\mu)}}
 \functionlabel3{a}
 \functionlabel6{h_\mu}
\colorinterval346{blue}
\colorinterval602{blue}
\colorinterval{3}{0}{4}{red}
\colorinterval{6}{2}{6}{red}
\colorinterval006{red}
\colorinterval9{0.5}{4.5}{red}
 \end{tikzpicture}
\caption{The composite $g_{c(h_{\mu})}ah_{\mu}$ in the proof of Lemma \ref{ultra_inj}.}\label{fig8}
\end{center}
\end{figure}
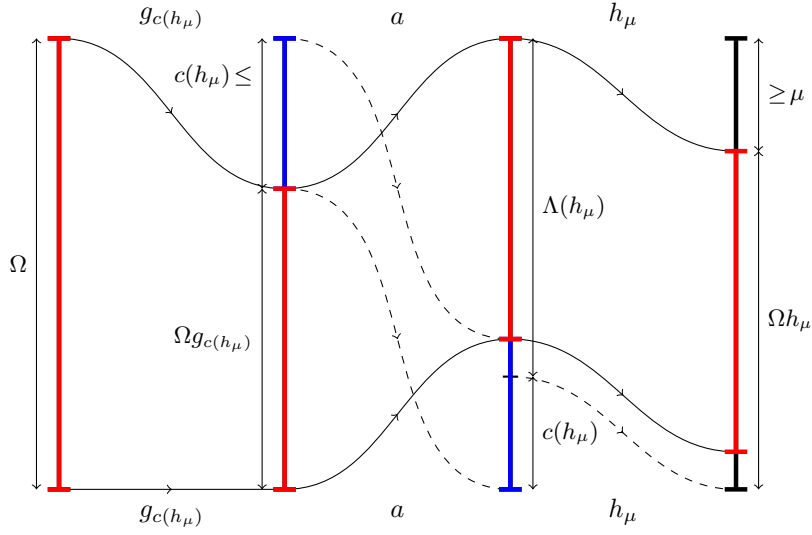

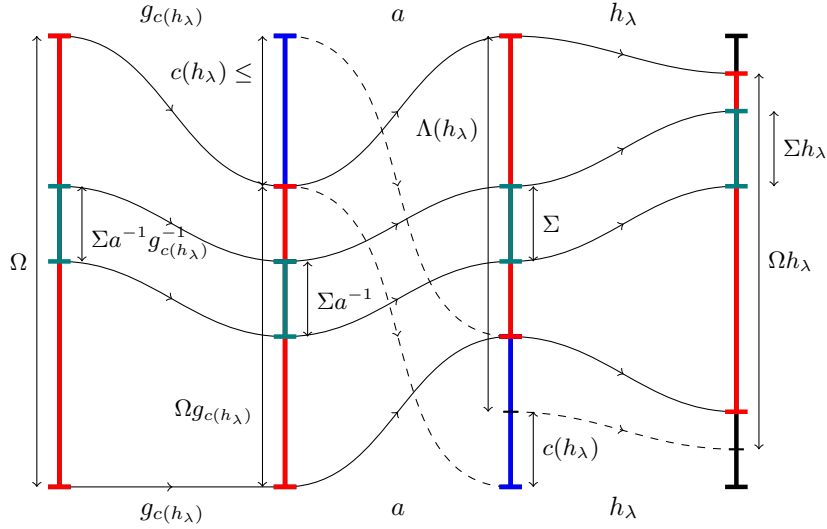
\begin{figure}
\begin{center}
 \begin{tikzpicture}[scale=1]
 \verticalline0
 \verticalline3
 \verticalline6
 \verticalline9
 \chip33
 \smallchip34

 \chip62

 \smallchip9{0.5}
  \smallchip9{5.5}
 \verticallinearrowleft006{$\Omega$}
 \verticallinearrowrightlow034{$\Sigma a^{-1}g^{-1}_{c(h_{\lambda})}$}
  \verticallinearrowright323{$\Sigma a^{-1}$}
  \verticallinearrowright634{$\Sigma$}
  \verticallinearrowlefthigh346{$c(h_{\lambda})\leq$}
    \verticallinearrowright601{$c(h_{\lambda})$}
  \verticallinearrowlefthigh616{$\Lambda(h_{\lambda})$}
  \verticallinearrowright{9.2}45{$\Sigma h_{\lambda}$}
  \verticallinearrowright9{0.5}{5.5}{$\Omega h_{\lambda}$}
  \verticallinearrowleftlow304{$\Omega g_{c(h_{\lambda})}$}
 \functionline000
 \functionline064
 \functionline346
% \functionline600
 \functionlinedash362
  \functionlinedash340
  \functionline043
   \functionline032
  \functionline302
  \functionline323
  \functionline334 
 \functionline66{5.5}
  \functionlinedash61{.5}
  \functionline634
  \functionline645
  \functionline621
 \functionlabel0{g_{c(h_{\lambda})}}
 \functionlabel3{a}
 \functionlabel6{h_{\lambda}}
\colorinterval346{blue}
\colorinterval602{blue}
\colorinterval006{red}
\colorinterval91{5.5}{red}
\colorinterval{6}{2}{6}{red}
\colorinterval{3}{0}{4}{red}
 \colorinterval{0}{3}{4}{teal}
  \colorinterval{3}{2}{3}{teal}
  \colorinterval{6}{3}{4}{teal}
  \colorinterval{9}{4}{5}{teal}
  \smallchip61

 \end{tikzpicture}

\caption{The composite $g_{c(h_{\lambda})}ah_{\lambda}$ in the proof of Lemma \ref{ultra_inj}.}\label{fig9}
\end{center}
\end{figure}

\begin{figure}
\begin{center}
 \begin{tikzpicture}[scale=1]
 \verticalline0
 \verticalline3
 \verticalline6
 \verticalline9
 \chip33
 \smallchip34
 \smallchip63
 \smallchip6{5}
 \chip93
 \smallchip9{2.5}
 \smallchip92
 \smallchip9{1.5}
 \smallchip9{.75}
 \verticallinearrowleft006{$\Omega$}
 \verticallinearrowleft304{$\Omega f_0$}
 \verticallinearrowright303{$\Sigma$}
 \verticallinearrowrighthigh903{$\Sigma$}
 \verticallinearrowright346{$\geq\!\kappa$}
 \verticallinearrowleft603{$\Omega f_0c$}
 \verticallinearrowleft656{$\Omega f_0c$}
 \verticallinearrowleft6{4.75-1}{5.5-1}{$\Gamma f_1^{-1}$}
 \verticallinearrowright603{$\Sigma$}
 \verticallinearrowleft90{.75}{$\Sigma f_1$}
 \verticallinearrowlefthigh92{2.5}{$\Gamma$}
 \verticallinearrowright{9.2}0{1.5}{$\Omega f_0cf_1$}
 \functionline000
 \functionline064
 \functionline333
 \functionline335
 \functionline346
 \functionline300

 \functionlineleft66{1.5}
 \functionlineleft65{0.75}
 \functionlinealeftdash6{4.5}{2.5}
  \functionlinearightdash6{4.5}{2.5}

 \functionlinealeftdash6{3.75}2
  \functionlinearightdash6{3.75}2
 \functionlineright66{1.5}
 \functionlineright65{0.75}
 \functionline63{.75}
 \functionline600
 \functionlinedash365
  \functionlinedash343
 \functionlabel0{f_0}
 \functionlabel3{c}
 \functionlabel6{f_1}
 \draw(4.5,1)node{{\small fixes $\Sigma$}};
 \draw(4.5,.6)node{{\small pointwise}};
\colorinterval346{blue}
\colorinterval635{blue}
\colorinterval006{red}
 \colorinterval304{red}
 \colorinterval603{red}
 \colorinterval656{red}
 \colorinterval303{red}
\colorinterval90{0.75}{red}
\colorinterval9{0.75}{1.5}{red}

 \end{tikzpicture}

\caption{The composite $f_0cf_1$ in the proof of Lemma \ref{ultra_inj}.}\label{lemma6_2}
\end{center}
\end{figure}

\begin{lem}\label{ultra_inj}
Let $U$ be a subset of $\Omega^\Omega$ containing the stabiliser $\sym(\Omega)_{\{\F\}}$ of $\F$ but which is not 
contained in $U_2(\F, \mu)$, $S_2$, or $S_4(\nu)$  for any cardinals $\mu,\nu$ such that 
$\aleph_0\leq\nu\leq \kappa< \mu\leq |\Omega|^+$, and let $\Lambda$ be an 
assignment of transversals  (as defined 
in Definition \ref{assign}) for $U$ such that $\Lambda(u)\in \F$ for all $u\in U\setminus U_2(\F, \mu)$. Then there exists an 
injective $f\in \mutt{U}{\Lambda}$ such that $\Omega f\not\in \F$. 
\end{lem}
\proof
If $\Sigma\subseteq \Omega$ such that $|\Sigma|<\kappa$, then $\Sigma\not\in \F$ and so every $a\in \sym(\Omega)$ 
such that $|\supp(a)|<\kappa$ belongs to $\sym(\Omega)_{\{\F\}}$ and hence to $U$. 
Thus by Lemma \ref{sym_inj} there exists an injective $f_0\in \mutt{U}{\Lambda}$ such that $d(f_0)\geq \kappa$. 
We start by showing that there exists an injective $f_1\in \mutt{U}{\Lambda}$ and $\Sigma\in \F$ such that 
$\Sigma f_1\not\in \F$.

If there exists $\Sigma\in \F$ such that $\Sigma f_0\not\in \F$, then $f_1:=f_0$ is the required function. 
Hence we may assume that $\Sigma f_0\in \F$ for all $\Sigma \in \F$. 
For every cardinal $\mu$ such that $\kappa< \mu \leq |\Omega|^+$, let $h_\mu$ be an element of 
$U\setminus U_2(\F, \mu)$. 
Then the following hold:  $c(h_\mu)< \mu$, $\Lambda(h_{\mu})\in \F$, and either 
 $d(h_{\mu})\geq\mu$ or $\Sigma h_{\mu}\not\in \F$ for some $\Sigma\in \F$. 
Note that $d(h_{|\Omega|^+})\leq |\Omega|<|\Omega|^+$ and so there exists $\Sigma\in \F$ such that 
$\Sigma h_{|\Omega|^+}\not\in \F$.
Thus we may define
$$\lambda=\min\set{\mu}{\kappa< \mu\leq |\Omega|^+\text{ and }(\exists \Sigma\in \F)(\Sigma h_{\mu}\not\in\F)}.$$ 

 We will  show, by transfinite 
induction, that for every cardinal $\mu$ strictly less than $\lambda$:
\begin{equation}\label{ultra_induction}
\text{ there exists an injective $g_{\mu} \in \mutt{U}{\Lambda}$ such that   $d(g_{\mu})\geq \mu$ and 
$\Sigma g_{\mu}\in \F$ for all $\Sigma \in \F$}.
\end{equation}
By assumption, $f_0$ satisfies \eqref{ultra_induction}, for all $\mu\leq \kappa$. 
 So let $\mu$ be any cardinal such that $\kappa<\mu<\lambda$  and assume that (\ref{ultra_induction}) holds for all 
 cardinals strictly less than $\mu$.

By the inductive assumption there exists an injective $g_{c(h_{\mu})}\in \mutt{U}{\Lambda}$ such that 
 $d(g_{c(h_{\mu})})\geq c(h_{\mu})$ and $\Sigma g_{c(h_{\mu})}\in \F$ for all $\Sigma \in \F$. In particular, 
 $\Omega g_{c(h_{\mu})}\in \F$ and so by the comments preceding the lemma there exists $a\in \sym(\Omega)_{\{\F\}}$ 
 such that 
  $\Omega g_{c(h_{\mu})}a\subseteq \Lambda(h_{\mu})$; see Figure \ref{fig8}. 
We define $g_{\mu}:=g_{c(h_{\mu})}a h_{\mu}$.
Then by construction $g_{\mu}\in \mutt{U}{\Lambda}$,
$g_{\mu}$ is  injective, and $d(g_{\mu})\geq  d(h_{\mu})\geq \mu$. 
Also $\Sigma g_{\mu}\in \F$ for all $\Sigma \in \F$ since this property holds for $g_{c(h_{\mu})}$, $a$, and $h_{\mu}$ 
(since $\mu<\lambda$). 
 Hence (\ref{ultra_induction}) holds for all $\mu < \lambda$. 

Since $c(h_{\lambda})<\lambda$, there exists an injective $g_{c(h_{\lambda})} \in \mutt{U}{\Lambda}$ such that 
$d(g_{c(h_{\lambda})})\geq c(h_{\lambda})$ and $\Sigma g_{c(h_{\lambda})}\in \F$ for all $\Sigma\in \F$. Then as above 
there exists $b\in \sym(\Omega)_{\{\F\}}$  such that  $(\Omega g_{c(h_{\lambda})})b\subseteq \Lambda(h_{\lambda})$; 
see Figure \ref{fig9}. 
Then  $g_{c(h_{\lambda})}b h_{\lambda}\in \mutt{U}{\Lambda}$ is injective.
By the definition of $\lambda$ there exists $\Sigma\in \F$ such that $\Sigma h_{\lambda}\not\in \F$. Hence 
$\Sigma a^{-1}\in \F$ and so $\Omega \setminus \Sigma a^{-1}\not\in \F$. It follows that
$(\Omega\setminus \Sigma a^{-1})g^{-1}\not\in \F$ and, since $g$ is injective, we have that $\Sigma a^{-1}g^{-1}\in \F$.  
Thus if we let $f_1=g_{c(h_{\lambda})} a h_{\lambda}$ and $\Sigma'=\Sigma a^{-1}g^{-1} \in \F$, then 
$\Sigma' f_1=\Sigma h_{\lambda}\not\in \F$,  which completes this 
part of the proof.  

If $\Omega f_1\not\in \F$, then $f_1$ satisfies the conclusion of the lemma. If $\Omega f_1\in \F$, then there exists 
$\Gamma\subseteq \Omega f_1\setminus \Sigma' f_{1}$ such that $\Gamma\in \F$ and $|\Gamma|=\kappa$ (the least
cardinality of a set in $\F$). Hence $\Gamma f_1^{-1}\cap \Sigma'=\emptyset$ and $|\Gamma f_1^{-1}|=\kappa$, since 
$f_1$ is injective.  Thus there exists $c\in \sym(\Omega)_{(\Sigma)}\leq  \sym(\Omega)_{\{\F\}}$ such that 
$\Gamma f_1^{-1}\subseteq \Omega\setminus\Omega f_0 c$; see Figure \ref{lemma6_2}.  The required function is then $f=f_0cf_1$ since 
$\Omega f \subseteq \Omega\setminus\Gamma\not\in \F$ and so $\Omega f\not\in \F$. Finally, $f\in\mutt{U}{\Lambda}$ since $f_0, a, f_1\in \mutt{U}{\Lambda}$. 
\qed

%%%%%%%%%%%%%%%%%%%%%%%%%

\begin{lem}\label{ultra_surj}
Let $U$ be a subset of $\Omega^\Omega$ containing the stabiliser $\sym(\Omega)_{\{\F\}}$ of $\F$ but which is not 
contained in $U_1(\F, \mu)$, $S_1$, or $S_3(\nu)$  for any cardinals $\mu,\nu$ such that 
$\aleph_0\leq\nu\leq \kappa< \mu\leq |\Omega|^+$. Then there exists a surjective $g\in \genset{U}$ with a transversal 
$\Lambda(g)$ which does not 
belong to $\F$. 
\end{lem}
\proof
If $u\in U$ is arbitrary, then we denote an arbitrary inverse for $u$ by $u'$. 
We denote $\set{u'\in \Omega^\Omega}{u\in U}$ by $U'$ and we set $\Lambda: U'\to \P(\Omega)$ to be the 
assignment of transversals for $U'$ defined by $\Lambda(u')=\Omega u$.  

Since $\sym(\Omega)_{\{\F\}}\subseteq U\not\subseteq U_2(\F, \mu)\cup S_1\cup S_3(\nu)$, it follows that 
$\sym(\Omega)_{\{\F\}}\subseteq U'\not\subseteq U_2(\F, \mu)\cup S_2\cup S_4(\nu)$ for any cardinals $\mu,\nu$ such 
that  $\aleph_0\leq\nu\leq \kappa< \mu\leq |\Omega|^+$. If $u'\not\in U_2(\F, \mu)$ for some 
$u\in U$, then $u\not\in U_1(\F, \mu)$ and so $\Lambda(u')=\Omega u\in \F$. Thus by 
Lemma \ref{ultra_inj} there exists an injective $f\in \mutt{U'}{\Lambda}$ such that 
$\Omega f \not\in \F$. Then, by Corollary \ref{dual_cor}, $f$ has an inverse $g\in \genset{U}$. 
Then $g$ is surjective and 
$\Omega f\not\in \F$ is a transversal of $g$, as required. \qed\vspace{\baselineskip}

%%%%%%%%%%%%%%%%%%%%%%%%%

\proofref{ultra}  It is easy to check that $U_1(\F, \mu)$ and $U_2(\F, \mu)$ are semigroups, and that neither  is contained in 
the other, nor in any of the semigroups listed in Theorem \ref{dicks}.  
Let $M$ be any subsemigroup of~$\Omega^\Omega$ containing $\sym(\Omega)_{\{\F\}}$.
As in the proof of Theorem \ref{stab_finite}, it suffices to prove that if $M$ is not contained in any of the semigroups from 
Theorems \ref{dicks} or \ref{ultra}, then $M=\Omega^\Omega$. 

By Lemmas \ref{ultra_inj} and \ref{ultra_surj}, there exist $f,g\in M$ such that $f$ is injective, $\Omega f\not\in \F$, $g$ is 
surjective and $g$ has a transversal $\Lambda(g)\not\in \F$. Since $\sym(\Omega)_{\{\F\}}$ contains the pointwise 
stabilisers in $\sym(\Omega)$ of the complements of $\Omega f$ and  $\Lambda(g)$, it follows from Lemmas \ref{sym_inj} 
and \ref{sym_surj} that there exist $f^*, g^*\in M$ with $f^*$ injective, $g^*$ surjective, $d(f^*)=c(g^*)=|\Omega|$,  
$\Omega f^*\subseteq \Omega f$ and a transversal $\Lambda(g^*)\subseteq \Lambda(g)$ for $g^*$. Also 
since $\Omega f, \Lambda(g)\not\in \F$ it follows that $\Omega f^*, \Lambda(g^*)\not \in \F$.  
Since $\sym(\Omega)_{\{\F\}}$  is contained in $M$ and it is transitive on moieties not belonging to $\F$, it follows that 
every element of $\sym(\Omega)$ can be given in the form $f^*a g^*$ for some $a\in \sym(\Omega)_{\{\F\}}$. In particular, 
$\sym(\Omega)\subseteq M$, and so, by Theorem \ref{dicks}, $M=\Omega^\Omega$. 
\qed\vspace{\baselineskip}

%%%%%%%%%%%%%%%%%%%%%%%%%
%%%%%%%%%%%%%%%%%%%%%%%%%

\section{The almost stabiliser of a finite partition - the proof of Theorem  \ref{almost_stab}}\label{proof_almost_stab}

Recall that a finite partition of $\Omega$ is a partition of $\Omega$ into finitely many moieties. 
Throughout this section we denote the finite partition $\{\Sigma_0, \Sigma_1, \ldots, \Sigma_{n-1}\}$ of $\Omega$  with
 $n\geq 2$ by $\P$, and we write
$$\stab(\mathcal{P})=\set{g\in\sym(\Omega)}{(\forall\ i)(\exists\ j)\Sigma_ig= \Sigma_j}$$
for the stabiliser of $\mathcal P$.

A \emph{binary relation} on an arbitrary set $\Lambda$ is just a subset of $\Lambda \times \Lambda$. 
If $\rho$ and $\sigma$ are binary relations on $\Lambda$, then the composition $\rho\sigma$ of $\rho$ and $\sigma$ is 
defined to be $$\rho\sigma=\set{(\alpha,\beta)\in\Lambda\times\Lambda}{(\exists \gamma)(\alpha,\gamma)\in\rho\text{ and }(\gamma,\beta)\in \sigma}.$$
Composition of binary relations is associative and so we may refer to the semigroup generated by a set of binary relations.  A relation $\rho$ on $\Lambda$ is \emph{total} if $\alpha\rho=\set{\beta\in\Lambda}{(\alpha,\beta)\in\rho}\not=\emptyset$ for all $\alpha\in\Lambda$.

%%%%%%%%%%%%%%%%%%%%%%%%%

Recall that if $f\in\Omega^\Omega$, then  $\rho_f$ is the binary relation on $n=\{0,1,\ldots, n-1\}$ defined in \eqref{rho_f} as 
$$\rho_f=\set{(i,j)}{|\Sigma_if\cap \Sigma_j|=|\Omega|}.$$

The purpose of this section is to prove the following theorem.

%%%%%%%%%%%%%%%%%%%%%%%%%

\begin{main}
Let $\Omega$ be any infinite set and let $\mathcal{P}=\{\Sigma_0, \Sigma_1, \ldots,\Sigma_{n-1}\}$, $n\geq 2$, be a finite 
partition of $\Omega$. Then the maximal subsemigroups of $\Omega^\Omega$ containing $\stab(\P)$ but not 
$\sym(\Omega)$ are:
\begin{align*}
A_1(\P)&=\set{f\in \Omega^\Omega}{\rho_f\in \sym(n)\text{ or }\rho_f\text{ is not total }};\\
A_2(\P)&=\set{f\in \Omega^\Omega}{\rho_f\in \sym(n)\text{ or }\rho_f^{-1}\text{ is not total }}.
\end{align*}
\end{main}

%%%%%%%%%%%%%%%%%%%%%%%%%
We start by showing that $A_1(\P)$ and $A_2(\P)$ in Theorem \ref{almost_stab} are semigroups.
\begin{prop}
The sets $A_1(\P)$ and $A_2(\P)$ as defined in Theorem \ref{almost_stab} are subsemigroups of $\Omega^\Omega$ 
and neither is a subset of the other nor of any of the semigroups in Theorem \ref{dicks}. 
\end{prop}
\proof 
It is easy to verify that neither $A_1$ nor $A_2$ is contained in the other, nor in any of the semigroups listed in Theorem 
\ref{dicks}. 
We only prove that $A_1(\P)$ is a subsemigroup of $\Omega^\Omega$; the proof that $A_2(\P)$ is a subsemigroup 
follows by a dual argument. 

Let $f, g\in A_1(\P)$. Then, certainly, $\rho_{fg}\subseteq \rho_f\rho_g$. Hence, if $\rho_f$ is not total, then 
$\rho_{f}\rho_{g}$ is not total, and so $\rho_{fg}$ is not either, whence $fg\in A_1(\P)$. 
Assume that $\rho_f\in \sym(n)$. Then either $\rho_f\rho_g\in \sym(n)$  or $\rho_f\rho_g$ is not total, depending on 
whether $\rho_g\in\sym(n)$ or $\rho_g$ is not total. Hence $\rho_{fg}\in \sym(n)$ or $\rho_{fg}$ is not total and in either 
case $fg\in A_1(\P)$. 
 \qed\vspace{\baselineskip}

We prove Theorem \ref{almost_stab} in a sequence of lemmas. If $\Sigma\subseteq \Omega$, then we denote by $\sym(\Sigma)$ the pointwise stabiliser of $\Omega\setminus \Sigma$ 
in $\sym(\Omega)$.

%%%%%%%%%%%%%%%%%%%%%%%%%

\begin{lem}\label{rho}
Let $f, g\in\Omega^\Omega$. Then there exists $a\in\stab(\P)$ such that $\rho_{fag}=\rho_{f}\rho_{g}$. 
\end{lem}
\proof 
Let $i\in\n$ be arbitrary.  If $j\in i\rho_f^{-1}$, then $|\Sigma_jf\cap \Sigma_i|=|\Omega|$ and so $\Sigma_jf\cap \Sigma_i$ can be partitioned into $|i\rho_g|+1$ moieties. If $k\in i\rho_g$, then $g$ has a \hyperlink{transversal}{transversal} that intersects $\Sigma_kg^{-1}\cap \Sigma_i$ in a set $\Gamma_k$ where $|\Gamma_k|=|\Omega|$. Hence $\Gamma_k$ can be partitioned into $|i\rho_f^{-1}|+1$ moieties.  Let $a_i\in \sym(\Sigma_i)$ be any element  mapping one of the moieties partitioning $\Sigma_jf\cap \Sigma_i$ to one of the moieties partitioning $\Gamma_k$ for all $j\in i\rho_f^{-1}$ and for all $k\in i\rho_g$. 
The required $a\in \stab(\P)$ is then just $a_0\cdots a_{n-1}$. 
\qed

%%%%%%%%%%%%%%%%%%%%%%%%%

\begin{lem}\label{b_fin}
Let $\rho$ and $\sigma$ be  (not necessarily distinct) binary relations on $\{0,1,\ldots, n-1\}$ such that $\rho$ and $\sigma^{-1}$ are total but $\rho, \sigma \not\in\sym(n)$. Then the semigroup $\genset{\sym(n), \rho, \sigma}$  contains the total relation $n\times n$. 
\end{lem}
\proof 
We prove that there exists $\tau_0\in \genset{\sym(n), \rho, \sigma}$ such that $0\tau_0=\n$. If this is the case, then by replacing $\rho$ by $\sigma^{-1}$ and $\sigma$ by $\rho^{-1}$, there exists $\tau_1\in \genset{\sym(n), \sigma^{-1}, \rho^{-1},}$ such that $0\tau_1=\n$.  Hence $\tau_1^{-1}\in \genset{\sym(n), \rho, \sigma}$ and $\tau_1^{-1}\tau_0=n\times n$, as required.

We may assume without loss of generality that $0\rho=\set{i}{(0,i)\in \rho}\not=\n$.  
Let $A$ be a subset of $\n$ with least cardinality such that $$A\sigma =\set{j}{(\exists i\in A)(i,j)\in \sigma}=\n.$$ 
Since $\sigma\not\in \sym(n)$, it follows that $|A|<n$ and without loss of generality that  $0\in A$ and $|0\sigma|>1$.
Also by the minimality of $A$,  for all $i\in A$ there exists $j\in i\sigma$ such that $j\not \in (A\setminus \{i\})\sigma$. 

If $|0\rho|\geq |A|$, then let  $a_0\in \sym(n)$ be any permutation such that $A\subseteq 0\rho a_0$. In this case, $0\rho a_0\sigma=\n$, as required.   If $|0\rho|<|A|$, then let $a_0\in \sym(n)$ be any permutation such that $0\in 0\rho a_0$ and $0\rho a_0\subsetneq A$. In this case, $|0\rho a_0\sigma|\geq |0\rho|+1>|0\rho|$.  
By repeating this argument we find $a_1, a_2, \ldots, a_m\in \sym(n)$ such that $0\rho a_0\sigma a_1\sigma \cdots a_m\sigma=\n$, as required.  
\qed

%%%%%%%%%%%%%%%%%%%%%%%%%

\begin{lem}\label{almost_stab_inj}
Let $f\in\Omega^\Omega$ be injective such that $d(f)>0$. Then there exists an injective $f^*\in \genset{\stab(\P), f}$ such that $|\Sigma_i\setminus \Omega f^*|\geq d(f)$ for all $i$ with $0\leq i\leq n-1$.  If $d(f)$ is infinite, then $|\Sigma_i\setminus\Omega f^*|=d(f)$ for all $i$.
\end{lem}
\proof
Let $\mu=d(f)$ and let $g=f^{2n}$.  By Lemma \ref{tech}(iii) and (iv), $g$ is injective and $d(g)=2n\mu$.  In particular, 
there exists $0\leq i\leq n-1$ such that $|\Sigma_i\setminus\Omega g|\geq2\mu$.  If 
$|\Sigma_j\setminus\Omega g|\geq\mu$ for all $0\leq j\leq n-1$, then the proof is completed by setting $f^*=g$. Suppose 
that there exists $j$ such that $0\leq j\leq n-1$ and $|\Sigma_j\setminus\Omega g|<\mu$.  It follows that 
$j\rho_g^{-1}\not=\emptyset$ and  so there exists $a\in\stab(\P)$ such that $i\rho_a\subseteq j\rho_g^{-1}$ and 
$|(\Sigma_i\setminus\Omega g)a\cap\Sigma_jg^{-1}|\geq\mu$.  Hence 
$$(\Sigma_i\setminus\Omega g)ag\cap\Sigma_j\subseteq(\Omega\setminus\Omega g)ag\cap\Sigma_j\subseteq(\Omega\setminus\Omega gag)\cap\Sigma_j=\Sigma_j\setminus\Omega gag$$ and so
$$|\Sigma_j\setminus\Omega gag|\geq|(\Sigma_i\setminus\Omega g)ag\cap\Sigma_j| \geq |(\Sigma_i\setminus\Omega g)ag\cap\Sigma_j g^{-1}g| = |(\Sigma_i\setminus\Omega g)a\cap\Sigma_jg^{-1}|\geq\mu.$$  Also, for all $0\leq k\leq n-1$ such that $|\Sigma_k\setminus\Omega g|\geq\mu$, we have $|\Sigma_k\setminus\Omega gag|\geq\mu$.  Thus, by repeating this process at most $n-1$ times, we obtain the required $f^*$.

If $d(f)$ is infinite, and $h\in\genset{\stab(\P),f}$, then either $d(h)=0$ or $d(h)=d(f)$ by Lemma \ref{tech}(ii) and (iii).  The final statement follows immediately.
\qed

\newcommand{\grayshade}[5]{
\draw[fill=red!15!white] (#1,#2) to [out=0, in=180] (#1+3,#3) to [out=90, in=270] (#1+3,#4) to [out=180, in=0] (#1,#5) to [out=270, in=90] (#1,#2);
}

\begin{figure}
\begin{center}
 \begin{tikzpicture}[scale=1]
 \verticalline0
 \verticalline3
 \verticalline6
 \verticalline9

 \chip34
 \chip64

 \verticallinearrowleft002{$\Sigma_i$}
% \verticallinearrowright3{0.5}{2}{$\geq\!\mu$}
 \verticallinearrowleftlow3{0}{1.45}{$\Sigma_i\setminus\Omega g$}
 \verticallinearrowleft6{2.75}6{$\Sigma_jg^{-1}$}
 \verticallinearrowrighthigh{9}46{$\Sigma_j$}

% \verticallinearrowright9{0.5}2{$\geq\!\mu$}
  \verticallinearrowright6{2.8}{3.5}{$\geq\!\mu$}
    \verticallinearrowright{9.2}4{5}{$\geq\!\mu$}
    \verticallinearrowleft3{1.5}{5.5}{$\Omega g$}

 \functionline00{1.5}
 \functionline06{5.5}
 \functionlinedash366
 \functionline3{5.5}{5.5}
 \functionline344
 \functionlinearightdash302
 \functionlinearightdash324
 \functionlinearight342
 \functionlinearight320
 \functionlinealeftdash302

\functionlinealeftdash3{1.5}{3.5}
\functionlinearightdash3{1.5}{3.5}
 \functionlinealeftdash324
 \functionlinealeft342
 \functionlinealeft320
 \functionline623
 %\functionline6{2.75}{3.75}
 \functionlinedash6{2.75}4
 \functionline60{1.5}
% \functionline623
 %\functionline62{3.25}
 \functionline6{3.5}5
 \functionline6{5.5}{5.25}
 \functionlinedash66{5.5}
% \functionline6{2.5}{3.5}
%
 \functionlabel0{g}
 \functionlabel3{a}
 \functionlabel6{g}
  \colorinterval935{blue}
 \colorinterval9{1.5}3{red}
  \colorinterval95{5.25}{red}
 \colorinterval30{1.5}{blue}
  \colorinterval62{3.5}{blue}
 \colorinterval006{red}
  \colorinterval602{red}
 \colorinterval3{1.5}{5.5}{red}
 \colorinterval6{3.5}{5.5}{red}
 \colorinterval344{red}
 \colorinterval644{red}
 \colorinterval322{red}
  \chip92
 \chip94
  \chip02
 \chip04
 \smallchip6{2.75}
 \smallchip9{5.5}
 
 \end{tikzpicture}
\caption{The composite $gag$ in the proof of Lemma \ref{almost_stab_inj}.}
\end{center}
\end{figure}

%%%%%%%%%%%%%%%%%%%%%%%%%

\begin{lem} \label{almost_stab_surj}
Let $g\in \Omega^\Omega$ be surjective such that  $c(g)>0$. Then there exists  $g^*\in\genset{\stab(\P), g}$ and a
 transversal $\Gamma$ of $g^*$ such that $|\Sigma_i\setminus \Gamma|\geq c(g)$ for all $i$ such that $0\leq i\leq  n-1$.  If 
 $c(g)$ is infinite, then $|\Sigma_i\setminus\Gamma|=c(g)$ for all $i$.
\end{lem}
\proof
If $f$ is any inverse of $g$, then $f$ is injective and $d(f)=c(g)>0$, and so by Lemma 
\ref{almost_stab_inj} there exists $f^*\in \genset{\stab(\P), f}$  such that $|\Sigma_i\setminus \Omega f^*|\geq d(f)$ for all
 $i$ with $0\leq i\leq n-1$. But every element of $\genset{\stab(\P), f}$ is injective, and so $\Omega$ is the unique 
 transversal of every element in $\genset{\stab(\P), f}$. In particular, if $\Lambda$ is the unique assignment of transversals 
 for  $\genset{\stab(\P), f}$, then $\genset{\stab(\P), f}=\mutt{\genset{\stab(\P), f}}{\Lambda}$ and so 
 $f^*\in\mutt{\genset{\stab(\P), f}}{\Lambda}$.  Thus, by Corollary \ref{dual_cor}, $f^*$ has an inverse $g^*$ in 
 $\genset{\stab(\P), g}$.  Moreover, if $\Gamma=\Omega f^*$, then $\Gamma$ is a transversal of $g^*$ and 
 $$|\Sigma_i\setminus \Gamma|=|\Sigma_i\setminus \Omega f^*|\geq d(f)=c(g),$$
 for all $i$.  
 \qed 

%%%%%%%%%%%%%%%%%%%%%%%%%

\begin{lem} \label{extract}
Let $U$ be a subsemigroup of $\Omega^\Omega$ containing $\stab(\P)$ such that there exist $f,g,t\in U$ and the following hold:
\begin{itemize}
\item[\rm (i)] $f$ is injective, $g$ is surjective, and $d(f)=c(g)=|\Omega|$;
\item[\rm (ii)] $\rho_t=n\times n$.
\end{itemize}
Then $\sym(\Omega)$ is contained in $U$. 
\end{lem}
\proof 
We start by showing that there are $f_0, g_0\in U$ such that $f_0$ is injective, $\Omega f_0$ is a moiety of  $\Sigma_0$, and $\Omega f_0g_0=\Omega$.   
 By Lemma \ref{almost_stab_inj} there exists $f^*\in \genset{\stab(\P), f}$ such that  $|\Sigma_i\setminus \Omega f^*|=|\Omega|$ for all $i$ such that $0\leq i\leq n-1$.  Since $f$ is injective, every element of $\genset{\stab(\P), f}$ is injective, and so, in particular, $f^*$ is injective.
Let $0\leq i\leq n-1$ be arbitrary.  By assumption, $\Omega f^*\cap \Sigma_i$ is contained in a moiety of $\Sigma_i$. Also since $\rho_t=n\times n$, it follows that  $\Sigma_0t^{-1}\cap \Sigma_i$ is  a moiety of $\Sigma_i$. In particular, there exists a transversal $\Gamma_i$ of $t|_{\Sigma_0t^{-1}\cap \Sigma_i}$ such that $\Gamma_i$ is a moiety of $\Sigma_i$.  Hence there exists 
$a_0\in \stab(\P)$ such that $(\Omega f^*\cap \Sigma_i)a_0\subseteq \Gamma_i$ for all $i$. Then $\Omega f^*a_0t\subseteq \Sigma_0$ 
and so $\Omega (f^*a_0t)^2$ is a moiety of $\Sigma_0$. 
Thus $f_0=(f^*a_0t)^2$ is the required mapping.

For each $i$, let $\Delta_i$ be a transversal of $t|_{\Sigma_it^{-1}\cap\Sigma_0}:\Sigma_0\to \Sigma_i$.  
So each $\Delta_i$ is a moiety of 
$\Sigma_0$. Let $a_1\in\stab(\P)$ be any permutation such that 
 $\Omega f_0a_1\cap \Delta_i$ is a moiety of $\Delta_i$ for all $i$. 
Then $|\Omega f_0 a_1 t\cap \Sigma_i|=|\Omega|$ for all $i$. By Lemma \ref{almost_stab_surj}, there exists 
$g^*\in\genset{\stab(\P), g}$  and a transversal $\Lambda$ of $g^*$ such that  $|\Sigma_i\setminus \Lambda|=|\Omega|$ for 
all $i$ such that $0\leq i\leq  n-1$.   In other words, $\Lambda\cap \Sigma_i$ is contained in a moiety of $\Sigma_i$ for all 
$i$. Since $g$ is surjective, every element of $\genset{\stab(\P), g}$ is surjective, and so $g^*$ is surjective. Therefore there 
exists $a_2\in \stab(\P)$ such that $\Omega f_0 a_1t a_2$ contains $\Lambda$. Hence 
$\Omega f_0 a_1ta_2g^*=\Omega g^*=\Omega$ and $g_0=a_1ta_2g^*$ is the required function.

To conclude, let $b\in \sym(\Omega)$ be arbitrary. Then if $\Gamma$ is a transversal of $g_0$ contained in $\Omega f_0$, there exists  $a_3\in \stab(\P)$ such that $\alpha f_0a_3\in \alpha bg_0^{-1}\cap \Gamma$ for all $\alpha\in \Omega$.  But then $b=f_0a_3g_0$, and so $\sym(\Omega)$ is contained in  $U$, as required. \qed\vspace{\baselineskip}

%%%%%%%%%%%%%%%%%%%%%%%%%

At this stage it is straightforward to classify the maximal subsemigroups of $\Omega^\Omega$ containing the almost 
stabiliser of a finite partition using the results proved so far.  Since the stabiliser is a subgroup of the almost stabiliser 
this classification is actually a corollary of Theorem \ref{almost_stab}. To prove the more general Theorem 
\ref{almost_stab} we require two further lemmas which are similar to Lemmas \ref{sym_inj} and \ref{sym_surj}. 

%%%%%%%%%%%%%%%%%%%%%%%%%

\begin{cor}\label{almost_stab_cor2}
Let $\Omega$ be any infinite set and let $\mathcal{P}=\{\Sigma_0, \Sigma_1, \ldots,\Sigma_{n-1}\}$, $n\geq 2$, be a finite partition of $\Omega$. Then the maximal subsemigroups of $\Omega^\Omega$ containing $\astab(\P)$ but not $\sym(\Omega)$ are:
\begin{align*}
A_1(\P)&=\set{f\in \Omega^\Omega}{\rho_f\in \sym(n)\text{ or }\rho_f\text{ is not total }};\\
A_2(\P)&=\set{f\in \Omega^\Omega}{\rho_f\in \sym(n)\text{ or }\rho_f^{-1}\text{ is not total }}.
\end{align*}
\end{cor}
\proof 
Let $M$ be a subsemigroup of $\Omega^\Omega$ containing $\astab(\P)$ but which is not contained in any of the 
semigroups listed in Theorems \ref{dicks} or \ref{almost_stab}.  As in the proof of Theorem \ref{stab_finite}, it suffices to 
show that $M=\Omega^\Omega$.

Since $M\not\subseteq A_1(\P), A_2(\P)$, there exist $f,g\in M$ such that $\rho_f$ and $\rho_g^{-1}$ are total but 
$\rho_f, \rho_g\not\in \sym(n)$. 
Hence, by Lemmas \ref{rho} and \ref{b_fin}, there exists $t\in M$ such that $\rho_t=n\times n$.  Since $\astab(\P)$ contains $\set{a\in \sym(\Omega)}{|\supp(a)|<|\Omega|}$, it follows by Lemmas 
\ref{sym_inj} and~\ref{sym_surj} that there exist $f^*,g^*\in M$ such that $f^*$ is injective, $g^*$ is surjective, and $d(f^*)=c(g^*)=|\Omega|$. 
Thus, by Lemma \ref{extract}, $\sym(\Omega)$ is contained in $M$, and so, by Theorem \ref{dicks}, $M=\Omega^\Omega$, as required.
\qed\vspace{\baselineskip}

We return to the proof of Theorem \ref{almost_stab}. 

%%%%%%%%%%%%%%%%%%%%%%%%%

\begin{lem}\label{inj_almost_stab}
Let $U$ be a subset of $\Omega^\Omega$, which contains $\stab(\P)$ but which is not contained in  $S_2$, or 
$S_4(\mu)$ for any infinite $\mu\leq |\Omega|$ and let $\Lambda$ be any assignment of transversals for $U$ (as defined 
in Definition \ref{assign}). Then there exists an injective 
$f\in \mutt{U}{\Lambda}$ such that $d(f)=|\Omega|$.
\end{lem}
\proof 
We prove by transfinite induction that for all cardinals $\mu \leq | \Omega|$ 
\begin{equation}\label{hypo_inj_astab}
\text{ there exists an injective }f_{\mu}\in \mutt{U}{\Lambda}\text{ with }d(f_{\mu})\geq \mu. 
\end{equation}
Since $U$ is not contained in $S_2$, there exists an injective $h_0\in U\subseteq \mutt{U}{\Lambda}$ such that 
$d(h_0)>0$. By taking powers of $h_0$ (which also belong to $\mutt{U}{\Lambda}$) and applying Lemma \ref{tech}(iii) 
and (iv), it follows that \eqref{hypo_inj_astab} holds for all finite $\mu$. 

Let $\mu$ be any cardinal such that $\aleph_0\leq\mu\leq|\Omega|$ and assume that \eqref{hypo_inj_astab} holds for 
every cardinal $\nu<\mu$. Since $U\not\subseteq S_4(\mu)$, there exists $h_1\in U$ such that 
$c(h_1)<\mu\leq d(h_1)$.  
By our inductive hypothesis, there exists an injective $f_{c(h_1)}\in \mutt{U}{\Lambda}$ such that 
$d(f_{c(h_1)})\geq c(h_1)$.  

\begin{figure}
\begin{center}
 \begin{tikzpicture}[scale=1]
 \verticalline0
 \verticalline3
 \verticalline6
 \verticalline9

 \smallchip31
 %\smallchip32
 \smallchip33
 \chip34
 \smallchip35
 \smallchip61
 \chip62
 \smallchip63
 \smallchip64
 \smallchip65
 \smallchip6{5.5}
 \smallchip6{2.5}
 \smallchip6{1.5}
 \smallchip9{2.75}
 \smallchip9{2.5}
 \smallchip9{.25}
 \smallchip9{5}
 \verticallinearrowleft046{$\Sigma_0$}
 \verticallinearrowleft024{$\Sigma_1$}
 \verticallinearrowleft002{$\Sigma_2$}
 \verticallinearrowright{2.9}45{$\geq\!c(h_1)$}
 \verticallinearrowright{2.9}34{$\geq\!c(h_1)$}
 \verticallinearrowright{2.9}01{$\geq\!c(h_1)$}
 \verticallinearrowright6{5.5}6{$\leq\!c(h_1)$}
 \verticallinearrowright62{2.5}{$\leq\!c(h_1)$}
 \verticallinearrowright6{1.5}2{$\leq\!c(h_1)$}
 \verticallinearrowright6{2.5}{5.5}{$\Lambda(h_1)$}
 \verticallinearrowright60{1.5}{$\Lambda(h_1)$}
 \verticallinearrowright956{$\geq\!\mu$}
 \functionline001
 \functionline066
 \functionline0{3.5}5
 \functionline0{3.5}3
 \functionlinealeft365
 \functionlinealeftdash356
 \functionlinealeft354
 \functionlinealeftdash345
 \functionlinealeftdash343
 \functionlinealeft334
 \functionlinealeftdash332
 \functionlinealeft323
 \functionlinealeft321
 \functionlinealeftdash312
 \functionlinealeft310
 \functionlinealeftdash301
 \functionlinearight365
 \functionlinearightdash356
 \functionlinearight354
 \functionlinearightdash345
 \functionlinearightdash343
 \functionlinearight334
 \functionlinearightdash332
 \functionlinearight323
 \functionlinearight321
 \functionlinearightdash312
 \functionlinearight310
 \functionlinearightdash301
 \functionlinedash6{5.5}5
 \functionlinedash6{2.5}{2.75}
 \functionlinedash6{1.5}{2.5}
 \functionline60{.25}
 \functionline65{4.5}
 \functionline63{3.5}
 \functionline61{1.5}
 \functionlabel0{f^*}
 \functionlabel3{a}
 \functionlabel6{h_1}
 \colorinterval006{red}
 \colorinterval356{red}
 \colorinterval313{red}
 \colorinterval645{red}
 \colorinterval601{red}
  \colorinterval634{red}
  \colorinterval9{0.25}{1.5}{red}
    \colorinterval9{3.5}{4.5}{red}
\colorinterval322{red}
\colorinterval644{red}

  \chip92
  \chip94
  %\colorinterval323{red}
 %\colorinterval324{blue}
 %\colorinterval646{red}
% \colorinterval602{green}t
% \colorinterval624{blue}
% \colorinterval946{red}
 %\colorinterval902{green}
 %\colorinterval924{blue}
 %
 %
 \chip02
 \chip04

 \end{tikzpicture}
\caption{The composite $f^*ah_1$ in the proof of Lemma \ref{inj_almost_stab}.}\label{fig10}
\end{center}
\end{figure}
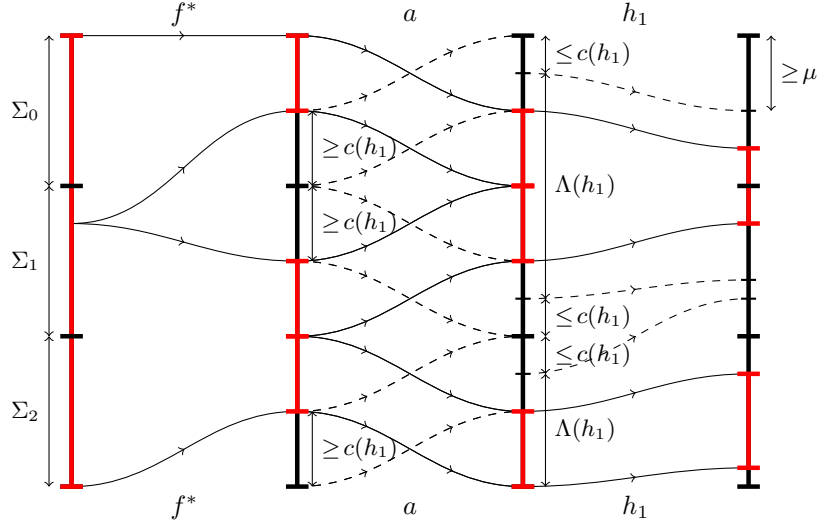

By Lemma \ref{almost_stab_inj}, there exists an injective $f^*\in\genset{\stab(\P),f_{c(h_1)}}$ such that 
$|\Sigma_i\setminus\Omega f^*|\geq d(f_{c(h_1})$ for all $0\leq i\leq n-1$.  Since every element of $\genset{\stab(\P),f_{c(h_1)}}$ 
is injective, it follows that $\genset{\stab(\P),f_{c(h_1)}}\subseteq \mutt{U}{\Lambda}$ and so $f^*\in \mutt{U}{\Lambda}$. 
Then, since 
$$|\Sigma_i\setminus\Lambda(h_1)|\leq|\Omega\setminus\Lambda(h_1)|=c(h_1)\leq d(f_{c(h_1)})\leq|\Sigma_i\setminus\Omega f^*|$$
for all $0\leq i\leq n-1$, there exists $a\in\stab(\P)$ such that $(\Omega f^*\cap\Sigma_i)a\subseteq\Lambda(h_1)\cap\Sigma_i$ 
for all $i$; see Figure \ref{fig10}.  Hence, if we set $f_\mu=f^*ah_1$, then, since $f_{\mu}, a\in \mutt{U}{\Lambda}$ and by the definition of $a$, 
it follows that $f_{\mu}\in \mutt{U}{\Lambda}$, $f_\mu$ is injective, and $d(f_\mu)\geq d(h_1)\geq\mu$ by Lemma 
\ref{tech}(ii), as required. 
\qed

%%%%%%%%%%%%%%%%%%%%%%%%%

\begin{lem}\label{surj_almost_stab}
Let $U$ be a subset of $\Omega^\Omega$, which contains $\stab(\P)$ but which is not contained in  $S_1$, or 
$S_3(\mu)$ for any infinite $\mu\leq |\Omega|$. Then there exists a surjective $g\in U$ such that $c(g)=|\Omega|$.
\end{lem}
\proof 
Let $u'\in \Omega^\Omega$ be an arbitrary inverse for $u$ for all  
$u\in U$. We denote $\set{u'\in \Omega^\Omega}{u\in U}$ by $U'$ and 
we set $\Lambda: U'\to \P(\Omega)$ to be the assignment 
of transversals for $U'$ defined by $\Lambda(u')=\Omega u$.  Recall that $c(u)=d(u')$  and  $d(u)=c(u')$ 
for all $u\in U$. 

We prove that $U'$ satisfies the conditions of Lemma \ref{almost_stab_inj}. Since 
$U\not\subseteq S_1$, $U\not\subseteq S_3(\mu)$, it follows that $U'\not\subseteq S_2$ and 
$U'\not\subseteq S_4(\mu)$ for all 
infinite  $\mu\leq |\Omega|$. 
Thus by Lemma \ref{almost_stab_inj} there exists an injective $f^*\in \mutt{U'}{\Lambda}$ such 
that $d(f^*)=|\Omega|$.  By Corollary \ref{dual_cor}, $\genset{U}$ 
contains an inverse $g^*$ of $f^*$.  Therefore $c(g^*)=d(f^*)=|\Omega|$. \qed\vspace{\baselineskip}

%%%%%%%%%%%%%%%%%%%%%%%%%

\proofref{almost_stab}
 Let $M$ be a subsemigroup of $\Omega^\Omega$ containing $\stab(\P)$ but not contained in any of the semigroups listed in Theorems \ref{dicks} or \ref{almost_stab}. 
Since $M\not\subseteq A_1(\P), A_2(\P)$, there exists $f,g\in M$ such that $\rho_f$ and $\rho_g^{-1}$ are total but $\rho_f, \rho_g\not\in \sym(n)$. 
Hence, by Lemmas \ref{rho} and \ref{b_fin}, there exists $t\in M$ such that $\rho_t=n\times n$.  Also by Lemmas \ref{inj_almost_stab} and~\ref{surj_almost_stab} there exist $f^*,g^*\in M$ such that $f^*$ is injective, $g^*$ is surjective, and $d(f^*)=c(g^*)=|\Omega|$. 
Thus, by Lemma \ref{extract}, $\sym(\Omega)$ is contained in $M$, and so, by Theorem \ref{dicks}, $M=\Omega^\Omega$, as required.
 \qed

%%%%%%%%%%%%%%%%%%%%%%%%%
%%%%%%%%%%%%%%%%%%%%%%%%%

\section{Maximal subsemigroups of the symmetric group}\label{subsemigroups}

In this section we prove that the stabiliser of a non-empty finite set, the almost stabiliser of a finite partition, and the stabiliser of an 
ultrafilter  are maximal subsemigroups of the symmetric group and not just maximal subgroups.

Let $T$ be a subsemigroup of $\sym(\Omega)$, and let $G$ denote the group generated by $T$ . If 
$G\not=\sym(\Omega)$ and $T\not=G$ then, for any $f\in G\setminus T$, the semigroup generated by $T$ and $f$ is a 
subsemigroup of~$G$. In particular, $\genset{T,f}\not=\sym(\Omega)$  and hence $T$ is not maximal. (We remind the 
reader that~$\genset{U}$ always denotes the \emph{semigroup} generated by $U$.) Hence the group generated by any 
maximal subsemigroup of $\sym(\Omega)$ that is not a subgroup is $\sym(\Omega)$.

%%%%%%%%%%%%%%%%%%%%%%%%%

\begin{thm} \label{fin_semi}
Let $\Omega$ be any infinite set and let $\Sigma$ be a non-empty finite subset of $\Omega$. Then the setwise stabiliser $\sym(\Omega)_{\{\Sigma\}}$ of $\Sigma$ in $\sym(\Omega)$ is a maximal subsemigroup of $\sym(\Omega)$. 
\end{thm}

\proof
Let $f\in\sym(\Omega)\setminus \sym(\Omega)_{\{\Sigma\}}$.  We must show that 
$\genset{\sym(\Omega)_{\{\Sigma\}},f}=\sym(\Omega)$, i.e. that the semigroup generated by 
$\sym(\Omega)_{\{\Sigma\}}$ and $f$ is $\sym(\Omega)$. 
.  Since $\sym(\Omega)_{\{\Sigma\}}$ is a maximal subgroup of 
$\sym(\Omega)$, it suffices to find $g\in\genset{\sym(\Omega)_{\{\Sigma\}},f}$ such that $g$ has finite order and 
$g\not\in\sym(\Omega)_{\{\Sigma\}}$.  By postmultiplying by an element of $\sym(\Omega)_{\{\Sigma\}}$ if necessary, we 
may assume without loss of generality that every nontrivial cycle of $f$ contains an element of~$\Sigma$.  Since $\Sigma$ is 
finite, if every cycle of $f$ is finite, then $f$ itself has finite order, and setting $g=f$ concludes the proof in this case.  So 
suppose $f$ has at least one infinite cycle.  There exists $m\in\nat$ such that $f^m$ has only infinite cycles, each of which 
contains at most one element of $\Sigma$.  
Again we may assume without loss of generality that every nontrivial cycle of $f^m$ contains precisely one element of 
$\Sigma$. Then $f^m=c_1\cdots c_r$, where $c_1,\ldots,c_r$ are disjoint infinite cycles. We may write 
$c_i=(\ldots,\alpha_{i,-1},\alpha_{i,0},\alpha_{i,1},\alpha_{i,2},\ldots)$ where $\alpha_{i,0}\in\Sigma$. We let  
$d_i=(\ldots,\alpha_{i,2},\alpha_{i,1},\alpha_{i,-1},\alpha_{i,-2},\ldots)$ and $h=d_1\cdots d_r$. Then  $h\in\sym(\Omega)_{\{\Sigma\}}$ and 
 $$g=hf^m=(\alpha_{1,0},\alpha_{1,1})\cdots(\alpha_{r,0},\alpha_{r,1})\in\genset{\sym(\Omega)_{\{\Sigma\}},f}\setminus\sym(\Omega)_{\{\Sigma\}}$$ has order~$2$, which completes the proof.
\qed\vspace{\baselineskip}

%%%%%%%%%%%%%%%%%%%%%%%%%

If $H$ and $K$ are subgroups of  a group $G$, then the subsemigroup generated by $H$ and $K$ equals the group generated by $H$ and $K$.  Thus the following two lemmas are immediate consequences of the corresponding results about subgroups given in  \cite{Dixon1986aa} and \cite[Note 3(iii) of \S 4]{Neumann1988ab}, respectively. 

%%%%%%%%%%%%%%%%%%%%%%%%%

\begin{lem} \label{dixon}
If $\Gamma_1, \Gamma_2\subseteq \Omega$ and $|\Gamma_1\cap \Gamma_2|=\min\{|\Gamma_1|, |\Gamma_2|\}$, then 
$\sym(\Gamma_1\cup \Gamma_2)$  equals the subsemigroup $\genset{\sym(\Gamma_1), \sym(\Gamma_2)}$  generated by the subgroups $\sym(\Gamma_1)$ and $\sym(\Gamma_2)$. 
\end{lem}

%%%%%%%%%%%%%%%%%%%%%%%%%

%A subsemigroup $S$ of $\sym(\Omega)$ is \emph{full} on $\Sigma\subseteq \Omega$ if for all $f\in \sym(\Sigma)$ there exists $g\in S$ such that $g|_{\Sigma}=f|_{\Sigma}$. 

\begin{lem} \label{full1}
Let $S$ be a subsemigroup of $\sym(\Omega)$. 
If $S$ contains $\sym(\Sigma)$ for all moieties $\Sigma$ of $\Omega$, then $S=\sym(\Omega)$. 
\end{lem}

%%%%%%%%%%%%%%%%%%%%%%%%%

\begin{lem} \label{full2}
Let $S$ be a subsemigroup of $\sym(\Omega)$. 
If $S$ contains $\sym(\Sigma)$ for some  moiety $\Sigma$ of $\Omega$ and $S$ is transitive on moieties of 
$\Omega$, then $S=\sym(\Omega)$. 
\end{lem}
\proof 
It suffices by Lemma \ref{full1} to show that $S$ contains $\sym(\Gamma)$ for every moiety $\Gamma$ of 
$\Omega$. Let $\Gamma$ be any moiety of $\Omega$ and let $f\in \sym(\Gamma)$ be arbitrary. There exist 
$g, h, k\in S$ such that
$$\Gamma g=\Sigma,\quad \Sigma h=\Omega\setminus \Sigma,\quad (\Omega\setminus \Sigma)k
=\Gamma.$$
Since $(\Sigma)g^{-1}fk^{-1}h^{-1}=\Sigma$, it follows that there 
 exists $a\in \sym(\Sigma)\subseteq S$ such that $a|_{\Sigma}=g^{-1}fk^{-1}h^{-1}|_{\Sigma}$.
%It follows that $gah|_{\Gamma}=f(kl)^{-1}|_\Gamma$ and since $a$ fixes $\Omega\setminus \Sigma$ pointwise, it follows that 
%$gah|_{\Omega\setminus \Gamma}=gh|_{\Omega\setminus \Gamma}$. 
Also  $\Sigma h^{-1}g^{-1}k^{-1}=\Sigma$, there exists $b\in \sym(\Sigma)\subseteq S$
 such that $b|_{\Sigma}=h^{-1}g^{-1}k^{-1}|_{\Sigma}$. 
 
 We will show that $f=gahbk\in S$. 
 If $\alpha\in \Omega\setminus \Gamma$ is arbitrary, then $\alpha g\in \Omega\setminus \Sigma$ and so $\alpha ga=\alpha g$, and
 $\alpha gh\in \Sigma$ and so $\alpha ghb=\alpha k^{-1}$. Therefore
 $$(\alpha)gahbk=(\alpha)ghbk=(\alpha)k^{-1}k=\alpha.$$
If $\beta \in \Gamma$, then $\beta g\in \Sigma$ and so $\beta g a=\beta fk^{-1}h^{-1}$. Thus 
 $$(\beta)gahbk=(\beta) fk^{-1}bk$$
 and since $\beta fk^{-1}\in \Omega\setminus \Sigma$ and $b$ fixes $\Omega\setminus \Sigma$ pointwise, it follows that 
 $$(\beta)gahbk=(\beta)f,$$
 as required. 
\qed \vspace{\baselineskip}

%%%%%%%%%%%%%%%%%%%%%%%%%

Let $\Omega$ be an infinite set, let $\mathcal{P}=\{\Sigma_0, \Sigma_1, \ldots, \Sigma_{n-1}\}$, $n\geq 2$, be a finite 
partition of $\Omega$, and let $f\in \Omega^\Omega$. Recall that the binary relation $\rho_f$ on $\n$ is defined in 
Equation \eqref{rho_f} in Section \ref{subsection_almost} as:
\begin{equation*}
\rho_f=\set{(i,j)}{|\Sigma_if\cap \Sigma_j|=|\Omega|}.
\end{equation*} 

\begin{thm} \label{almost_stab_subsemi}
Let $\Omega$ be any infinite set and let $\P=\{\Sigma_0, \Sigma_1, \ldots, \Sigma_{n-1}\}$,  $n\geq 2$, be a finite partition of $\Omega$. Then  $\astab(\P)$ is a maximal subsemigroup of $\sym(\Omega)$. 
\end{thm} 
\proof 
Let $f\in \sym(\Omega)\setminus \astab(\P)$ be arbitrary. Then by Lemmas \ref{rho} and \ref{b_fin} 
there exists $g\in  \genset{\astab(\P), f}$ such that $\rho_g=n\times n$. 

Let $h\in \sym(\Sigma_0g^{-1})$. We will show that $h=gbga$ for some $a, b\in \astab(\P)$.  (In fact, $a,b$ will belong to $\stab(\P)$.)

Since $\rho_g=n\times n$, both $\Sigma_ig^{-1}\cap \Sigma_j$ and $\Sigma_ig \cap \Sigma_j$ are moieties in 
$\Sigma_j$ 
for all $i, j\in \{0,1,\ldots, n-1\}$. It follows that there exists $a\in \astab(\P)$ such that 
$$(\Sigma_ig)a=\Sigma_ig^{-1}$$
for all $i$. 
Define $b\in \sym(\Omega)$ by $\alpha b=\alpha g^{-1}ha^{-1}g^{-1}$ if $\alpha\in \Sigma_0$ and 
$\alpha b=\alpha g^{-1}a^{-1}g^{-1}$ if $\alpha\not\in \Sigma_0$.
Since $h\in \sym(\Sigma_0g^{-1})$, it follows that 
$(\Sigma_0g^{-1})h=\Sigma_0g^{-1}$ and so 
$$\Sigma_0b=\Sigma_0g^{-1}ha^{-1}g^{-1}=\Sigma_0g^{-1}a^{-1}g^{-1}=\Sigma_0gg^{-1}=\Sigma_0$$
and if $i\not=0$, then 
$$\Sigma_i b=\Sigma_ig^{-1}a^{-1}g^{-1}=\Sigma_igg^{-1}=\Sigma_i.$$
Hence $b\in \astab(\P)$.  Let $\alpha\in \Omega$ be arbitrary. If $\alpha\in\Sigma_0g^{-1}$, then $\alpha g\in \Sigma_0$ 
and so 
$$\alpha gbga=\alpha g g^{-1}ha^{-1}g^{-1}ga=\alpha h.$$
If $\alpha\not\in \Sigma_0g^{-1}$, then 
$$\alpha gbga=\alpha g g^{-1}a^{-1}g^{-1}ga=\alpha=\alpha h$$
 and so $h=gbga$, as required. It follows that  $\sym(\Sigma_0g^{-1})\leq \genset{\astab(\P), f}$. Therefore, since 
 $\Sigma_0g^{-1}\cap \Sigma_0$ and $\Sigma_0g^{-1}\cap \Sigma_1$ are moieties in $\Sigma_0$ and $\Sigma_1$, 
 respectively, and by Lemma \ref{dixon}, 
 $$\sym(\Sigma_0\cup \Sigma_1)\leq \genset{\sym(\Sigma_0), \sym(\Sigma_1), \sym(\Sigma_0g^{-1})}\leq \genset{\astab(\P), f}.$$ 
 Since $\astab(\P)$ is 2-transitive on $\Sigma_0, \ldots, \Sigma_{n-1}$, we conclude that 
 $\genset{\astab(\P), f}=\sym(\Omega)$ and so $\astab(\P)$ is a maximal subsemigroup of $\sym(\Omega)$. 
\qed\vspace{\baselineskip}

%%%%%%%%%%%%%%%%%%%%%%%%%

\begin{thm} \label{ultra_subsemi}
Let $\F$ be an ultrafilter on $\Omega$. Then the stabiliser $\sym(\Omega)_{\{\F\}}$ of $\F$ is a maximal subsemigroup of 
$\sym(\Omega)$. 
\end{thm}
\proof 
Let $f\in \sym(\Omega)\setminus \sym(\Omega)_{\{F\}}$.  Then either:
\begin{itemize}
	\item[(i)] there is a subset $\Sigma$ of $\Omega$ such that $\Sigma\in\F$ and $\Sigma f\not\in\F$, or
	\item[(ii)] there is a subset $\Gamma$ of $\Omega$ such that $\Gamma\not\in\F$ and $\Gamma f\in\F$.
\end{itemize}
It is straightforward to verify that $\Sigma$ and $\Gamma$ can be chosen to be moieties of $\Omega$. 
If (i) holds, then (ii) holds with $\Gamma=\Omega\setminus\Sigma$.  If (ii) holds, then (i) holds with 
$\Sigma=\Omega\setminus\Gamma$.  So we may assume that both (i) and (ii) hold.  Let $\Lambda$ and $\Delta$ be 
moieties of $\Omega$.  If $\Lambda$ and $\Delta$ both belong to $\F$ or neither belongs to $\F$, then there exists 
$a_0\in\sym(\Omega)_{\{\F\}}$ such that $\Lambda a_0=\Delta$.  If $\Lambda\in\F$ and $\Delta\not\in \F$, then we choose 
$a_1,a_2\in\sym(\Omega)_{\{\F\}}$ such that $\Lambda a_1=\Sigma$ and $(\Sigma f)a_2=\Delta$, and note that 
$\Lambda a_1fa_2=\Delta$.  Similarly, if $\Lambda\not\in\F$ and $\Delta\in\F$, then there exists 
$a_3,a_4\in\sym(\Omega)_{\{\F\}}$ such that $\Lambda a_3fa_4=\Delta$.  We have shown that
 $\genset{\sym(\Omega)_{\{\F\}},f}$ is transitive on moieties.  Since $\sym(\Omega)_{\{\F\}}$ is full on every moiety
  $\Xi\not\in\F$, the result follows from Lemma \ref{full2}. 
\qed

%%%%%%%%%%%%%%%%%%%%%%%%%
%%%%%%%%%%%%%%%%%%%%%%%%%

\bibliography{maximal}{}
\bibliographystyle{plain}

\end{document}